%% file: safereg_arXiv2017.tex
\documentclass[11pt]{article}
\input{declarations.tex}

\usepackage[round,comma]{natbib}
\usepackage[ruled,lined]{algorithm2e}
\usepackage{hyperref}
\DeclareMathOperator{\InvGamma}{Inv-gamma}
\newcommand{\betao}{\bar{\beta}_{0}} %

\newcommand{\etacrit}{\eta_{\text{crit}}}
\newcommand{\model}{{\cal M}}
\newcommand{\dens}{{f}}
\newcommand{\nc}{p}
\newcommand{\pmax}{\nc_{\max}}
\newcommand{\red}{\text{\sc red}}
\newcommand{\logrisk}{\text{\sc risk}^{\text{log}}}
\newcommand{\Rlogrisk}{\text{\sc risk}^{\text{R-log}}}
\newcommand{\squaredrisk}{\text{\sc risk}^{\text{sq}}}

\renewcommand{\pbayes}{\bar{\dens}}

\DeclareRobustCommand{\VAN}[3]{#2 #1}
\begin{document}
\title{Inconsistency of Bayesian Inference for Misspecified Linear Models, and a Proposal for Repairing It}
\date{}

\author{Peter Gr\"unwald \\ CWI, Amsterdam and Leiden University \\ The Netherlands \and 
Thijs van Ommen \\ University of Amsterdam \\ The Netherlands }

\maketitle

\begin{abstract}
  We empirically show that Bayesian inference can be inconsistent
  under misspecification in simple linear regression problems, both in
  a model averaging/selection and in a Bayesian ridge regression
  setting. We use the standard linear model, which assumes
  homoskedasticity, whereas the data are heteroskedastic, and observe
  that the posterior puts its mass on ever more high-dimensional
  models as the sample size increases.  To remedy the problem, we
  equip the likelihood in Bayes' theorem with an exponent called the
  learning rate, and we propose the {\em Safe Bayesian\/} method to
  learn the learning rate from the data. SafeBayes tends to select
  small learning rates as soon the standard posterior is not
  `cumulatively concentrated', and its results on our data are quite encouraging.   
\end{abstract}
\noindent
{\em This arxiv publication is the very first, 2014, version of the paper:
\begin{quote}{\rm P. D. Gr\"unwald, and T. van Ommen. Inconsistency of Bayesian Inference for Misspecified Linear Models, and a Proposal for Repairing It. {\em Bayesian Analysis}, 12(4), December 2017.  }\end{quote}
In 2017, the second version was posted which is exactly the same as the original version except for this first page, and an updated bibliography. \\ The BA (Bayesian Analysis) 2017 version is quite different though:}
due to length constraints,
\begin{enumerate}
\item
the BA  version only reports on a small subset of the experiments done here, and refers extensively to the additional experimental experiments done in the present paper\vspace*{-2 mm}. 
\item the BA version contains less details about the analytic calculations of the generalized posterior for regression\vspace*{-2 mm}.
\item the BA version contains no discussion of mix loss.
\end{enumerate}
  Otherwise, the paper has gone through many modifications and improvements. We refer to the BA version for:
\begin{enumerate}
\item A much more concise and better explanation of why small learning rates can vastly improve standard Bayes  under misspecification, and why SafeBayes `works'\vspace*{-2 mm}. 
\item A much more precise treatment of `hypercompression', the phenomenon underlying problems with misspecification\vspace*{-2 mm}.
\item A much updated discussion section\vspace*{-2 mm}.
\item A demonstration that standard theorems for consistency of Bayesian inference under misspecification do not apply to our standard regression model.
\end{enumerate}
\emph{Version 3 (2018) corrects the expression for $b_{n,\eta}$ in (\ref{eq:anbn}). We thank Tom Viering for spotting this error!}

\pagebreak
\section{Introduction}
\begin{figure}
\hspace*{0.1\textwidth}
\includegraphics[width=0.8\textwidth]{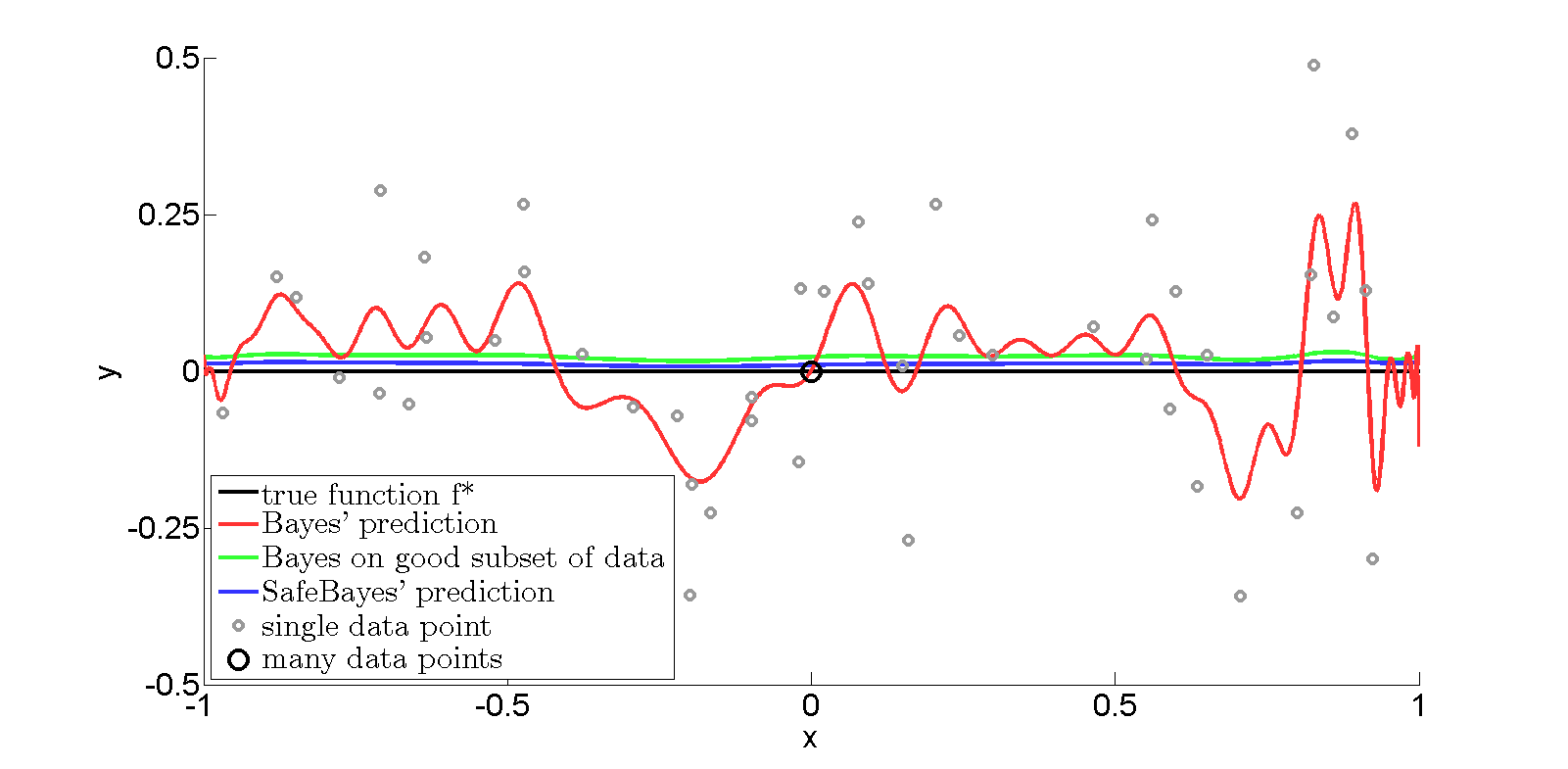}
\caption{\label{fig:polynomial} The conditional expectation
  $\Exp[Y|X]$ according to the full Bayesian posterior based on a
  prior on models $\cM_0, \ldots, \cM_{50}$ with polynomial basis
  functions, given 100 data points sampled i.i.d.~$\sim P^*$ (about 50
  of which are at $(0,0)$). Standard Bayes overfits, not as dramatically as
  maximum likelihood/unpenalized least squares, but still enough to
  show dismal predictive behavior as in Figure~\ref{fig:polynomialrisk}. In contrast,
  Safe Bayes (which chooses learning rate $\eta \approx 0.4$ here) and
  standard Bayes trained only at the points for which the model is
  correct (not $(0,0)$) both perform very well. }
\end{figure}

\begin{figure}
\hspace*{0.1\textwidth}
\includegraphics[width=0.8\textwidth]{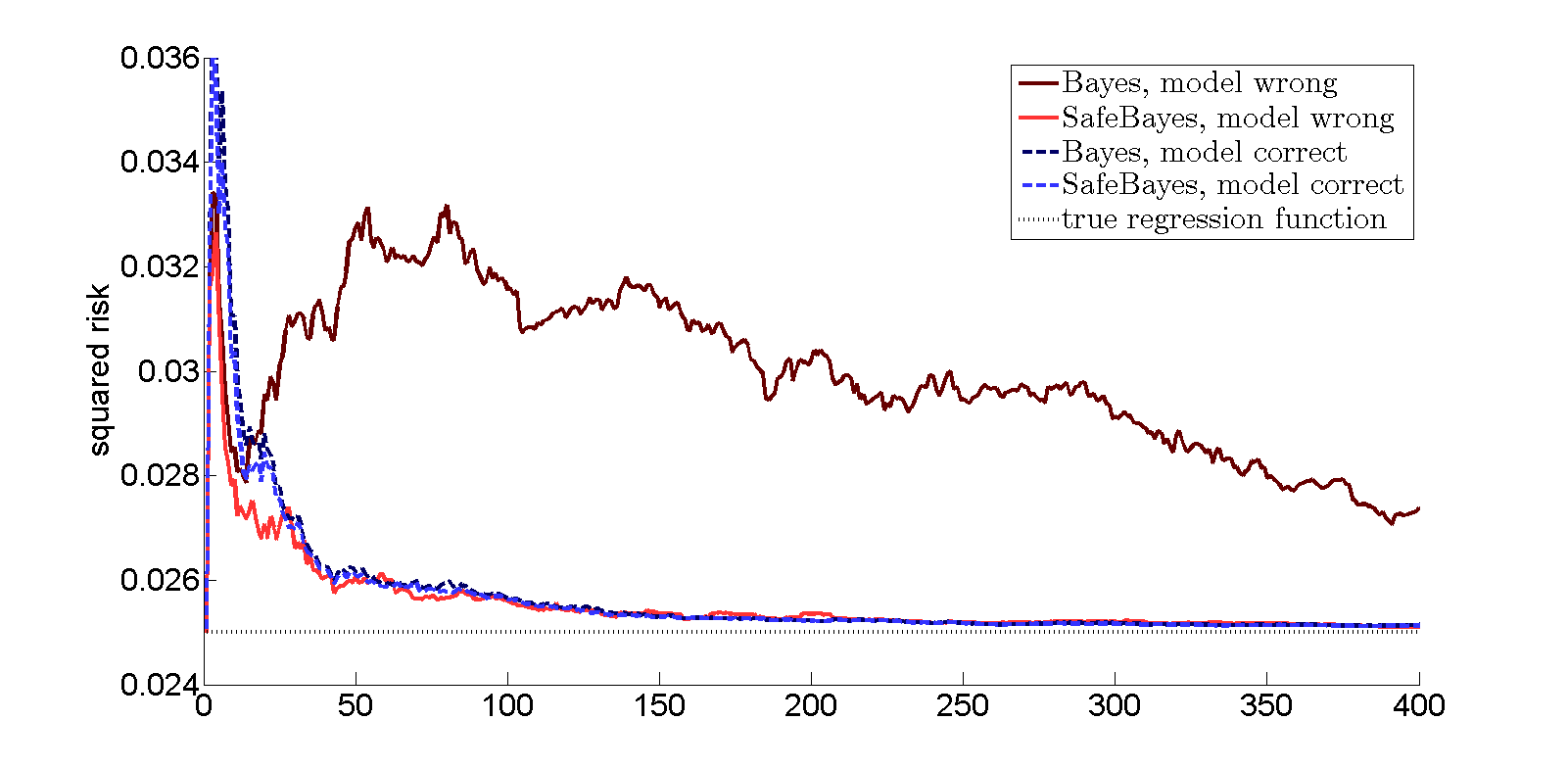}
\caption{\label{fig:polynomialrisk} The expected squared error risk obtained
  when predicting by the full Bayesian posterior (brown curve) and the
  Safe Bayesian posterior (red curve) and the optimal predictions
  (gray curve) as a function of
  sample size, for the setting of Figure~\ref{fig:polynomial}. SafeBayes is the $R$-log-version of SafeBayes defined in Section~\ref{sec:instantiationb}. Precise
definitions and further explanation in Section~\ref{sec:preparing} and Section~\ref{sec:statisticsshown}.}
\end{figure}

\paragraph{The Problem}
We empirically demonstrate the inconsistency of Bayes factor model
selection, model averaging and Bayesian ridge regression under model
misspecification on a simple linear regression problem with random
design. We sample data $(X_1,Y_1), (X_2, Y_2), \ldots$ i.i.d.\ from a
distribution $P^*$, where $X_{i} = (X_{i1}, \ldots, X_{i\pmax})$ are
high-dimensional vectors, and we allow $\pmax = \infty$. We use nested
models $\cM_0, \cM_1, \ldots$ where $\cM_p$ is a standard linear
model, consisting of conditional distributions $P(\cdot \mid
\beta,\sigma^2)$ expressing that
\begin{equation}\label{eq:basicequation}
Y_i = \beta_0 + \sum_{j=1}^p  \beta_j X_{ij} + \epsilon_i
\end{equation} 
is a linear function of $p \leq \pmax$ covariates with additive
independent Gaussian noise $\epsilon_i \sim N(0,\sigma^2)$. We equip
each of these models with standard priors on coefficients and the
variance, and also put a discrete prior on the models themselves.
$\cM := \bigcup_{\nc = 0..\pmax} \cM_{\nc}$ does not contain the
conditional `ground truth' $P^*(Y |X)$ (hence the model is
`misspecified'), but it does contain a $\tilde{P}$ that is `best' in
several respects: it is closest to $P^*$ in KL (Kullback-Leibler)
divergence, it represents the true regression function (leading to the
best squared error loss predictions among all $P \in \cM$) and it has
the true marginal variance (explained in Section~\ref{sec:optimalb}).
Yet, while $\tilde{P} \in \cM_0$ and $\cM_0$ receives substantial
prior mass, as $n$ increases, the posterior puts most of its mass on
complex $\cM_p$'s with higher and higher $p$'s, and, conditional on
these $\cM_p$'s, at distributions which are very far from $P^*$ both
in terms of KL divergence and in terms of $L_2$ risk, leading to bad
predictive behavior in terms of squared error.
Figure~\ref{fig:polynomial} and~\ref{fig:polynomialrisk} illustrate a
particular instantiation of our results, obtained when $X_{ij}$ are
polynomial basis functions, i.e.\ $X_{ij} = S_i^{j}$ and $S_i \in
[-1,1]$ uniformly i.i.d.  We also show comparably bad predictive
behavior for various versions of Bayesian ridge regression, involving
just a single, high-but-finite dimensional model. In that case Bayes
eventually recovers and concentrates on $\tilde{P}$, but only at a
sample size that is incomparably larger than what can be expected if
the model is correct.

These findings contradict the folk wisdom that, if the model is
incorrect, then ``Bayes tends to concentrate on neighborhoods of the
distribution(s) in $\cM$ that is/are closest to $P^*$ in KL
divergence.'' Indeed, the strongest actual theorems to this end that
we know of, \citep{KleijnV06,de2013bayesian,ramamoorthi2013posterior},
hold, as the authors emphasize, under regularity conditions that are
substantially stronger than those needed for consistency when the
model is correct (as by e.g.\ \cite{GhosalGV00} or
\cite{Zhang06density}), and our example shows that consistency may
fail to hold even in relatively simple problems.

\paragraph{The Solution: Generalized Posterior and Safe Bayes}
Bayesian updating can be enhanced with a {\em learning rate\/} $\eta$,
an idea put forward independently by several authors
\citep{Vovk90,McAllester02,BarronC91,WalkerH02,Zhang06density} and
suggested as a tool for dealing with misspecification by Gr\"unwald
(\citeyear{Grunwald11,Grunwald12}). $\eta$ trades off the relative
weight of the prior and the likelihood in determining the {\em
  $\eta$-generalized posterior}, where $\eta=1$ corresponds to
standard Bayes and $\eta=0$ means that the posterior always remains
equal to the prior. When choosing the `right' $\eta$, which in our
case is significantly smaller than $1$ but of course not $0$,
$\eta$-generalized Bayes becomes competitive again. In general, the
optimal $\eta$ depends on the underlying ground truth $P^*$, and the
problem has always been how to determine the optimal $\eta$
empirically, from the data. 

Recently,
\cite{Grunwald12} proposed the {\em Safe Bayesian\/} algorithm for
learning $\eta$, and theoretically showed that it achieves good
convergence rates in terms of KL divergence on a variety of
problems. Here we show empirically that Safe Bayes performs
excellently in our regression setting, being competitive with standard
Bayes if the model is correct and very significantly outperforming not
just standard Bayes, but also cross-validation and approaches such as
AIC when the model is incorrect. We do this by providing a wide range
of experiments, varying parameters of the problem such as the priors
and the true regression function and studying various performance
indicators such as the squared error risk, the posterior on the
variance etc.

We note that a Bayesian's (and our) first instinct would be to learn
$\eta$ itself in a Bayesian manner instead. Yet this does not solve
the problem, as we show in Section~\ref{sec:ridge}, where we consider
a setting in which $1/\eta$ turns out to be exactly equivalent to the
$\lambda$ regularization parameter in the Bayesian Lasso and ridge
regression approaches. We find that selecting $\eta$ by (empirical)
Bayes, as suggested by e.g.\ \cite{ParkC08}, does not nearly regularize
enough in our misspecification experiments. In the Bayesian ridge
regression setting with fixed variance, the Safe Bayesian algorithm
becomes very similar to learning $\lambda$ by cross-validation with
squared-error loss, as is standard in frequentist ridge regression
(cross-validation with a logarithmic score does {\em not\/} work
however). In the varying variance case, there is no such
straightforward interpretation of SafeBayes.

\paragraph{The Type of Misspecification}
The models are misspecified in that they make the standard assumption
of homoskedasticity --- $\sigma^2$ is independent of $X$ --- whereas
in reality, under $P^*$, there is heteroskedasticity, there being a
region of $X$ with low and a region with (relatively) high variance.
Specifically, in our simplest experiment the `true' $P^*$ is defined
as follows: at each $i$, toss a fair coin.  If the coin lands heads,
then sample $X_i$ from a uniform distribution on $[-1,1]$, and set
$Y_i = 0 + \epsilon_i$, where $\epsilon_i \sim N(0,\sigma_0^2)$.  If
the coin lands tails, then set $(X_i,Y_i) = (0,0)$, so that there is
no variance at all. The `best' conditional density $\tilde{P}$,
closest to $P^*(Y \mid X)$ in KL divergence, representing the true
regression function $Y=0$ and reliable in the sense of
Section~\ref{sec:optimalb}, is then given by (\ref{eq:basicequation})
with all $\beta$'s set to $0$ and $\tilde{\sigma}^2 = \sigma_0^2/2$.
In a typical sample of length $n$, we will thus have approximately
$n/2$ points with $X_i$ uniform and $Y_i$ normal with mean 0, and
approximately $n/2$ points with $(X_i,Y_i) = (0,0)$. These points seem
`easy' since they lie exactly on the regression function one would
hope to learn; but they really wreak severe havoc.

\paragraph{The In-Liers Cause the Problem}
While it is well-known that in the presence of outliers, Gaussian
assumptions on the noise lead to problems, both for frequentist and
Bayesian procedures, in the present problem we have {\em in-liers\/}
rather than outliers.  Also, if we slightly modify the setup so that
homoskedasticity holds, standard Bayes starts behaving excellently, as
again depicted in Figure~\ref{fig:polynomial}
and~\ref{fig:polynomialrisk}. Finally, while the figure shows what
happens for polynomials, we used independent multivariate $X$'s rather
than nonlinear basis functions in the main experiments below, getting
essentially the same results. All this indicates that the
inconsistency is really caused by misspecification, in particular the
presence of in-liers, and not by anything else.  The setup is inspired
by the work of \cite{GrunwaldL04,GrunwaldL07}, who gave a mathematical
proof that Bayesian inference can be inconsistent under
misspecification in a related but much more artificial classification
setting. Here we show that this can also happen in a much more natural
regression setting.  The setting being more natural, it is also harder
to analyze, and we only demonstrate the inconsistency empirically.
\subsection{Overview of this Paper}
\paragraph{KL-Associated Inference tasks}
Section~\ref{sec:preliminaries} introduces our setting and the main
concepts needed to understand our results. A crucial point here is
that, if Bayesian (or other likelihood-based methods) converge at all
to a distribution in the model $\cM$, this distribution (often called
the `pseudo-truth') is the $\tilde{P} \in \cM$ that minimizes
KL-divergence to the true distribution $P^*$. While the minimum KL
divergence point is often not of intrinsic interest, for some (not
all) models, $\tilde{P}$ can be of interest for other reasons as well
\citep{royall2003interpreting}: there may be {\em associated\/}
inference tasks for which $\tilde{P}$ is suitable as well. For
standard linear models with fixed $\sigma^2$, the main associated task
is squared error prediction: the KL-optimal $\tilde{P}$ is also
optimal, among all $P \in \cM$, in terms of squared error prediction
risk. If additionally $\sigma^2$ becomes a free parameter, then it is
also reliable, which roughly means that it is optimal in determining
its own squared error prediction quality (Section~\ref{sec:optimalb};
we have a lot more to say about associated inference tasks in
Section~\ref{sec:discussion}).  Thus, whenever one is prepared to work
with linear models and one is interested in squared risk or
reliability, then Bayesian inference would seem the way to go, even if
one suspects misspecification\textellipsis at least if there is consistency.

\paragraph{The Safe Bayesian Algorithm}
Section~\ref{sec:generalized} introduces the $\eta$-generalized
posterior and instantiates it to the linear
model. Section~\ref{sec:safebayes} introduces the `Safe Bayesian' algorithm, which
learns $\eta$ from the data. This is done via Dawid's
(\citeyear{Dawid84}) {\em prequential\/} view on Bayesian inference.
We then provide four instantiations of the SafeBayes method to linear models. 

Section~\ref{sec:mainexperiments} discusses our experiments.  We first
provide the necessary preparation in Section~\ref{sec:preparing}
and~\ref{sec:statisticsshown}. Section~\ref{sec:modsel} gives the
results of our first experiment, a comparison of Bayesian and
SafeBayesian model averaging and selection in two settings, one with a
correct model and one with a model corrupted by $50\%$ easy points as
above, but with independent Gaussian rather than polynomial inputs.
Section~\ref{sec:ridge} repeats these experiments for a Bayesian ridge
regression setting, Section~\ref{sec:summary} provides an `executive
summary'. In all experiments Safe Bayesian methods behave much better
in terms of squared error risk and reliability than standard Bayes if the
model is incorrect, and hardly worse (sometimes still better) than
standard Bayes if the model is correct.
\paragraph{Good vs.\ Bad Misspecification: Nonconcentration and Hypercompression}
In and of itself, the fact that one obtains inconsistency with
homoskedastic models and heteroskedastic data may not be very
surprising; and indeed, whether similar phenomena occur in real-world
data needs further study. The main strength of our example is rather
that it clearly shows what can happen in principle, and indicates how
one may go about solving it. We explain this in
Section~\ref{sec:explanation}, in particular on the basis of
Figure~\ref{fig:badmix} on page~\pageref{fig:badmix}, {\em the
  essential picture to understand the phenomenon}. Inconsistency can
only arise under a `bad' form of misspecification, depicted by the
figure. Under bad misspecification, the posterior may {\em fail to
  concentrate}, and this causes trouble. As a theoretical contribution
of this paper, we show in this section that, under some conditions, a
Bayesian strongly believes that her posterior will, in some sense,
concentrate fast. Indeed, SafeBayes will only select $\eta \ll 1$ if
the standard posterior is nonconcentrated, and may thus be (loosely)
viewed as a particular `prior predictive check'. 

Posterior nonconcentration in turn can lead to `hypercompression', the
phenomenon that the Bayes predictive distribution behaves {\em
  substantially better\/} under a logarithmic scoring rule than the
best distribution $\tilde{P} \in \cM$; this can happen because the
Bayes predictive distribution --- a mixture of elements of $\cM$ ---
behaves substantially differently from any of the elements of $\cM$.
Somewhat paradoxically (Section~\ref{sec:hypercompression}), Bayes'
overly good log-loss behavior is exactly what causes it to perform
badly for the associated inference tasks (squared error prediction and
reliability, in our case). Thus, there can be an inherent tension
between behavior under log-loss and behavior under its associated
tasks, a discrepancy which one can measure by the {\em mixability
  gap\/} (Section~\ref{sec:mixability}), a theoretical concept
introduced by \cite{ErvenGKR11} and \cite{Grunwald12}.  If one is
interested in log-loss, standard Bayes is just fine; the Safe Bayesian
algorithm should be used if one wants to optimize behavior against
the associated tasks.  Of course, whether such a task-dependent
modification of Bayes is desirable needs discussion, which we provide
in Section~\ref{sec:discussion}.

\paragraph{Additional Experiments}
The paper is followed by a long list of appendices where we provide a
battery of experiments to check the robustness of our
results. Specifically, we investigate what happens if we vary our
models and priors (using e.g.\ a fixed $\sigma^2$ and standard priors
used in the regression literature), our methods, and if we vary the
data generating distribution using e.g.\ `easy' points that are close
to, but not exactly $(0,0)$. Our main conclusion here is that, of the
four versions of SafeBayes which we propose, one is uncompetitive and
among the other three, there is no clear winner --- although they
consistently outperform Bayes under misspecification.  Furthermore we
show that AIC, BIC and cross-validation also have serious problems in
our regression setup. We also provide a proof for the theorem about
nonconcentration given in Section~\ref{sec:mixability}.

\section{Preliminaries: Setting, Optimal KL Distribution, Regression Function}
\label{sec:preliminaries}
\subsection{Setting, Logarithmic Risk, Optimal Distribution}
\label{sec:optimal}
In this paper we consider data $Z^n = Z_1, Z_2, \ldots, Z_n \sim$
i.i.d.~$P^*$, where each $Z_i = (X_i,Y_i)$ is an independently sampled
copy of $Z = (X,Y)$, $X$ taking values in some set $\cX$, $Y$ taking
values in $\cY$ and $\cZ = \cX \times \cY$. We are given a {\em
  model\/} $\cM = \{ P_{\theta} \mid \theta \in \Theta\}$
parameterized by (possibly infinite-dimensional) $\Theta$, and
consisting of conditional distributions $P_{\theta}(Y \mid X)$,
extended to $n$ outcomes by independence. For simplicity we assume
that all $P_{\theta}$ have corresponding conditional densities
$\dens_{\theta}$, and similarly, the conditional distribution $P^*(Y
\mid X)$ has a conditional $\dens^*$, all with respect to the same
underlying measure. While we do not assume $P^*(Y \mid X)$ to be in
(or even `close' to) $\cM$, we want to learn, from given data $Z^n$, a
`best' (in a sense to be defined below) element of $\cM$, or at least,
a distribution on elements of $\cM$ that can be used to make
predictions about future data. While our experiments focus on linear
regression, the discussion in this section holds for general
conditional density models. The logarithmic score, henceforth
abbreviated to {\em log-loss}, is defined in the standard manner: the loss
incurred when predicting $Y$ based on density $\dens(\cdot \mid x)$
and $Y$ takes on value $y$, is given by $ - \log \dens(y|x)$. A central
quantity in our setup is then the {\em expected log-loss\/} or {\em
  log-risk}, defined as
$$
\logrisk(\theta) := \Exp_{(X,Y) \sim P^*} [- \log \dens_{\theta}(Y \mid
X)],
$$
where here as in the remainder of this paper, $\log$ denotes natural
logarithm.

We let $P^*_X$ be the marginal distribution of $X$ under $P^*$.  The
Kullback-Leibler (KL) divergence $D(P^* \| P_{\theta})$ between $P^*$ and
conditional distribution $P_{\theta}$ is defined as the expected KL-divergence,
under $X \sim P^*_X$, of the KL divergence $D(P^*(\cdot \mid X ) \|
P_{\theta}(\cdot \mid X))$ between $P_{\theta}$ and the `true' conditional 
$P^*(Y|X)$: $D(P^* \| P_{\theta}) = E_{X \sim P^*_X}[ D(P^*(\cdot |X) \|
P_{\theta}(\cdot |X))]$. A simple calculation shows that for any
$\theta$, $\theta'$,
$$
D(P^* \| P_{\theta}) - D(P^* \| P_{\theta'}) = \logrisk(\theta) - \logrisk(\theta'),
$$
so that the closer $P_{\theta}$ is to $P^*$ in terms of KL divergence, the
smaller its log-risk, and the better it is, on average, when used for
predicting under the log-loss.

Now suppose that $\cM$ contains a unique distribution that is closest,
among all $P \in \cM$ to $P^*$ in terms of KL-divergence. We denote
such a distribution, if it exists, by $\tilde{P}$. Then $\tilde{P} =
P_{\theta}$ for at least one $\theta \in \Theta$; we pick any such
$\theta$ and denote it by $\tilde\theta$, i.e.~$\tilde{P} =
P_{\tilde\theta}$, and note that it also minimizes the log-risk:
\begin{equation}\label{eq:optimal}
\logrisk(\tilde{\theta}) =  \min_{\theta \in \Theta} \logrisk(\theta) = 
\min_{\theta \in \Theta} \Exp_{(X,Y) \sim P^*} [- \log \dens_{\theta}(Y \mid X)].
\end{equation}
We shall call such a $\tilde\theta$ {\em optimal}.

Since, in regions of about equal prior density, the log Bayesian
posterior density is proportional to the log likelihood ratio, we hope
that, given enough data, with high $P^*$-probability, the posterior
puts most mass on distributions that are close to $P_{\tilde\theta}$
in KL-divergence, i.e.\ that have log-risk close to optimal. Indeed,
all existing consistency theorems for Bayesian inference under
misspecification express concentration of the posterior around
$P_{\tilde\theta}$. 
\subsection{A Special Case: The Linear Model}
\label{sec:linear}
Fix some $\pmax \in \{0, 1, \ldots \} \cup \{ \infty \}$. We observe
data $Z_1, \ldots, Z_n$ where $Z_i = (X_i,Y_i)$, $Y_i \in \reals$ and
$X_i = (1,X_{i1}, \ldots, X_{i \pmax}) \in \reals^{\pmax+1}$. Note
that this is as in (\ref{eq:basicequation}) but from now on we adopt
the standard convention to take $X_{0i} \equiv 1$ as a dummy random
variable.  We denote by $\cM_\nc = \{ P_{\nc,\beta,\sigma^2} \mid
(\nc,\beta,\sigma^2) \in \Theta_{\nc} \}$ the standard linear model
with parameter space $\Theta_\nc := \{(\nc,\beta,\sigma^2) \mid \beta
= (\beta_0, \ldots, \beta_{\nc})^T \in \reals^{\nc+1}, \sigma^2 > 0\}$,
where the entry $\nc$ in $(\nc,\beta,\sigma^2)$ is redundant but included for notational
convenience.  We let $\Theta = \bigcup_{\nc= 0 .. \pmax} \Theta_{\nc}
$.  $\cM_{\nc}$ states that for all
$i$, (\ref{eq:basicequation}) holds, where $\epsilon_1, \epsilon_2,
\ldots \sim \text{i.i.d.~}N(0,\sigma^2)$.   When working with linear
models $\cM_{\nc}$, we are
usually interested in finding parameters $\beta$ that predict well in
terms of the {\em squared error loss function\/} (henceforth
abbreviated to {\em square-loss\/}): the square-loss on data
$(X_i,Y_i)$ is $(Y_i - \sum_{j = 0}^{\nc} \beta_j X_{ij})^2 = (Y_i -
X_i \beta)^2$. We thus want to find the
distribution minimizing the expected square-loss, i.e.\ {\em squared
  error risk\/} (henceforth abbreviated to `square-risk') relative to
the underlying $P^*$:
\begin{equation}\label{eq:squarederrorrisk}
\squaredrisk(\nc,\beta) := \Exp_{(X,Y) \sim P^*}(Y-
\Exp_{\nc,\beta,\sigma^2} [Y \mid X])^2 = \Exp_{(X,Y) \sim P^*}(Y - \sum_{j=0}^{\nc} \beta_j X_j)^2,
\end{equation}
where $\Exp_{\nc,\beta,\sigma^2} [Y \mid X]$ abbreviates $\Exp_{Y \sim
  P_{\nc,\beta,\sigma^2} \mid X}[Y]$. Since this quantity is
independent of the variance $\sigma^2$, $\sigma^2$ is not used as an argument
of $\squaredrisk$.

\subsection{KL-Associated Prediction Tasks for the Linear Model: The
  KL-Optimal $\tilde{\theta} = (\tilde\beta,\tilde\nc,\tilde\sigma^2)$
  is square-risk optimal and reliable}%
\label{sec:optimalb}
Suppose that an optimal $\tilde{P} \in \cM$ exists in the regression model. We denote by $\tilde{\nc}$ the smallest $\nc$ such that $\tilde{P} \in \cM_{\nc}$, and define $\tilde{\sigma}^2, \tilde{\beta}$ such that $\tilde{P} = P_{\tilde{\nc}, \tilde{\beta}, \tilde{\sigma}^2}$. A straightforward computation
shows that for all $(\nc,\beta,\sigma^2) \in \Theta$:
\begin{equation}\label{eq:THEidentity}
\logrisk( (\nc,\beta,\sigma^2) ) =  \frac{1}{2 \sigma^2} \squaredrisk(
(\nc,\beta) ) + \frac{1}{2} \log (2 \pi \sigma^2),  
\end{equation}
so that the $(\nc,\beta)$ achieving minimum log-risk for each fixed
$\sigma^2$ is equal to the $(\nc,\beta)$ with the minimum
square-risk. In particular, $(\tilde\nc,\tilde\beta,\tilde\sigma^2)$
must minimize not just log-risk, but also square-risk.  Moreover, the
conditional expectation $\Exp_{P^*}[Y \mid X]$ is known as the {\em
  true regression function}. It minimizes the square-risk among all
conditional distributions for $Y \mid X$. Together with
(\ref{eq:THEidentity}) this implies that, if there is some
$(\nc,\beta)$ such that $\Exp[Y \mid X] =  \sum_{j=0}^{\nc}
\beta_j X_j = X \beta$, i.e.\ $(\nc,\beta)$ represents the true regression
function, then $(\tilde{\nc},\tilde{\beta})$ also represents the true
regression function. In all our examples, this will be the case: the
model is misspecified only in that the true noise is heteroskedastic;
but the model does invariably contain the true regression function.

Moreover, for each fixed $(\nc,\beta)$, the $\sigma^2$
minimizing $\logrisk$ is, as follows by differentiation, given by
$\sigma^2 = \squaredrisk(\nc,\beta)$. In particular, this implies that 
\begin{equation}
\tilde{\sigma}^2 = \squaredrisk(\tilde\nc,\tilde{\beta}),
\end{equation} 
or in words: the KL-optimal model variance $\tilde{\sigma}^2$ is equal
to the true expected (marginal, not conditioned on $X$) square-risk
obtained if one predicts with the optimal
$(\tilde{p},\tilde\beta)$. This means that the optimal
$(\tilde{p},\tilde\beta,\tilde\sigma^2)$ is {\em reliable\/} in the
sense of \cite{Grunwald98b,Grunwald99a}: its self-assessment about its square-loss performance
is correct, independently of whether
$\tilde\beta$ is equal to the true regression function or not: {\em
  $(\tilde{p},\tilde\beta,\tilde\sigma^2)$ correctly predicts how well
  it predicts}.

Summarizing, for misspecified models,
$(\tilde\nc,\tilde\beta,\tilde{\sigma}^2)$ is optimal not just in
KL/log-risk sense, but also in terms of square-risk and in terms of
reliability; in our examples, it also represents the true regression
function. We say that, for linear models, square-risk optimality,
square-risk reliability and regression-function consistency are {\em
  KL-associated prediction tasks\/}: if (as we hope Bayes will do, but
as we will see sometimes won't) we can find the KL-optimal
$\tilde\theta$, we automatically behave well in these associated tasks
as well.
\section{The Generalized Posterior}\label{sec:generalized}
\paragraph{General Losses} The original generalized posterior is a
notion going back at least to \cite{Vovk90} and has been developed
mainly within the so-called (frequentist) {\em PAC-Bayesian\/}
framework
\cite{McAllester02,Seeger02,Catoni07,Audibert04,Zhang06loss}; see also
\cite{bissiri2016general} and the discussion in
Section~\ref{sec:discussion}.  It is defined relative to a prior on
{\em predictors\/} rather than probability distributions.  Depending
on the decision problem at hand, predictors can be e.g.\ classifiers,
regression functions or probability densities.  Formally, we are given
an abstract space of predictors represented by a set $\Theta$, which
obtains its meaning in terms of a loss function $\ell: \cZ \times
\Theta \rightarrow \reals$, writing $\ell_{\theta}(z)$ as shorthand
for $\ell(z,\theta)$. Following e.g.\ \cite{Zhang06loss}, for any
prior $\Pi$ on $\Theta$ with density $\pi$ relative to some underlying
measure $\rho$, we define the {\em generalized Bayesian posterior with
  learning rate $\eta$ relative to loss function $\ell$}, denoted as $
\Pi \mid Z^n, \eta $, as the distribution on $\Theta$ with density
\renewcommand{\L}{\ell}
\begin{equation}\label{eq:genpost}
\pi(\theta \mid z^n, \eta) := 
\frac{e^{- \eta \sum_{i=1}^n \L_{\theta}(z_i)} \pi(\theta)}{\int 
e^{- \eta \sum_{i=1}^n \L_{\theta}(z_i)} \pi(\theta) \rho(d \theta)}
= \frac{e^{- \eta 
\sum_{i=1}^n \L_{\theta}(z_i)
} \pi(\theta)
}{\Exp_{\theta \sim \Pi} [ e^{- \eta  
\sum_{i=1}^n \L_{\theta}(z_i)
}]}.
\end{equation}
Thus, if $\theta_1$ fits the data better than $\theta_2$ by a difference of $\epsilon$
according to loss function $\ell$, then their posterior ratio is
larger than their prior ratio by an amount exponential in $\epsilon$,
where the larger $\eta$, the larger the influence of the data as
compared to the prior. 

If $z_i = (x_i,y_i)$ with $y_i \in \reals$ and $x_i = (1,x_{i1}, \ldots,
x_{i\nc})$, and the
goal is to predict $y_i$ given $x_i$, then we may
take as our prediction model e.g.\ the set of linear predictors that
predict $y_i$ by $\sum \beta_j x_{ij} = x_i \beta$, and as our loss
function the squared error loss, $\ell_{\beta}(x_i,y_i) = (y_i -
x_i\beta)^2$. We may then study the behavior of such a procedure in
its own right, irrespective of a Bayesian misspecification
interpretation; the experiments we perform in Appendix~\ref{sec:mainfixedsigma} and~\ref{sec:ridgefixedsigma} can
be interpreted in this manner. 
\paragraph{Log-Loss and Likelihood} 
Now if the set $\Theta$ represents a model of (conditional)
distributions $\cM = \{ P_{\theta} \mid \theta \in \Theta \}$, we may
set, for $z_i = (x_i, y_i)$, $\ell_{\theta}(z_i) = - \log
\dens_{\theta}(y_i \mid x_i)$ to be the log-loss as defined above. In this
special case, the definition of $\eta$-generalized posterior
specializes to the definition of `generalized posterior' as known
within the Bayesian literature \citep{WalkerH02,Zhang06density}: 
\begin{equation}\label{eq:genpostb}
\pi(\theta \mid z^n, \eta) = 
 \frac{(\/\dens(y^n \mid x^n,\theta) \/)^{\eta} \pi(\theta)}{\int
(\/ \dens(y^n \mid x^n,\theta) \/)^{\eta} 
\pi(\theta) \rho(d \theta)}
= \frac{
(\/ \dens(y^n \mid x^n,\theta) \/)^{\eta}  \pi(\theta)
}{\Exp_{\theta \sim \Pi} [ (\/ \dens(y^n \mid x^n,\theta) \/)^{\eta} ] }.
\end{equation}
Again, the larger $\eta$, the larger the influence of the likelihood.
Obviously $\eta = 1$ corresponds to standard Bayesian inference,
whereas if $\eta = 0$ the posterior is equal to the prior and nothing
is ever learned. Our algorithm for learning $\eta$ will usually end up
with values in between. It has long been known that in model selection
and nonparametric settings, there is an issue with consistency proofs
for full Bayes, Bayes MAP and MDL if we take the standard $\eta = 1$,
and indeed, this is part of the reason why the generalized posterior
in the form (\ref{eq:genpostb}) was derived in the first place: for
example, \cite{BarronC91} give general consistency theorems for 2-part
MDL (closely related to Bayes MAP) and note that they hold for any
$\eta < 1$; but for $\eta = 1$, additional assumptions must be made.
\citet{Zhang06density} gives an explicit example in which the
posterior shows anomalous behavior at $\eta=1$.  A connection to
misspecification was first made by \cite{Grunwald11} (see
Section~\ref{sec:related}) and \cite{Grunwald12}.

\paragraph{Generalized Predictive Distribution}
We also define the predictive distribution based on the
$\eta$-generalized posterior (\ref{eq:genpostb}) as a generalization of the standard
definition as follows: for $m \geq 0, m'\geq m$, we set
\begin{align}\label{eq:bayespred}
\pbayes(y_i,\ldots, y_{i+ m} \mid x_i, \ldots, x_{i+m'}, z^{i-1}, \eta)& 
:= \Exp_{\theta \sim \Pi \mid
  z^{i-1}, \eta} [ \dens(y_i, \ldots, y_{i+ m} \mid x_i, \ldots, x_{i+
  m'}, \theta)] \nonumber \\ &\phantom{:}= 
\Exp_{\theta \sim \Pi \mid z^{i-1}, \eta} [
\dens(y_i, \ldots, y_{i+ m} \mid x_i, \ldots, x_{i+ m}, \theta)].
\end{align}
where the first equality is a definition and the second follows by our
i.i.d.\ assumption. 
We always use the bar-notation $\pbayes$ to indicate marginal and
predictive distributions, i.e.\ distributions on data that are arrived
at by integrating out parameters. If $\eta =1$ then $\pbayes$ and
$\pi$ become the standard Bayesian predictive density and
posterior, and if it is clear from the context that we consider $\eta
= 1$, we leave out the $\eta$ in the notation.
\begin{ownquote}
  The generalized posterior is created by exponentiating the
  likelihood according to individual elements $\theta \in \Theta = \bigcup_\nc
  \Theta_\nc$ in the model and renormalizing, which is not the same as
  exponentiating marginal likelihoods and renormalizing. In
  particular, $\pi(\nc \mid z^n, \eta)$ as given by
  (\ref{eq:modelposteriorb}) is in general {\em not\/} proportional to $
  (\/\pbayes(y^n \mid x^n, \nc)\/)^{\eta} \pi(\nc)$.  Similarly, for
  generalized marginal distributions, as soon as $\eta \neq 1$, we
  have that in general
$$\pbayes(y_i, y_{i+1} \mid x_i, x_{i+1}, z^{i-1}, \eta) \neq
\pbayes(y_i \mid x_i, z^{i-1}, \eta) \cdot \pbayes(y_{i+1} \mid
x_{i+1}, z^{i}, \eta), $$ unlike for the standard Bayesian marginal
distribution for which equality holds (in Section~\ref{sec:howitworks} we
encounter a further modification of the generalized posterior whose
marginals do satisfy this product rule).
\end{ownquote}

\subsection{Instantiation to Linear Model Selection and Averaging}
\label{sec:instantiation}
Now consider again a linear model $\cM_{\nc}$ as defined in
Section~\ref{sec:optimalb}. We instantiate the generalized posterior
and its marginals for this model.
With prior $\pi(\beta,\sigma^2  \mid \nc)$ taken relative to Lebesgue measure,
(\ref{eq:genpostb}) specializes to:
$$
\pi(\beta,\sigma \mid z^n, \nc, \eta) = 
\frac{(2 \pi \sigma^2)^{-n \eta /2}
e^{- \frac{\eta}{2 \sigma^2}  \sum_{i=1}^n (y_i - x_i \beta)^2 } \pi(\beta,\sigma 
\mid \nc)}{\int (2 \pi \sigma^2)^{-n \eta /2}
e^{- \frac{\eta}{2 \sigma^2}  \sum_{i=1}^n (y_i - x_i \beta)^2 } \pi(\beta,\sigma 
\mid \nc) \,d\beta \,d \sigma
}.
$$
Note that in the numerator $1/\sigma^2$ and $\eta$ are interchangeable
in the exponent, but not in the factor in front: their role is subtly
different. For Bayesian
inference with a sequence of models $\cM = \bigcup_{\nc=0..\pmax}
\cM_{\nc}$, with $\pi(\nc)$  a probability mass function on $\nc \in \{0, \ldots,
\pmax\}$, we get: 
\begin{align}\label{eq:genpostd}
\pi(\/ \theta\/ \mid z^n  ,\eta)
&=   
\frac{\dens(y^n \mid x^n, \theta )^{\eta} 
\pi(\theta)}{\int_{\theta \in \Theta} 
\dens(y^n \mid x^n,\theta)^{\eta} 
\pi(\theta) \rho( d\theta)} \text{\ \ with \ \  } \theta = (\beta,\sigma^2,\nc)
\nonumber \\ = \pi(\/ \beta,\sigma,\nc \/ \mid z^n  ,\eta)
&= \frac{(2 \pi \sigma^2)^{-n \eta /2}
e^{- \frac{\eta}{2 \sigma^2}  \sum_{i=1}^n (y_i - x_i \beta)^2 } \pi(\beta,\sigma 
\mid \nc)\pi(\nc)}{\sum_{\nc=0}^{\pmax} \int (2 \pi \sigma^2)^{-n \eta /2}
e^{- \frac{\eta}{2 \sigma^2}  \sum_{i=1}^n (y_i - x_i \beta)^2 } \pi(\beta,\sigma 
\mid \nc) \pi(\nc) \,d\beta \,d \sigma
}
\end{align}
The total generalized posterior probability of model $\cM_\nc$ then becomes: 
\begin{equation}\label{eq:modelposteriorb}
\pi (\nc  \mid z^n, \eta) = \int \pi(\beta,\sigma,\nc \mid z^n, \eta) \,d
\beta \,d \sigma .
\end{equation}
Analogously to (\ref{eq:bayespred}), for given $\nc$, we define (writing $a_i^j$ as shorthand for $a_i, \ldots, a_{j}$), the
$\eta$-generalized Bayesian predictive distribution as:
\begin{align}
\pbayes(y_i^{i+m} \mid x_i^{i+m'}, z^{i-1}, \nc, \eta) 
& := \Exp_{\beta,\sigma^2 \sim \Pi \mid
  z^{i-1}, \nc, \eta} [ \dens(y_i^{i+ m} \mid x_i^{i+
  m'}, \beta,\sigma^2, \nc)] \nonumber \\ & \phantom{:}= %
\Exp_{\beta,\sigma^2 \sim \Pi \mid z^{i-1},\nc, \eta} [
\dens(y_i^{i+ m} \mid x_i^{i+ m},
\beta,\sigma^2, \nc)]. \label{eq:ci}
\end{align}
The previous displays held for general priors. The experiments in this
paper adopt widely used priors (see e.g.\ \cite{raftery1997bayesian}):
normal priors on the $\beta$'s and inverse gamma priors on the
variance.  These conjugate priors allow explicit analytical formulas
for all relevant quantities for arbitrary $\eta$, provided below. We
only consider the simple case of a fixed $\cM_{\nc}$ here; the more
complicated formulas with an additional prior on $\nc$ are given in
Appendix~\ref{app:mixability}.
\paragraph{Fixed $p$ and $\sigma^2$}
Let ${\bf X}_n = (x_1^T, \ldots, x_n^T)^T$ be the design matrix. For a
linear model $\cM_\nc$ with fixed variance $\sigma^2$ and initial
Gaussian prior on $\beta$ given by $N(\betao, \sigma^2 \Sigma_0)$, the
generalized posterior on $\beta$ is again Gaussian with mean
\begin{equation}\label{eq:olvg}
\bar{\beta}_{n,\eta} := \Exp_{\beta \sim \Pi \mid z^n, \nc, \eta} \beta
= \Sigma_{n,\eta} (\Sigma_0^{-1} \betao + \eta {\bf X}_n^T y^n)
\end{equation}
and covariance matrix $\sigma^2 \Sigma_{n,\eta}$, where
$\Sigma_{n,\eta} = (\Sigma_0^{-1} + \eta {\bf X}_n^T {\bf X}_n)^{-1}$. 
\paragraph{Fixed $p$, varying $\sigma^2$}
Now consider linear models with a Gaussian prior on $\beta$
conditional on $\sigma^2$ as above, and a conjugate (inverse gamma)
prior on $\sigma^2$, i.e.\ $\pi(\sigma^2) = \InvGamma(\sigma^2 \mid
a_0, b_0)$ for some $a_0$ and $b_0$. Here we use the following
parameterization of the inverse gamma distribution:
\begin{equation}\label{eq:inversegamma}
\InvGamma(\sigma^2 \mid a, b)
= \sigma^{-2(a+1)} e^{-b / \sigma^2} b^a / \Gamma(a).
\end{equation}
The posterior $\pi(\sigma^2, z^n, p)$ is then given by $\InvGamma(\sigma^2 \mid
a_{n,\eta}, b_{n,\eta})$ where
\begin{equation}\label{eq:anbn}
a_{n,\eta} = a_0 + \eta n / 2\ \ ; \ \  
b_{n,\eta} = b_0 + \frac{1}{2} \betao^T \Sigma_0^{-1} \betao + \frac{\eta}{2} (y^n)^T y^n - \frac{1}{2} \bar\beta_{n,\eta}^T \Sigma_{n,\eta}^{-1} \bar\beta_{n,\eta}.
\end{equation}
The posterior expectation of $\sigma^2$ can be calculated as
\begin{equation}\label{eq:expsigma}
\bar{\sigma}^2_{n,\eta} := \frac{b_{n,\eta}}{a_{n,\eta} - 1}.
\end{equation}
Note that the posterior mean of $\beta$ given $\sigma^2$ does not depend on $\sigma^2$.
\section{The Safe Bayesian Algorithm}
\label{sec:safebayes}
\subsection{Introducing Safe Bayes via the Prequential View}
We introduce SafeBayes via Dawid's prequential interpretation of Bayes
factor model selection. As was first noticed by
\cite{Dawid84} and \cite{Rissanen84}, we can think of Bayes factor model
selection as picking the model with index $\nc$ that, when used for
sequential prediction with a logarithmic scoring rule, minimizes the
cumulative loss. To see this, note that for any distribution whatsoever, we have that, by definition of conditional probability, 
$$
- \log \dens(y^n) =  - \log \prod_{i=1}^n \dens(y_i \mid y^{i-1}) =
\sum_{i=1}^n - \log \dens(y_i \mid y^{i-1}). 
$$
In particular, for the standard Bayesian marginal distribution
$\pbayes(\cdot \mid \nc) = \pbayes(\cdot \mid \nc, \eta = 1)$ as
defined above, for each fixed $\nc$, we have
\begin{equation}
\label{eq:snack}
- \log \pbayes(y^n \mid x^n, \nc) = \sum_{i=1}^n - \log \pbayes(y_i \mid x^n, y^{i-1}, \nc)
= \sum_{i=1}^n - \log \pbayes(y_i \mid x_i, z^{i-1}, \nc) ,
\end{equation}
where the second equality holds by (\ref{eq:ci}).
If we assume a uniform prior on model index $\nc$, then Bayes factor
model selection picks the model maximizing $\pi(\nc \mid z^n)$, which
by Bayes' theorem coincides with the model minimizing (\ref{eq:snack}),
i.e.\ minimizing cumulative log-loss. Similarly, in `empirical Bayes'
approaches, one picks the value of some nuisance parameter $\rho$ that
maximizes the marginal Bayesian probability $\pbayes(y^n \mid x^n,
\rho)$ of the data. By (\ref{eq:snack}), which still holds with $\nc$
replaced by $\rho$, this is again equivalent to the $\rho$ minimizing
the cumulative log-loss. This is the {\em prequential\/}
interpretation of Bayes factor model selection and empirical Bayes
approaches, showing that Bayesian inference can be interpreted as a
sort of {\em forward\/} (rather than cross-) validation
\citep{Dawid84,Rissanen84,Hjorth82}.

We will now see whether we can use this approach with $\rho$ in the
role of the $\eta$ for the $\eta$-generalized posterior that we want to
learn from the data. 
We continue to rewrite (\ref{eq:snack}) as follows (with $\rho$
instead of $p$
that can either stand for a continuous-valued parameter or for a model
index but not yet for $\eta$), using the fact that the Bayes predictive distribution given
$\rho$ and $z^{i-1}$ can be rewritten as a posterior-weighted average of
$\dens_{\theta}$:
\begin{align}\label{eq:sunshine}
\breve{\rho} & := \arg \max_{\rho}  \pbayes(y^n \mid x^n,  \rho) = \arg \min_{\rho}
\sum_{i=1}^n \left(\/ - \log \pbayes(y_i \mid x_i,z^{i-1}, \rho)
  \/\right) \nonumber \\ & \phantom{:}= 
\arg \min_{\rho} \sum_{i=1}^n \left(\/ - \log \Exp_{\theta \sim \Pi \mid z^{i-1}, {\rho}}
[ \dens(y_i \mid x_i,\theta)] \/\right). 
\end{align}
This choice for $\breve{\rho}$ being entirely consistent with the
Bayesian approach, our first idea is to choose $\hat\eta$ in the same
way: we simply pick the $\eta$ achieving (\ref{eq:sunshine}), with
$\rho$ substituted by $\eta$. However as Figure~\ref{fig:etafunction} will show (the blue line there depicts (\ref{eq:sunshine}) for one of our experiments),
this will tend to pick $\eta$ close to $1$ and does not improve
predictions under misspecification. Indeed, we introduced $\eta$ to
deal with the case in which the Bayesian model assumptions are
violated, so we cannot expect that learning it in a Bayes-like
way such as (\ref{eq:sunshine}) will resolve the issue. But it turns
out that a {\em slight\/} modification of (\ref{eq:sunshine}) does the
trick: we simply interchange the order of logarithm and expectation in
(\ref{eq:sunshine}) and pick the $\eta$ minimizing
\begin{equation}\label{eq:kibbeling}
\sum_{i=1}^n \Exp_{\theta \sim \Pi \mid z^{i-1}, \eta} \left[ - \log \dens(y_i \mid x_i,\theta) \right].
\end{equation}
In words, we pick the $\eta$ minimizing the {\em P\/}osterior-{\em E\/}xpected
{\em P\/}osterior-{\em R\/}andomized log-loss, i.e.\ the log-loss we expect to obtain,
according to the $\eta$-generalized posterior, if we actually sample from this
posterior. This modified loss function has also been called {\em Gibbs
  error\/} \citep{Cuong13}, and while the abbreviation {\em PEPR\/}-log-loss
would be more correct, we simply call it the $\eta$-$R$-log-loss from now on.

A detailed explanation of why this works will have to wait until
Section~\ref{sec:hypercompression} and~\ref{sec:mixability}; for now we just notice that by Jensen's
inequality, for any fixed $\eta$, for every sequence of data we must
have
\begin{equation}\label{eq:jensen}
  \Exp_{\theta \sim \Pi \mid z^{i-1}, \eta} \left[ - \log \dens(y_i \mid x_i,\theta) \right] \geq - \log \Exp_{\theta \sim \Pi \mid z^{i-1}, \eta} \left[ \dens(y_i \mid x_i,\theta) \right],\end{equation} 
yet, the difference between both sides is small if the posterior is
{\em concentrated\/} for $(x_i,y_i)$, i.e.\ for small $\epsilon$ and small
positive $\delta$, it puts $1- \delta$ of its
mass on distributions which assign the same density to $y_i$ given
$x_i$ up to a factor $1 + \epsilon$ --- clearly, if $\delta = \epsilon = 0$ then both
sides are the same. Thus, at values for $\eta$ at which the
generalized posterior is `cumulatively concentrated', i.e.\ concentrated
at most sample points, the objective function will be similar to the
standard Bayesian one. This is the clue to further analysis of the algorithm to follow later. 

In practice, it is computationally infeasible to try all values of
$\eta$ and we simply have to try out a number of values.  For
convenience we give a detailed description of the resulting algorithm
below, copied from \cite{Grunwald12}. In this paper, we will
invariably apply it with $z_i = (x_i,y_i)$ as before, and
$\ell_{\theta}(z_i)$ set to the (conditional) log-loss as defined
before, although it sometimes also has a second interpretation with
$\ell_{\theta}$ as square-loss.
  
\begin{algorithm}[h]\SetAlgoLined \label{alg:algorithm}
  \KwIn{data $z_1, \ldots, z_n$, model $\model = \{ \dens(\cdot \mid
    \theta) \mid \theta \in \Theta \}$, prior $\Pi$ on $\Theta$,
    step-size $\kappa_{\text{\sc step}}$, max.\ exponent $\kappa_{\max}$, loss function
    $\ell_{\theta}(z)$}
\KwOut{Learning rate $\hat\eta$} 
${\cal S}_n := \{ 1, 2^{-\kappa_{\text{\sc step}}},
2^{-2\kappa_{\text{\sc step}}}, 2^{-3 \kappa_{\text{\sc step}}}, \ldots, 2^{-
  \kappa_{\max}}\}$\; \For{all $\eta \in {\cal S}_n$}{ $s_{\eta} :=
  0$\; \For {$i = 1 \ldots n$}{ Determine generalized posterior
    $\Pi(\cdot \mid z^{i-1}, \eta)$ of
    Bayes with learning rate $\eta$. \\
    Calculate ``posterior-expected posterior-randomized loss'' of
    predicting actual next outcome:
\begin{equation}\label{eq:rloss}
r := 
\L_{\Pi \mid z^{i-1}, \eta}(z_i)
 = E_{\theta \sim \Pi\mid z^{i-1}, \eta}\left[ 
\L_{\theta}(z_i)
\right] 
\end{equation}
$s_{\eta}  := s_{\eta} + r$\;
}}
Choose $\hat{\eta} := \arg \min_{\eta \in {\cal S}_n } \{ s_{\eta} \}$ 
(if $\min$ achieved for several  $\eta\in \cS_n$, pick  largest)\;
\caption{The ($R$-) Safe Bayesian Algorithm.}
\end{algorithm}
\paragraph{Variation}
As we will see in Section~\ref{sec:mixability}, the crucial property
to make inference about $\eta$ work is that the expression inside the
sum in (\ref{eq:sunshine}) is replaced by
\begin{equation}\label{eq:haring}
\Exp_{\theta \sim \Pi'} [ -
\log f_{\theta}(Y_i \mid X_i)]
\end{equation} where $\Pi'$ should be chosen such that 
the resulting log-loss is as small as possible. In
(\ref{eq:kibbeling}) we set $\Pi' = \Pi$, but $\Pi'$ is allowed to be {\em any\/}
distribution on $\theta$ under which the expected log-loss is small. The heuristic analysis of Section~\ref{sec:mixability} 
suggests that the smaller the loss that can be formed this way (see also
under `Open Problems' in Section~\ref{sec:discussion}), the
better the resulting method is expected to work. 

Now by Jensen's inequality, the $\eta$-{\em in-model-log-loss\/} or just
$\eta$-{\em $I$-log-loss}, defined as,
\begin{equation}\label{eq:jibbeling}
\sum_{i=1}^n  \left[ - \log \dens(y_i \mid x_i,\Exp_{\theta \sim \Pi \mid z^{i-1}, \eta} [\theta]) \right],
\end{equation}
is always smaller than (\ref{eq:kibbeling}) for the linear models that
we consider. This means that, instead of finding the $\eta$ minimizing
(\ref{eq:kibbeling}), we may want to find the $\eta$ minimizing
(\ref{eq:jibbeling}), which is of the form (\ref{eq:haring}) with
$\Pi'$ equal to a point mass on $\bar{\theta}_{i,\eta} := \Exp_{\theta
  \sim \Pi \mid z^{i-1}, \eta} \dens_{\theta}$.  We call the version
of SafeBayes which minimizes the alternative objective function
(\ref{eq:jibbeling}) {\em in-model SafeBayes}, abbreviated to
$I$-SafeBayes, and from now on use $R$-SafeBayes for the original
version based on the $R$-log-loss.  We did not realize
the potential benefits of using in-model SafeBayes at the time of
writing \cite{Grunwald12}, and while the theoretical results of
\cite{Grunwald12} can be adjusted to deal with such modifications, we
cannot get any better theoretical convergence bounds as yet, but this may be
an artifact of our proof techniques. A secondary goal of the
experiments in this paper is thus to see whether one can really improve
SafeBayes by using the `in-model' version.

\subsection{Instantiating SafeBayes to the Linear Model}
\label{sec:instantiationb}
Our experiments concern four instantiations of SafeBayes: 
$R$-SafeBayes and $I$-SafeBayes for models with fixed variance,
denoted {\em $R$-square-SafeBayes and $I$-square-SafeBayes\/} for
reasons that will become clear below, are the topic of experiments in
Appendix~\ref{sec:mainfixedsigma} and~\ref{sec:ridgefixedsigma}. The
main text instead investigates, in Section~\ref{sec:mainexperiments},
 $R$-SafeBayes and $I$-SafeBayes for models with varying variance,
denoted {\em $R$-log-SafeBayes and $I$-log-SafeBayes}. Below we give explicit formulas for
each when conditioned on a fixed model $\cM_{\nc}$; the case with a
posterior on $\nc$ itself can easily be derived from these.
\paragraph{Fixed $\sigma^2$: $R$-square- and $I$-square-SafeBayes} When
conditioned on a fixed $\nc$ and $\sigma^2$ (a situation with which we
experiment in Section~\ref{sec:ridgefixedsigma}), SafeBayes tries to
minimize the $R$-log-loss, which, as an easy calculation shows, is just  the sum, from $i=0$ to
$n-1$, of
\begin{multline}\label{eq:fronkie}
  \Exp_{\beta \sim \Pi \mid z^i, \nc, \eta}
  \left[ -\log \dens(y_{i+1} \mid x_{i+1}, \beta, \sigma^2) \right]\\
  = \frac{1}{2} \log(2 \pi \sigma^2) + \frac{1}{2\sigma^2}(y_{i+1} -
  x_{i+1} \bar{\beta}_{i,\eta})^2 + \frac{1}{2} x_{i+1}
  \Sigma_{i,\eta} x_{i+1}^T,
\end{multline}
where $\bar{\beta}_{i,\eta}$ and $\Sigma_{i,\eta}$ are given as in and below
(\ref{eq:olvg}).  Note that $\bar{\beta}_{i,\eta}$ depends on $\eta$ but not
on $\sigma$, and note also that, since ${\bf X}_n^T {\bf X}_n$ (as in
(\ref{eq:olvg})) tends to increase linearly in $n$ and $p$, the final
term is of order $\nc/(n \eta)$. 

In the corresponding
in-model version of SafeBayes, we use the in-model-loss as given by $- \log
\dens(y_{i+1} \mid x_{i+1}, \bar{\beta}_{i,\eta}, \sigma^2)$, which is
equal to (\ref{eq:fronkie}) without the final term. Since the first
term of (\ref{eq:fronkie}) does not depend on the data, this version
of SafeBayes thus amounts to picking the $\hat{\eta}$ minimizing just
the sum of square-loss prediction errors, {\em which does not depend
  on the chosen $\sigma^2$}. It thus becomes a standard version of
`prequential model selection' as based on the square-loss, which in
turn is similar to (though having different asymptotics than)
leave-one-out cross validation based on the square-loss.  

Indeed, the fixed $\sigma^2$ versions of SafeBayes can be interpreted
in two ways: first, as we did until now, in terms of SafeBayes with
$\ell_{\theta}$ in (\ref{eq:rloss}) set to the log-loss, i.e.\ as a
tool for dealing with misspecification. Second, with $\ell_{\theta}$
in (\ref{eq:rloss}) set proportionally to the square-loss, as a
generic tool to learn good square-loss predictors (not distributions)
in a pseudo-Bayesian way.
More precisely, $I$-SafeBayes with the log-loss for fixed
$\sigma^2$ is equivalent to the version of $I$-SafeBayes we would get
if we set $\ell_{\beta,\sigma^2}(x,y) := C (y- x \beta)^2$, for any
constant $C > 0$. Similarly, $R$-SafeBayes with the log-loss for
fixed $\sigma^2$ is equivalent to the version of $R$-SafeBayes we
would get if we set $\ell_{\beta,\sigma^2}(x,y) := C(y- x \beta)^2$,
although now equivalence only holds if we set $C = 1/2 \sigma^2$. For this reason we will now refer to them as {\em $I$-square-SafeBayes and $R$-square-SafeBayes}, respectively.

\paragraph{Varying $\sigma^2$: $R$-log- and $I$-log-SafeBayes} Next
consider the situation with fixed $p$ and varying $\sigma^2$, with
posterior on $\sigma^2$ an inverse gamma distribution with parameters $a_{n,\eta}$
and $b_{n,\eta}$ as given by (\ref{eq:anbn}). Then the $R$-log-loss is
given by
\begin{multline}\label{eq:nodig}
\Exp_{\sigma^2, \beta \sim \Pi \mid z^i, \nc, \eta}
\left[ -\log \dens(y_{i+1} \mid x_{i+1}, \beta, \sigma^2) \right]\\
= \frac{1}{2} \log 2 \pi b_{i,\eta}
- \frac{1}{2} \psi(a_{i,\eta})
+ \frac{1}{2} \frac{(y_{i+1} - x_{i+1} \bar{\beta}_{i,\eta})^2}{b_{i,\eta} / a_{i,\eta}}
+ \frac{1}{2} x_{i+1} \Sigma_{i,\eta} x_{i+1}^T  \\
= \frac{1}{2} \log 2 \pi \bar{\sigma}^2_{i,\eta} +  \frac{1}{2} \frac{(y_{i+1} - x_{i+1} \bar{\beta}_{i,\eta})^2}{\bar{\sigma}^2_{i,\eta}} + \frac{1}{2} x_{i+1} \Sigma_{i,\eta} x_{i+1}^T 
+r(i,\eta),
\end{multline}
where $\psi$ is the digamma function, $\bar{\sigma}^2_{i,\eta}$ is the
$\eta$-posterior expectation of $\sigma^2$ as given by
(\ref{eq:expsigma}) and $r(i,\eta)$ is a remainder function which is
$O(1/i)$ whenever $\sum_{i=1}^n (y_i - x_i \beta_{n,\eta})^2$
increases linearly in $i$. This final approximation follows by
(\ref{eq:expsigma}) and because we have $\psi(x) \in [\log(x-1),\log
x]$. $R$-SafeBayes for varying $\sigma^2$ minimizes (\ref{eq:nodig}),
and, because there is now only a log-loss and not a direct square-loss
interpretation, we will call it {\em $R$-log-SafeBayes\/} from now on.

To calculate the corresponding in-model version of SafeBayes, {\em $I$-log-SafeBayes},  note that it minimizes the sum of 
\begin{equation}\label{eq:bluenote}
  - \log f(y_{i+1} \mid x_{i+1}, \bar{\beta}_{i,\eta}, \bar{\sigma}^2_{i,\eta}) = 
  \frac{1}{2} \log 2 \pi \bar{\sigma}^2_{i,\eta} +  \frac{1}{2} \frac{(y_{i+1} - x_{i+1} \beta_{i,\eta})^2}{\bar{\sigma}^2_{i,\eta}}.
\end{equation}
Comparing the four versions of SafeBayes, we see that the
both $R$-SafeBayeses have an additional term which decreases in $\eta$,
increases in model dimensionality $\nc$ (via the size of the matrix
$\Sigma_{i,\eta}$) but becomes negligible for $n \gg \nc$.

\subsection{SafeBayes learns to predict as well as the Optimal Distribution}
We first define the {\em Ces\`aro-averaged\/} posterior given data
$Z^n$ by setting, for any subset $\Theta' \subset\Theta$,
\begin{equation}\label{eq:cesaro}
\Pi_{\text{\sc Ces}}(\Theta' \mid Z^n,\eta) := \frac{1}{n}\sum_{i=1}^n \Pi(\Theta' \mid Z^ i, \eta)
\end{equation}
to be the posterior probability of $\Theta'$ averaged over the $n$
posterior distributions obtained so far. Predicting based on
Ces\`aro-averaged posteriors was introduced independently by several
authors \citep{barron1987bayes,helmbold1992some,Yang00regression,Catoni97} and
has received a lot of attention in the machine learning literature in
recent years, also under the name ``on-line to batch conversion of
Bayes'' or {\em progressive mixture rule\/}
\citep{audibert2007progressive} or {\em mirror averaging\/}
\citep{juditsky2008learning,DalalyanT12}, but is of course unnatural
from a Bayesian perspective.

The main result of \cite{Grunwald12} essentially states the following:
suppose that, under $P^*$, the density ratios are uniformly bounded,
i.e.\ there is a finite $v$ such that for all
$\theta, \theta' \in \Theta$,
$P^*(\dens_{\theta}(Y \mid X)/\dens_{\theta'}(Y \mid X) \leq v) =
1$. Suppose further that the prior $\Pi$ assigns `sufficient mass' in
KL-neighborhoods of $P_{\tilde\theta}$. Then $\Pi_{\text{\sc Ces}}$
applied with the $\hat{\eta}$ learned by the Safe Bayesian algorithm
concentrates on the optimal $P_{\tilde{\theta}}$. That is, let
$\Theta_{\delta}$ be the subset of all $\theta \in \Theta$ with
$D(P^* \| P_{\theta}) \geq D(P^* \| P_{\tilde{\theta}}) +
\delta$. Then for all $\delta > 0$, with $P^*$-probability 1, as
$n \rightarrow \infty$, we have that
$\Pi_{\text{\sc Ces}}(\Theta_{\delta} \mid Z^n, \hat{\eta})$ goes to
$0$. \cite{Grunwald12} goes on to show that in several settings, one
can design priors such that the rate at which the posterior
concentrates is minimax optimal, i.e.\ no algorithm can do better in
general. On the negative side, the requirement of bounded density
ratio is strong, and the replacement of the standard posterior by the
Ces\`aro one is awkward. On the positive side, the theorem has no
further conditions and can be applied to parametric and nonparametric
cases alike. It is very likely that the requirement of bounded density
ratios can be dropped, cf. the developments of
\cite{grunwald2016fast}; it is not so easy to drop the
C\`esaro-replacement, but we suspect that this is an artifact of the
proof technique. To see whether there is any practical difference,
below we include experimental results both for the Ces\`aro-averaged
$\eta$-generalized posterior
$\Pi_{\text{\sc Ces}}(\cdot \mid Z^n, \hat{\eta})$ and for the
standard $\eta$-generalized posterior
$\Pi(\cdot \mid Z^n, \hat{\eta})$.

\pagebreak

\section{Main Experiment: Varying $\sigma^2$}
\label{sec:mainexperiments}
In this section we provide our main experimental results, based on
linear models $\cM_\nc$ as defined in Section~\ref{sec:linear} with a
prior on both the mean and the variance.
Figure~\ref{fig:mainexperimenta}--\ref{fig:mainexperimentd} depict,
and Section~\ref{sec:modsel} discusses the results of model
selection/averaging experiments, which choose/average between the
models $0, \ldots, \pmax$, where we consider first an incorrectly and
then a correctly specified model, both with $\pmax = 50$ and later
with $\pmax = 100$. Section~\ref{sec:ridge} contains and interprets additional experiments
on Bayesian ridge regression, with a fixed $p$; a multitude of
additional experiments is provided in the
appendices. Section~\ref{sec:summary} in this section summarizes the
relevant findings of these additional experiments. 
\subsection{Preparing Main  Experiments: Model, Priors, Method, `Truth'}
\label{sec:preparing}
In this subsection we prepare the experiments: Section~\ref{sec:priors}
describes our priors $\pi$; Section~\ref{sec:truth} concerns the
sampling (`true') distributions $P^*$ with which we experiment; and
finally, Section~\ref{sec:statisticsshown} describes the data
statistics that we will report.
\subsubsection{The Priors}
\label{sec:priors}
\paragraph{Prior on Models}
In our model selection/averaging experiments, we use a fat-tailed
prior on the models given by
\[
\pi(\nc) \propto \frac{1}{ (\nc+2) (\log(\nc+2))^2}.
\]
This prior was chosen because it remains well-defined for an infinite
collection of models, even though we only use finitely many in our
experiments. 

{\em Variation}. As a sanity check we did repeat some of our experiments
with a uniform prior on $0,\ldots, \pmax$ instead; the results were
indistinguishable.

\paragraph{Prior on Parameters given Models}

Each model $\cM_{\nc}$ has parameters $\beta, \sigma^2$, on which we
put the standard conjugate priors as described in
Section~\ref{sec:instantiation}. We set the mean of the prior on
$\beta$ to $\betao = \mathbf 0$, and its covariance matrix to
$\sigma^2 \Sigma_0$. Our main experiments below are based on an {\em
  informative\/} instantiation of $\Sigma_0$, using the identity
matrix $\Sigma_0 = {\mathbf I}_{\nc+1}$; this prior equals the
posterior we would get by starting with an improper Jeffreys' prior on
$\beta$ and then observing, for each coefficient $\beta_j$, one extra
point $z = (x, 0)$ with $x_j = 1$ and $x_i = 0$ for $i \neq j$.

{\em Variations} We also ran experiments with a `slightly informative'
$\Sigma_0$, where we set $\Sigma_0 = 1000 \cdot {\mathbf I}_{\nc+1}$,
comparable to observing points $z = (x, 0)$ with $x_j =
1/\sqrt{1000}$. Finally, following the standard reference
\citet{raftery1997bayesian}, we also used a prior with a level of
informativeness depending on the submodel, described in more
detail in Appendix~\ref{sec:priorvariations}. 

As to the prior on $\sigma^2$: Jeffreys' prior is obtained for the
choice $a_0 = b_0 = 0$ in (\ref{eq:inversegamma}). We do not use this
improper prior, because of the well-known issues with Bayes factors
under improper priors \citep{OHagan95}. Moreover, to calculate the
posterior's reliability (defined in Section~\ref{sec:statisticsshown}
and shown in Figure~\ref{fig:mainexperimenta}) and also for the
$I$-log-loss, we need to calculate the posterior
expectation of the variance $\sigma^2$ quantity as given by
(\ref{eq:expsigma}), which is only well-defined and
finite for $a_n > 1$. We want to make $\pi(\sigma^2)$ as uninformative
as possible while ensuring that (for any positive learning rate) this
variance exists for the posterior based on at least one sample. This
is accomplished by choosing $a_0 = 1$: for standard Bayes,
the posterior after one observation has $a_1 = a_0 + 1/2$; for
generalized Bayes, $a_1 = a_0 + \eta/2$. To
set $b_0$, we use that $b_0 / a_0$ represents the sample variance of a virtual initial data sequence \citep[Section 14.8]{BDA3}.
We choose $b_0 = 1/40$
so that $b_0 / a_0 = 1/40$, the true variance of the noise in our
data, as we describe next.

\subsubsection{The ``Truth'' (Sampling Distribution)}
\label{sec:truth}
Our experiments fall into two categories: correct-model and
wrong-model experiments.
\paragraph{Correct-Model Experiments}
Here $X_1, X_2, \ldots$ are sampled i.i.d., with, for
each individual $X_i = (X_{i1}, \ldots, X_{i \pmax})$, $X_{i1} \ldots,
X_{i \pmax}$ i.i.d.~$\sim N(0,1)$.  Given each $X_i$, $Y_i$ is
generated as
\begin{equation}\label{eq:pseudotruereg}
Y_i = .1 \cdot (X_{i1} + \ldots + X_{i4}) + \epsilon_i, 
\end{equation}
where the $\epsilon_i$ are i.i.d.~$\sim N(0,\sigma^{*2})$ with 
variance $\sigma^{*2} = 1/40$. 
\paragraph{Wrong-Model Experiments}
Now at each time point $i$, a fair coin is tossed independently of
everything else. If the coin lands heads, then the point is `easy',
and $(X_i,Y_i) := ({\bf 0},0)$. If the coin lands tails, then $X_i$ is
generated as for the correct model, and $Y_i$ is generated as
(\ref{eq:pseudotruereg}), but now the noise random variables have
variance $\sigma_0^2 = 2\sigma^{*2} = 1/20$. Thus, $Z_i = (X_i,Y_i)$
is generated as in the true model case but with a larger variance;
this larger variance has been chosen so that the marginal variance of
each $Y_i$ is the same value $\sigma^{*2}$ in both experiments.

From the results in Section~\ref{sec:optimalb} we immediately see
that, for both experiments, the optimal model is $\cM_{\tilde{\nc}}$
for $\tilde{\nc} = 4$, and the optimal distribution in $\cM$ and
$\cM_{\tilde{\nc}}$ is parameterized by $\tilde\theta =
(\tilde{\nc},\tilde{\beta},\tilde{\sigma}^2)$ with $\tilde{\nc}=4$, $\tilde{\beta} = (\tilde\beta_0, \ldots \tilde\beta_4) =
(0,.1,.1,.1,.1), \tilde\sigma^2 = 1/40$ (in the correct model experiment, $\tilde\sigma^2 = \sigma^{*2}$; in the wrong model experiment, since $\tilde{\sigma}^2$ must be reliable, it must be equal to the square-risk obtained with $(\tilde{p},\tilde\beta)$, which is $(1/2) \cdot (1/20) = 1/40$).  $f(x) := x \tilde\beta$ is then  equal to the {\em true\/} regression function $\Exp_{P^*}[Y \mid X]$.

{\em Variations}.  We have already seen a variation of Experiments 1
and 2 depicted in Figure~\ref{fig:polynomial} and
\ref{fig:polynomialrisk}. In the correct model version of that
experiment, $P^*$ is defined by setting $X_{j} = S^j$, and let $S$ be
uniformly distributed on $[-1,1]$ and set $Y = 0 + \epsilon$, where
$\epsilon \sim N(0,\sigma^{*2})$, with $\sigma^{*2} = 1/40$;
$(X_1,Y_1), \ldots$ are then sampled as i.i.d.~copies of $(X,Y)$. Note
that the true regression function is $0$ here. In
Appendix~\ref{sec:wildvariations} we briefly consider this and several
other variations of these ground truths.

\subsection{The Statistics We Report}
\label{sec:statisticsshown}
Figure~\ref{fig:mainexperimenta} reports the results of the
wrong-model, $\nc = 50$ experiment; Figure~\ref{fig:mainexperimentb}
shows correct-model, $\nc=50$; Figure~\ref{fig:mainexperimentc} is
about wrong-model, $\nc=100$ and Figure~\ref{fig:mainexperimentd}
depicts the correct-model, $\nc= 100$ setting.  For all four
experiments we measure three aspects of the performance of Bayes and
SafeBayes, each summarized in a separate graph. First, we show the
behavior of several prediction methods based on Safe Bayes relative to
square-risk; second, we measure whether the methods provide a good
assessment of their own predictive capabilities in terms of
square-loss, i.e.\ whether they are reliable and not `overconfident'. Third, we check a
form of model identification consistency. Below we explain these three
performance measures in detail. We postpone all experiments with
log-loss rather than square-loss to Section~\ref{sec:mixability}. We
also provide a fourth graph in each case indicating what $\hat\eta$'s
are typically selected by the two versions of SafeBayes.
\paragraph{Square-Risk}\label{page:regressionbasedonw} For a given distribution $W$ on
$(\nc,\beta,\sigma^2)$, the {\em regression function based on $W$}, a
function mapping covariate $X$ to $\reals$, abbreviated to $\Exp_W[Y
\mid X]$, is defined as
\begin{equation}\label{eq:regfun}
\Exp_W[Y \mid X] := \Exp_{(\nc,\beta,\sigma)  \sim W} \Exp_{Y \sim
  P_{\nc,\beta,\sigma} \mid X}[Y] = \Exp_{(\nc,\beta,\sigma)  \sim W}
\left[  \sum_{j=0}^{\nc} \beta_j X_j  \right].
\end{equation}
If we take $W$ to be the $\eta$-generalized posterior, then
(\ref{eq:regfun}) is also simply called the $\eta$-posterior
regression function. The {\em square-risk\/} relative to $P^*$ based
on predicting by $W$ is then defined as an extension of
(\ref{eq:squarederrorrisk}) as
\begin{equation}\label{eq:kabeljauw}
\squaredrisk(W) := \Exp_{(X,Y) \sim P^*}(Y- \Exp_W[Y \mid X] )^2.
\end{equation}
In the experiments below we measure the square-risk relative to
$P^*$ at sample size $i-1$ achieved by, respectively, (1), the
$\eta$-generalized posterior, (2), the $\eta$-generalized posterior
conditioned on the MAP (maximum a posteriori) model, and, (3), the
$\eta$-generalized Ces\`aro-averaged posteriors, i.e.
\begin{multline}
\Exp_{Z^{i-1}  \sim P^*} [\squaredrisk(W)], \text{\ with\ } \\
W= \Pi \mid Z^{i-1}, \eta \ ; \ W = \Pi \mid Z^{i-1}, \eta, \breve{\nc}_{\text{map}(Z^{i-1}, \eta)} \ ; \ W= \Pi_{\text{\sc Ces}} \mid Z^{i-1}, \eta,
\end{multline}
respectively, where the MAP (maximum a posteriori) model
$\breve{\nc}_{\text{map}(Z^{i-1}, \eta)}$ is defined as the $\nc$
achieving $\max_{\nc \in 0..\pmax} \pi(\nc \mid Z^{i-1}, \eta)$, with
$\pi(\nc \mid Z^{i-1}, \eta)$ defined as in (\ref{eq:modelposteriorb}),
and $\Pi_{\text{\sc Ces}}$ is the Ces\`aro-averaged posterior as
defined as in (\ref{eq:cesaro}). We do this for three values of
$\eta$: (a) $\eta = 1$, corresponding to the standard Bayesian
posterior, (b), $\eta := \hat{\eta}(Z^{i-1})$ set by the $R$-log Safe
Bayesian algorithm run on the past data $Z^{i-1}$, and (c) $\eta$ set
by the $I$-log Safe Bayesian algorithm. In the figures of
Section~\ref{sec:modsel}, 1(a) is abbreviated to {\em Bayes}, 1(b) is
{\em R-log-SafeBayes}, 1(c) is {\em I-log-SafeBayes}, 2(a) is {\em
  Bayes MAP}, 2(b) is {\em R-log-SafeBayes MAP}, 2(c) is {\em
  I-log-SafeBayes MAP}, and results with Ces\`aro-averaging are
discussed but not explicitly shown. In Section~\ref{sec:ridge},
additionally 3(a) is {\em Bayes Ces\`aro}, 3(b) is {\em R-log-SafeBayes
  Ces\`aro}, and 3(c) is {\em I-log-SafeBayes Ces\`aro}.

Concerning the three square-risks that we record: The first choice is
the most natural, corresponding to the prediction (regression
function) according to the `standard' $\eta$-generalized posterior;
the second corresponds to the situation where one first selects a
single submodel $\breve{\nc}_{\text{map}}$ and then bases all
predictions on that model; it has been included because such methods
are often adopted in practice. The
third choice, the {\em Ces\`aro-averaged generalized posterior\/} is
included because, when $\eta =\hat{\eta}$ is set by Safe Bayes, this
is the choice that \cite{Grunwald12} provides theoretical convergence
results for. But we are also interested in the results for
the Ces\`aro-average for 
$\eta = 1$, because this 
has been proposed earlier --- albeit somewhat implicitly and with
different models  --- to stabilize Bayesian
predictions in adversarial circumstances \citep{helmbold1992some}, so
we include these as well.

In Figure~\ref{fig:mainexperimenta} and subsequent figures below, we
depict these quantities by sequentially sampling data $Z_1, Z_2,
\ldots, Z_{\max}$ i.i.d.\ from a $P^*$ as defined above in
Section~\ref{sec:truth}, where $\max$ is some large number. At each
$i$, after the first $i-1$ points $Z^{i-1}$ have been sampled, we
compute the three square-risks given above. We repeat the whole
procedure a number of times (called `runs'); the graphs show the
average risks over these runs.

\paragraph{MAP-model identification/Occam's Razor} 
When the goal of inference is model identification, `consistency' of a
method is often defined as its ability to identify the smallest model
$\cM_{\tilde{p}}$ containing the `pseudo-truth'
$(\tilde\beta,\tilde\sigma^2)$. To see whether standard Bayes and/or
Safe Bayes are consistent in this sense, we check whether the MAP
model $\breve{\nc}_{\text{map}(Z^{i-1}, \eta)}$ is equal to
$\tilde{\nc}$.

\paragraph{Reliability vs.\ Overconfidence}
Does Bayes learn how good it is in terms of squared error?   To answer this question, we define, for a predictive
distribution $W$ as in (\ref{eq:kabeljauw}) above, $U_i^{[W]}$ (a function of $X_i, Y_i$ and (through $W$) of $Z^{i-1}$),  as
$$U_i^{[W]} = (Y_i - \Exp_{W}[Y_i \mid X_i])^2.$$
This is the error we make if we predict $Y_i$ using the regression
function based on prediction method $W$. In the graphs in the next
sections we plot the {\em self-confidence ratio\/} $\Exp_{X_i,Y_i \sim
  P^*}[U_i^{[W]}]/\Exp_{X_i \sim P^*}\Exp_{Y_i \sim W \mid
  X_i}[U_i^{[W]}]$ as a function of $i$ for the three prediction
methods/choices of $W$ defined above. We may think of this as the
ratio between the actual expected prediction error (measured in
square-loss) one gets by using a predictor who based predictions on
$W$ and the marginal (averaged over $X$) subjectively expected
prediction error by this predictor.  We previously, in
Section~\ref{sec:optimalb}, showed that the KL-optimal
$(\tilde{p},\tilde{\beta},\tilde{\sigma}^2)$ is {\em reliable\/}: this
means that, if we would take $W$ the point mass on
$(\tilde{p},\tilde{\beta},\tilde{\sigma}^2)$ and thus, irrespective of
past data $Z^{i-1}$, would predict by
$\Exp_{(\tilde{p},\tilde{\beta},\tilde{\sigma}^2)}[Y_i \mid X_i] =
\sum_{j=0}^{\tilde{p}} \tilde{\beta}_j X_{ij}$, then the ratio would
be $1$. For the $W$ learned from data considered above, a value larger
than $1$ indicates that $W$ does not implement a `reliable' method in
the sense of Section~\ref{sec:optimalb}, but rather overconfident: it
predicts its predictions to be better than they actually are, in terms
of square-risk.

\subsection{Main Model Selection/Averaging Experiment}
\label{sec:modsel}
\begin{figure}[htp]{\hspace*{0.2\textwidth}
\includegraphics[width=0.6\textwidth]{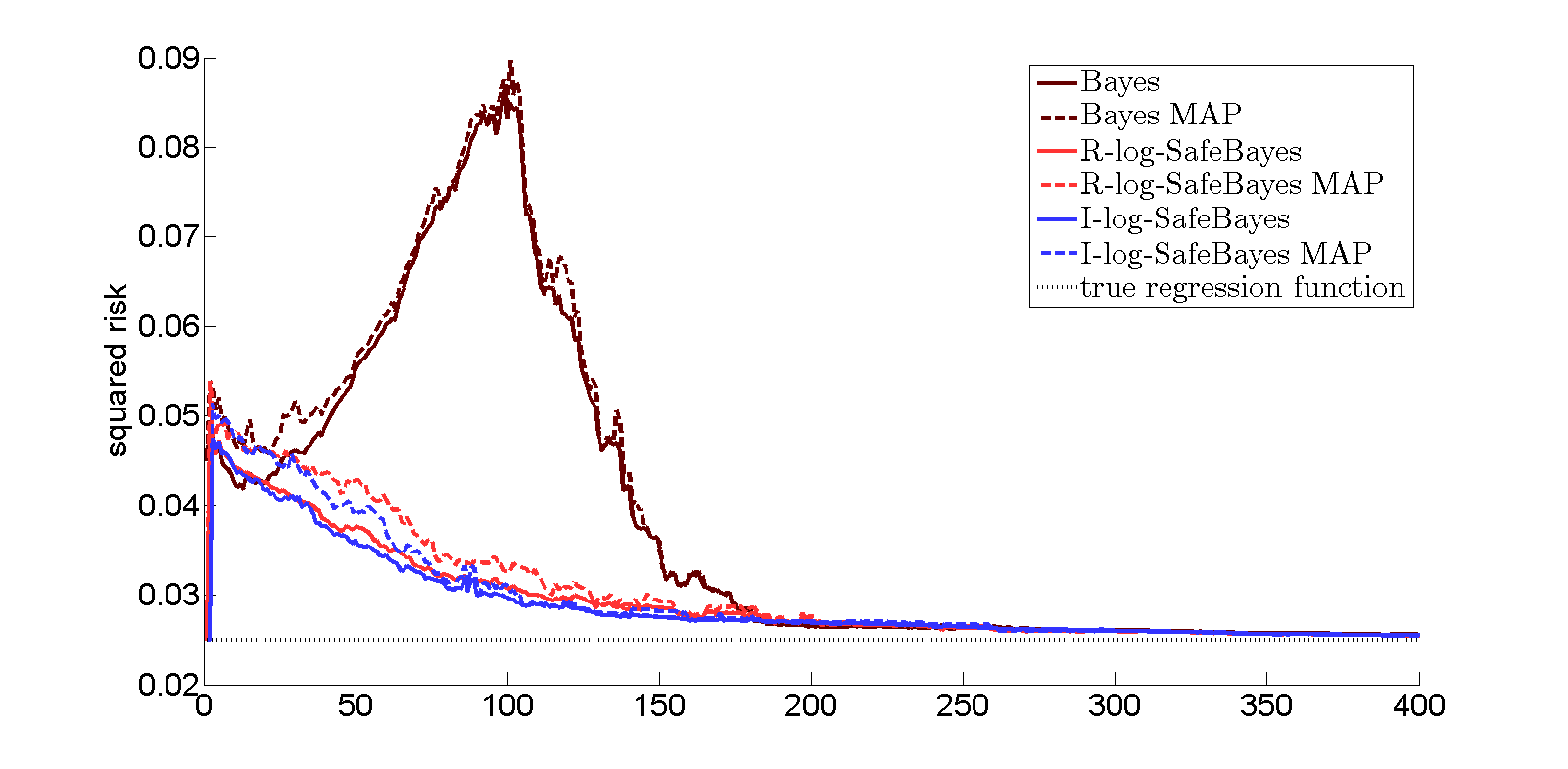} \\
\hspace*{0.2\textwidth}
\includegraphics[width=0.6\textwidth]{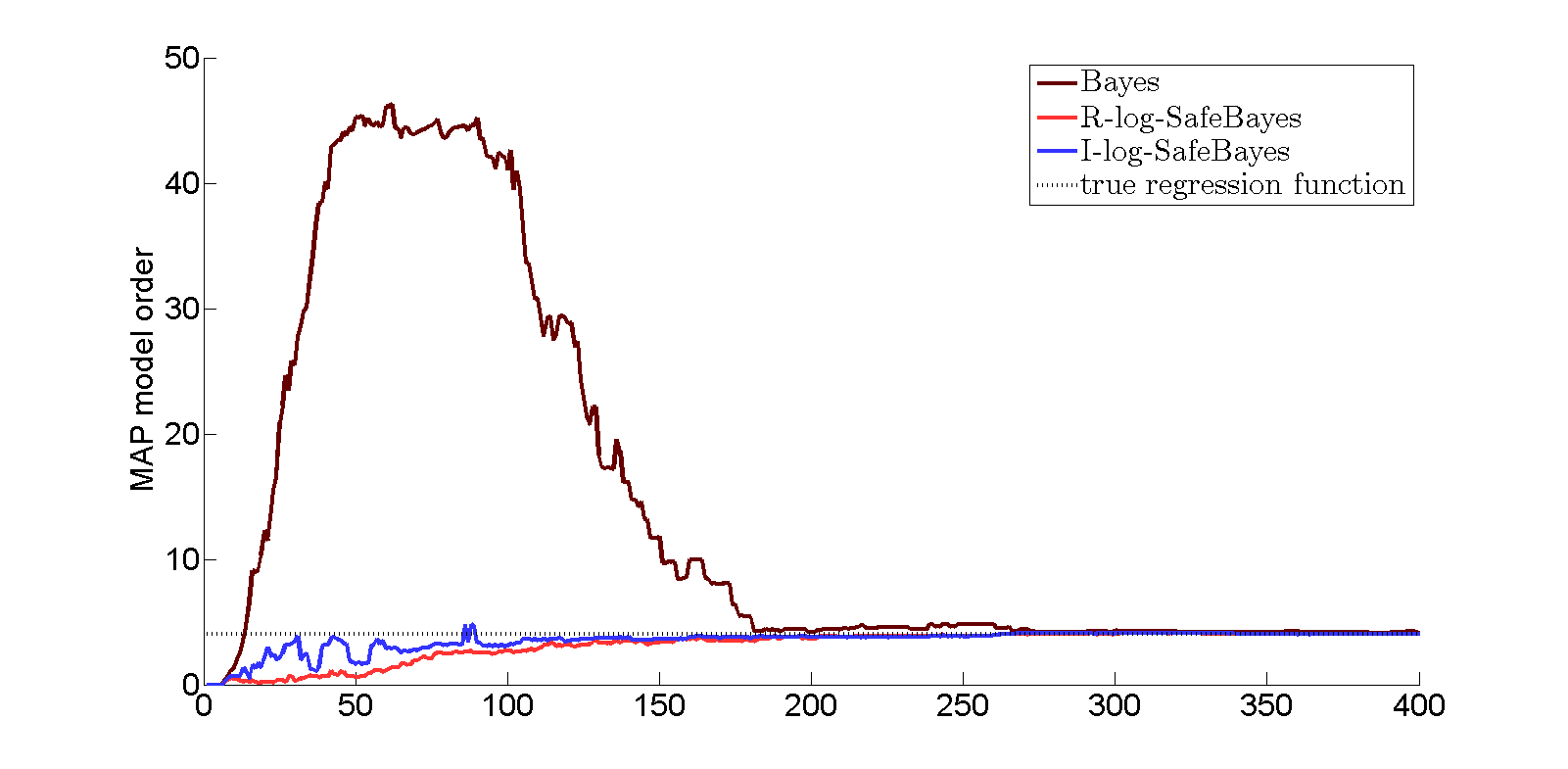} \\
\hspace*{0.2\textwidth}
\includegraphics[width=0.6\textwidth]{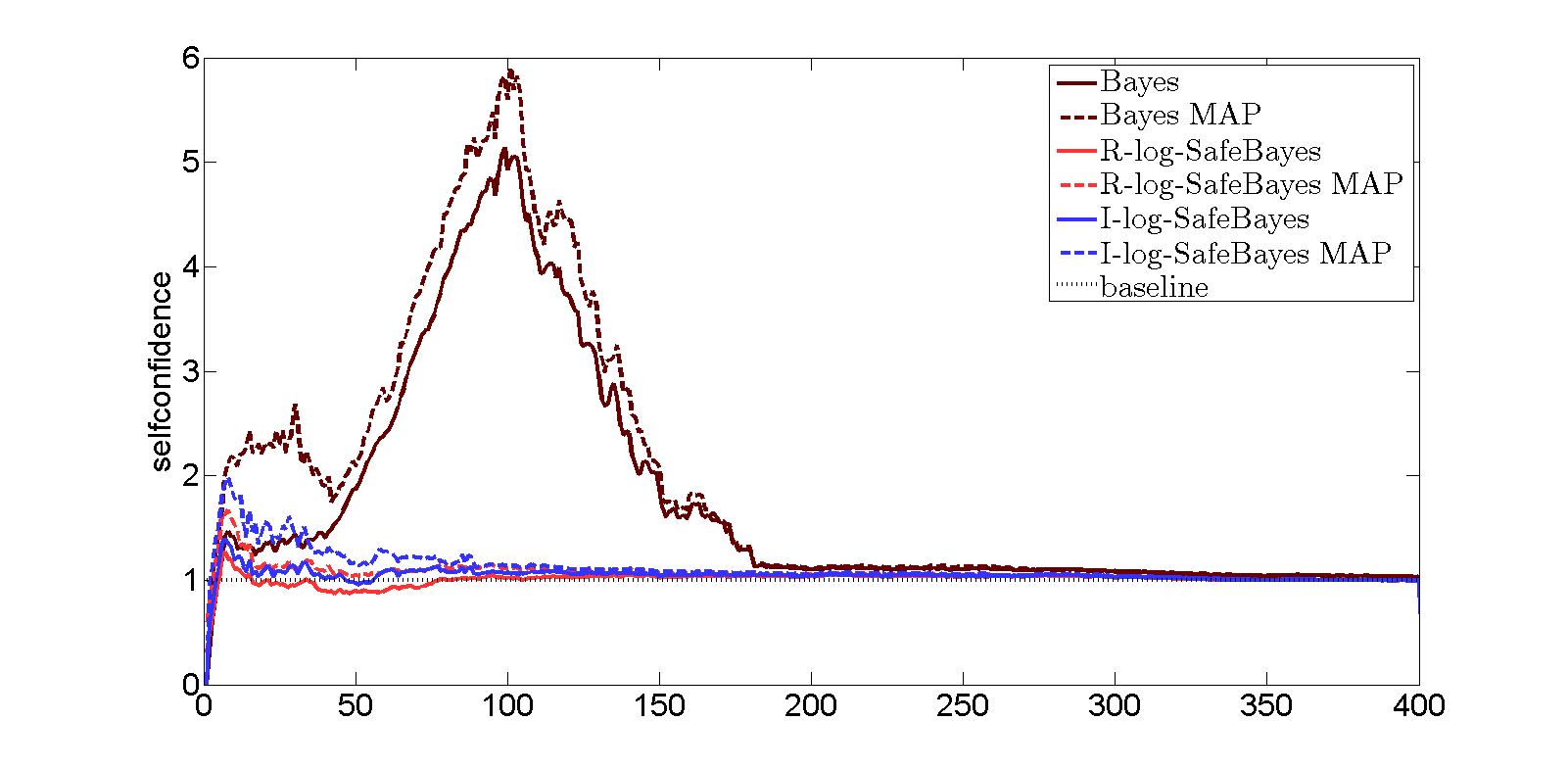} \\
\hspace*{0.2\textwidth}
\includegraphics[width=0.6\textwidth]{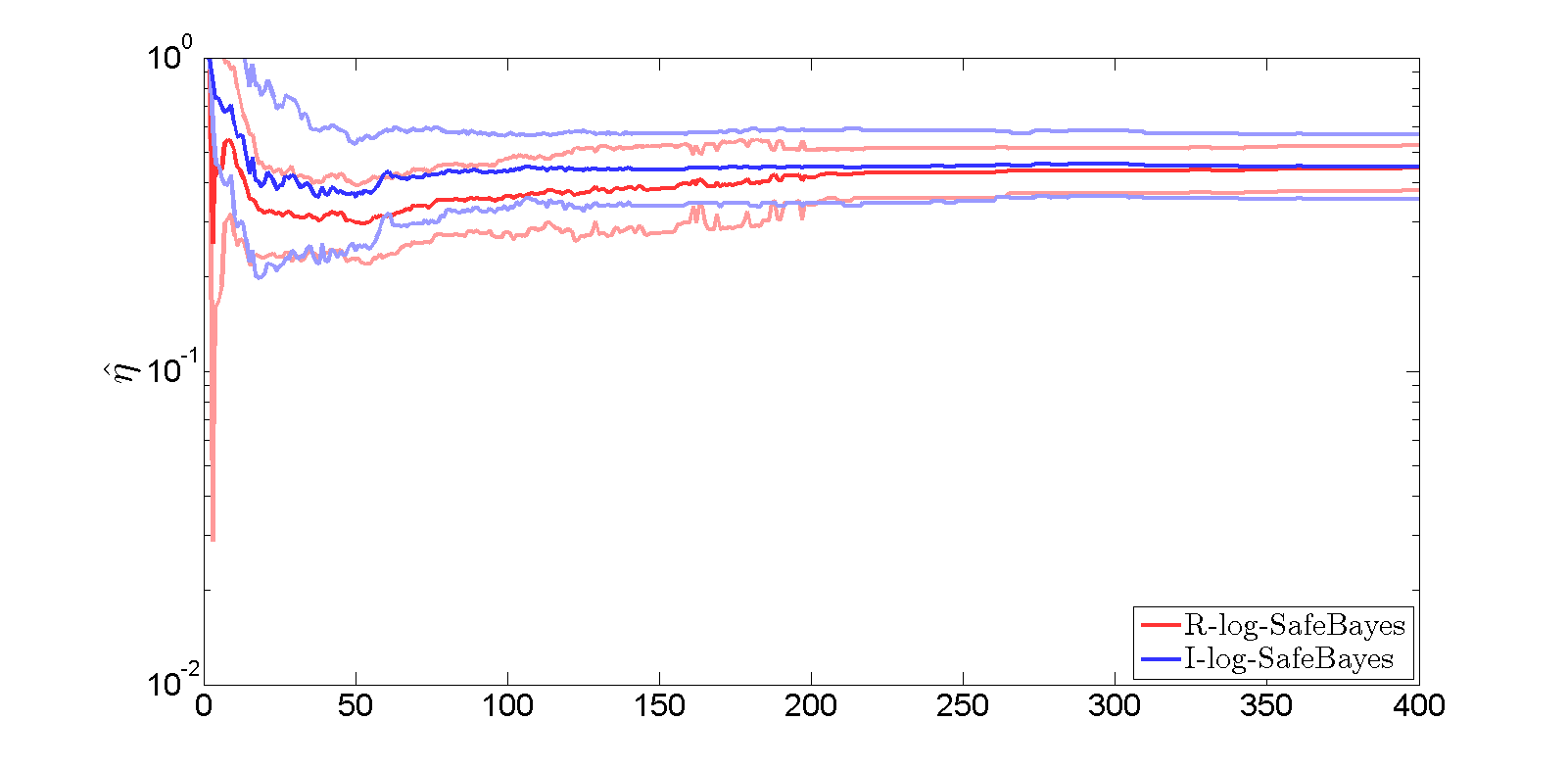}
\caption{\label{fig:mainexperimenta} Four graphs showing respectively
  the square-risk, MAP model order, overconfidence (lack of
  reliability), and selected $\hat\eta$ at each sample size, each
  averaged over 30 runs, for the wrong-model experiment with $\pmax =
  50$, for the methods indicated in Section~\ref{sec:statisticsshown}.  For the
  selected-$\hat\eta$ graph, the pale lines are one standard
  deviation apart from the average; all lines in this graph were
  computed over $\hat\eta$ indices (so that the lines depict the
  geometric mean over the values of $\hat\eta$ themselves).
}}
\end{figure}

\begin{figure}[h!]{\hspace*{0.2\textwidth}
\includegraphics[width=0.6\textwidth]{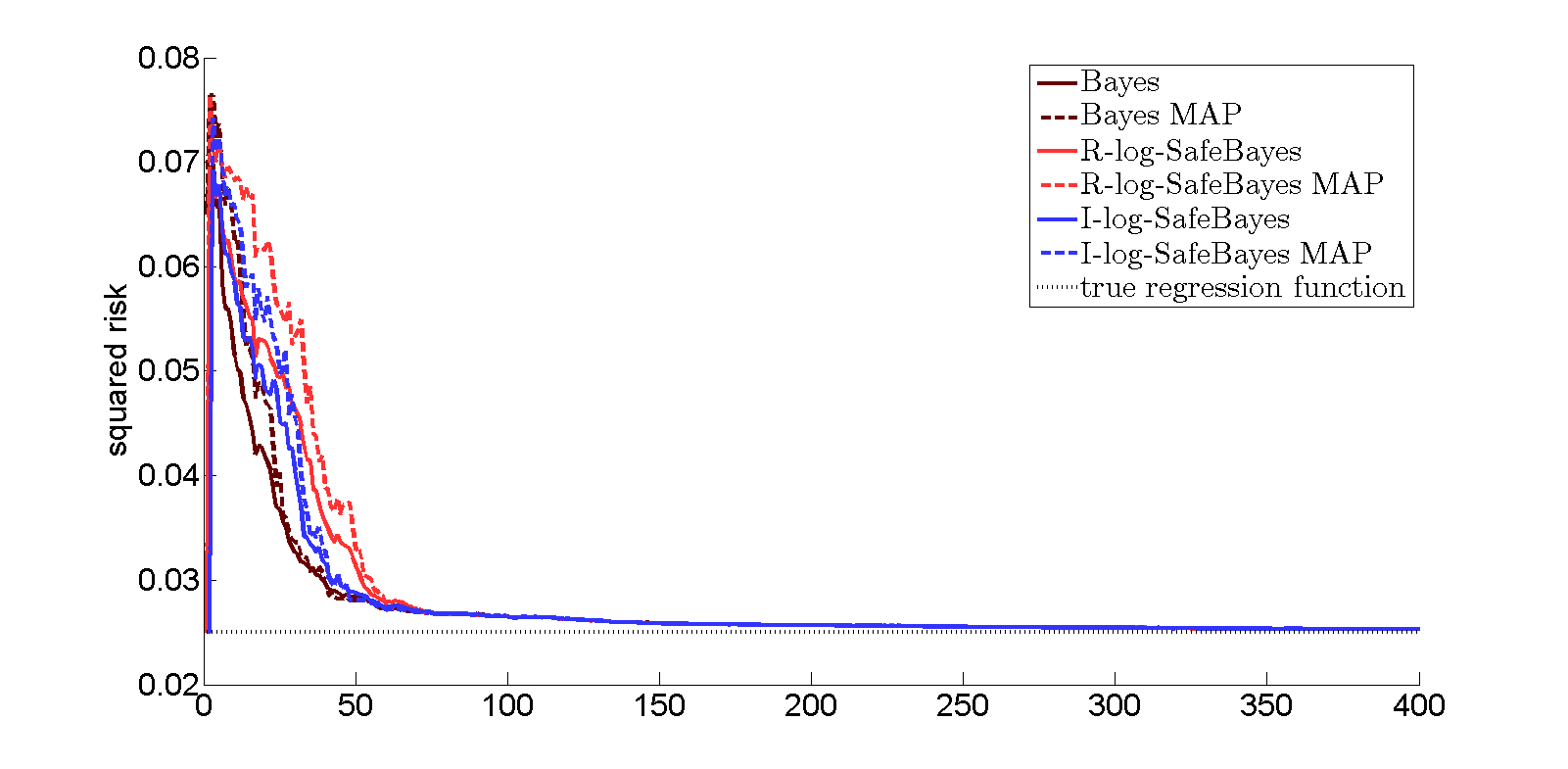} \\
\hspace*{0.2\textwidth}
\includegraphics[width=0.6\textwidth]{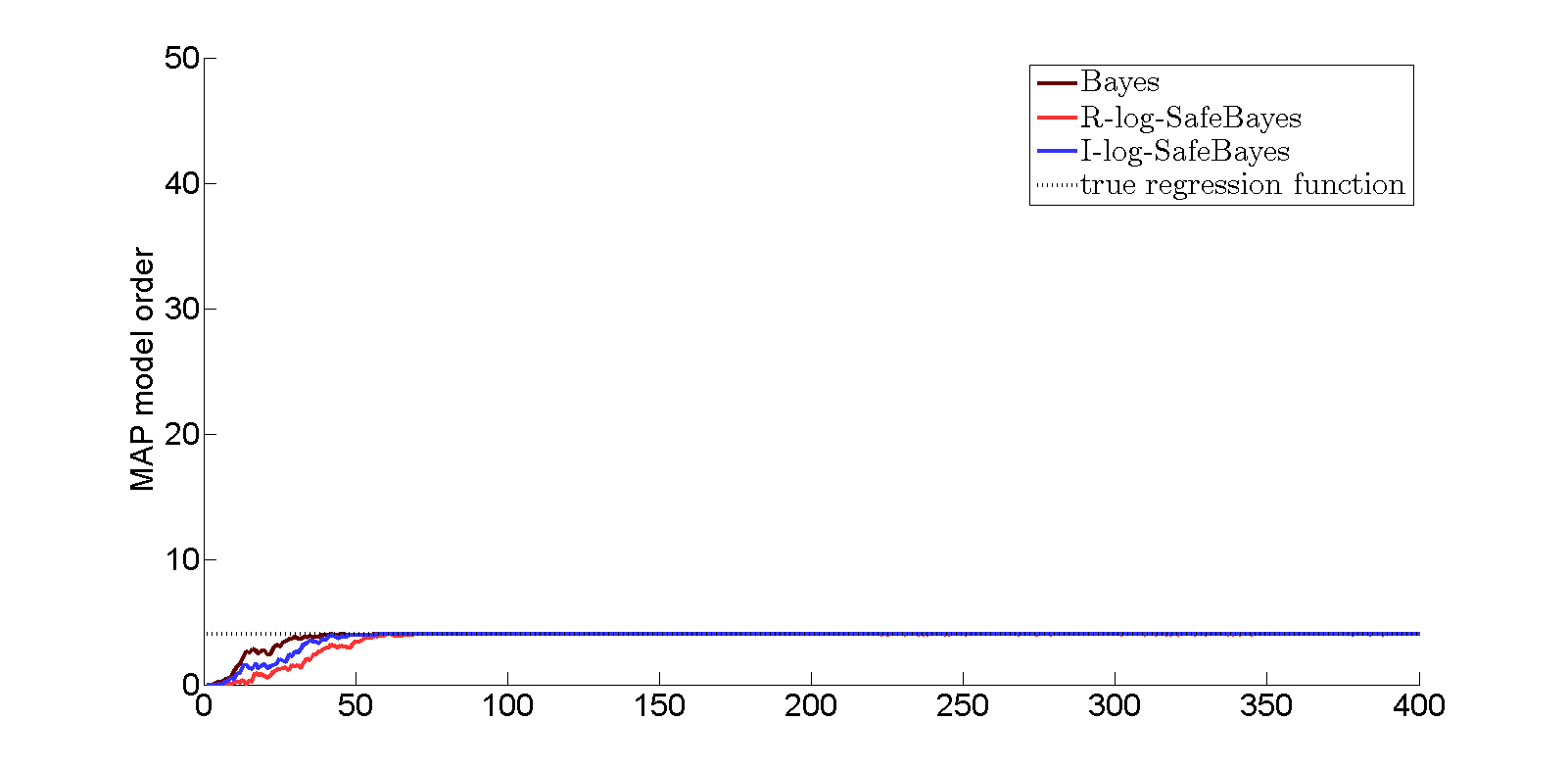} \\
\hspace*{0.2\textwidth}
\includegraphics[width=0.6\textwidth]{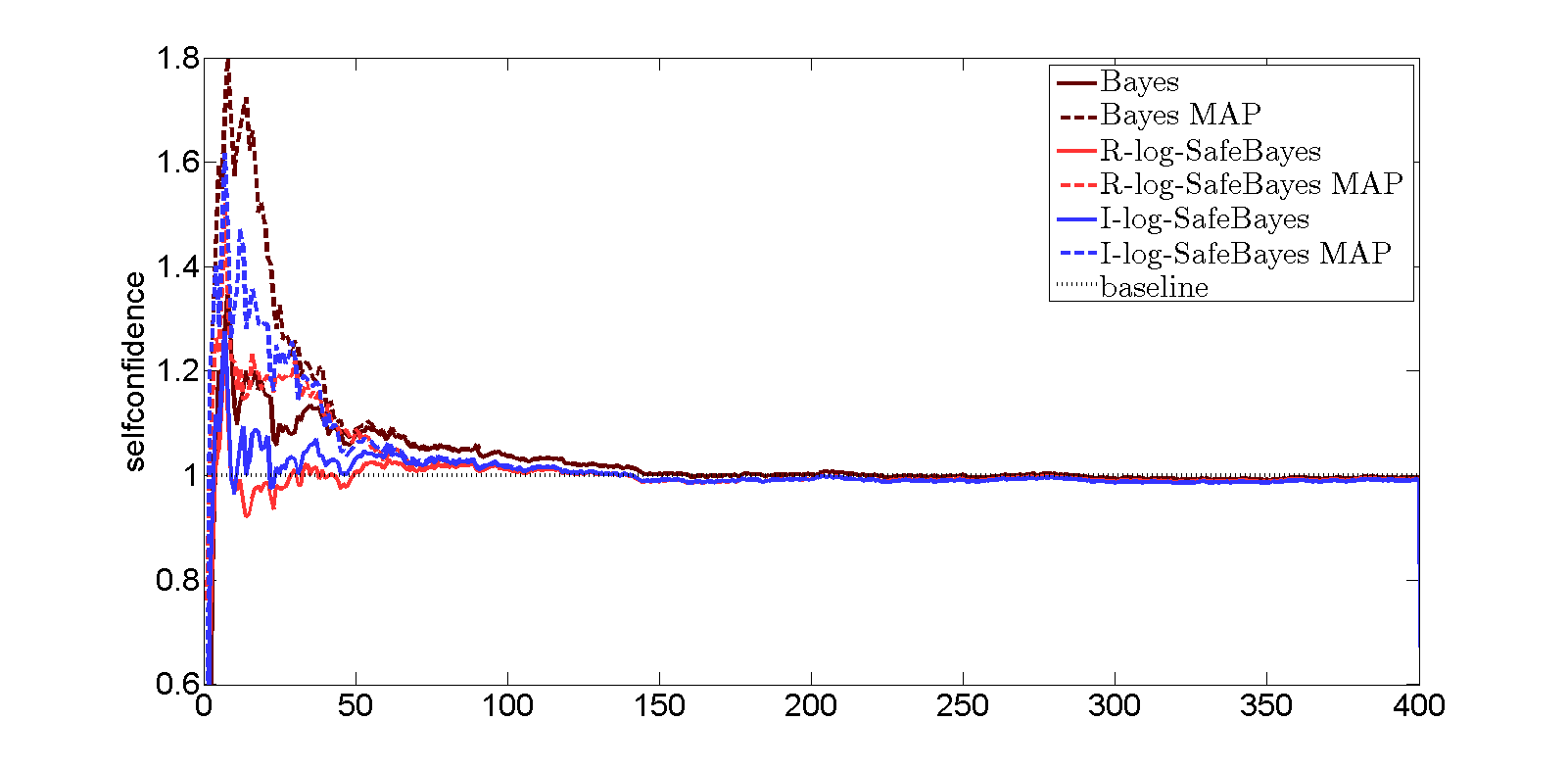} \\
\hspace*{0.2\textwidth}
\includegraphics[width=0.6\textwidth]{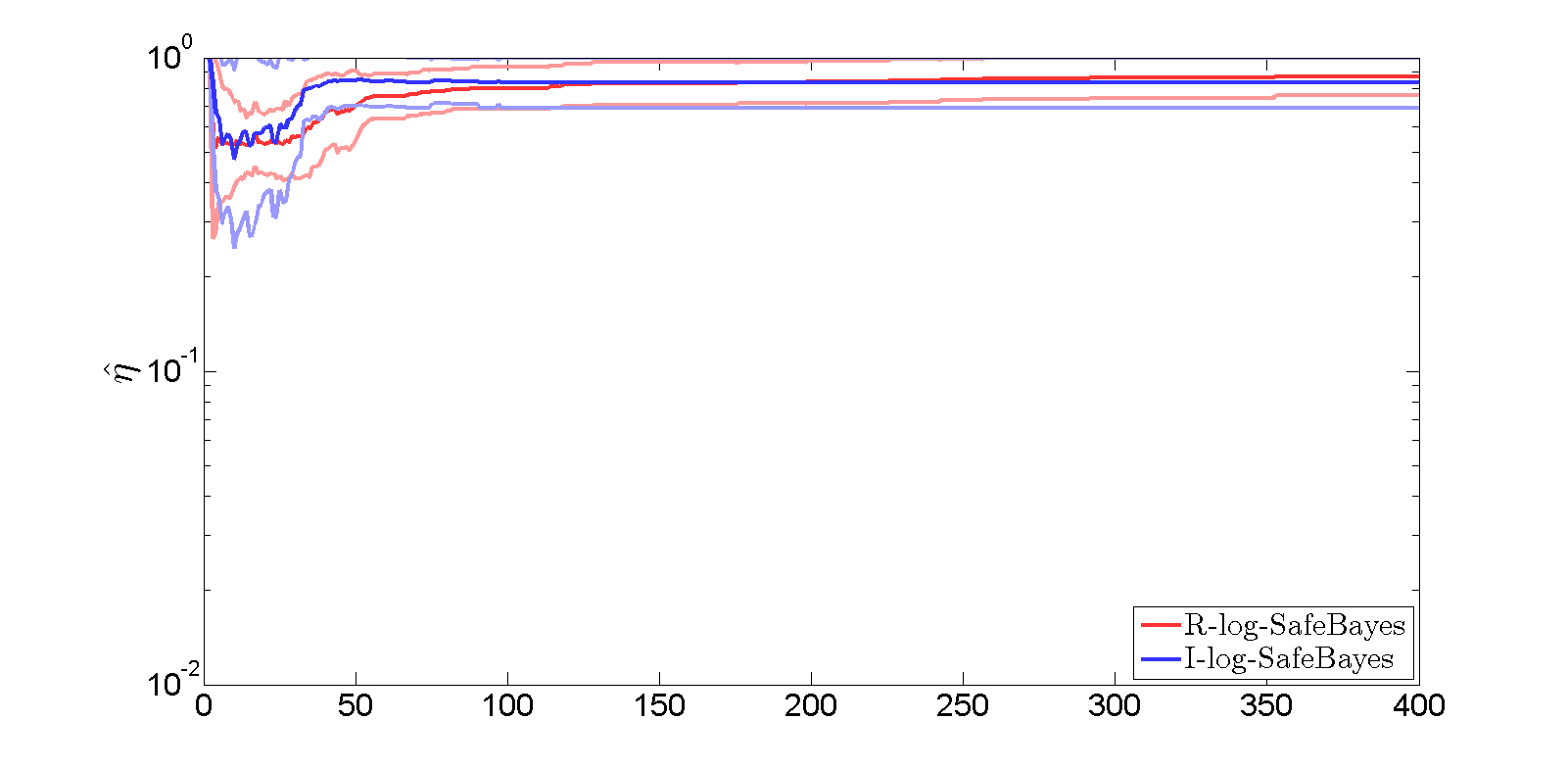}
\caption{\label{fig:mainexperimentb}
Same graphs as in Figure~\ref{fig:mainexperimenta} for the correct-model experiment with $\pmax = 50$.
}}
\end{figure}

\begin{figure}[htp]{\hspace*{0.2\textwidth}
\includegraphics[width=0.6\textwidth]{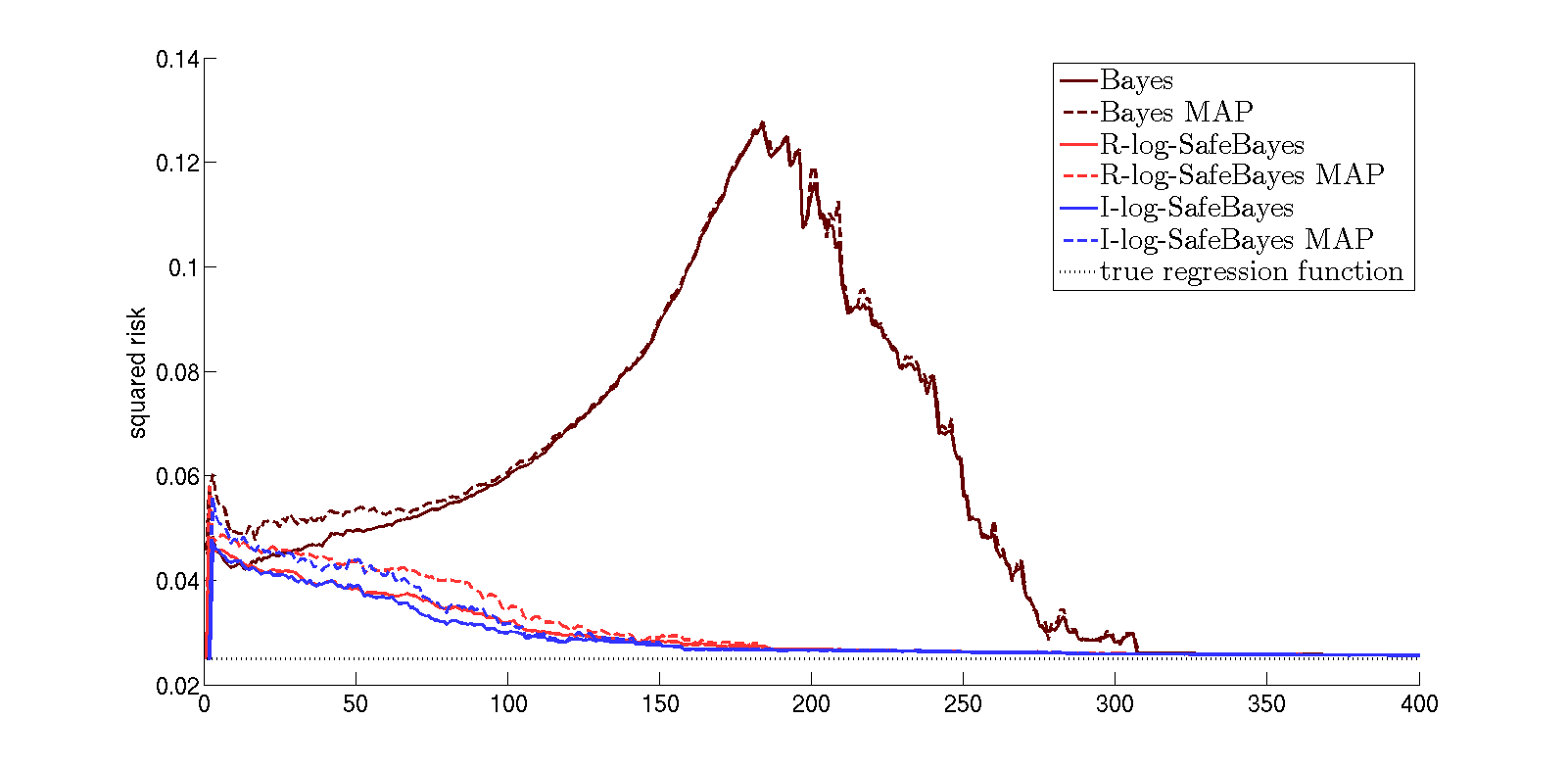} \\
\hspace*{0.2\textwidth}
\includegraphics[width=0.6\textwidth]{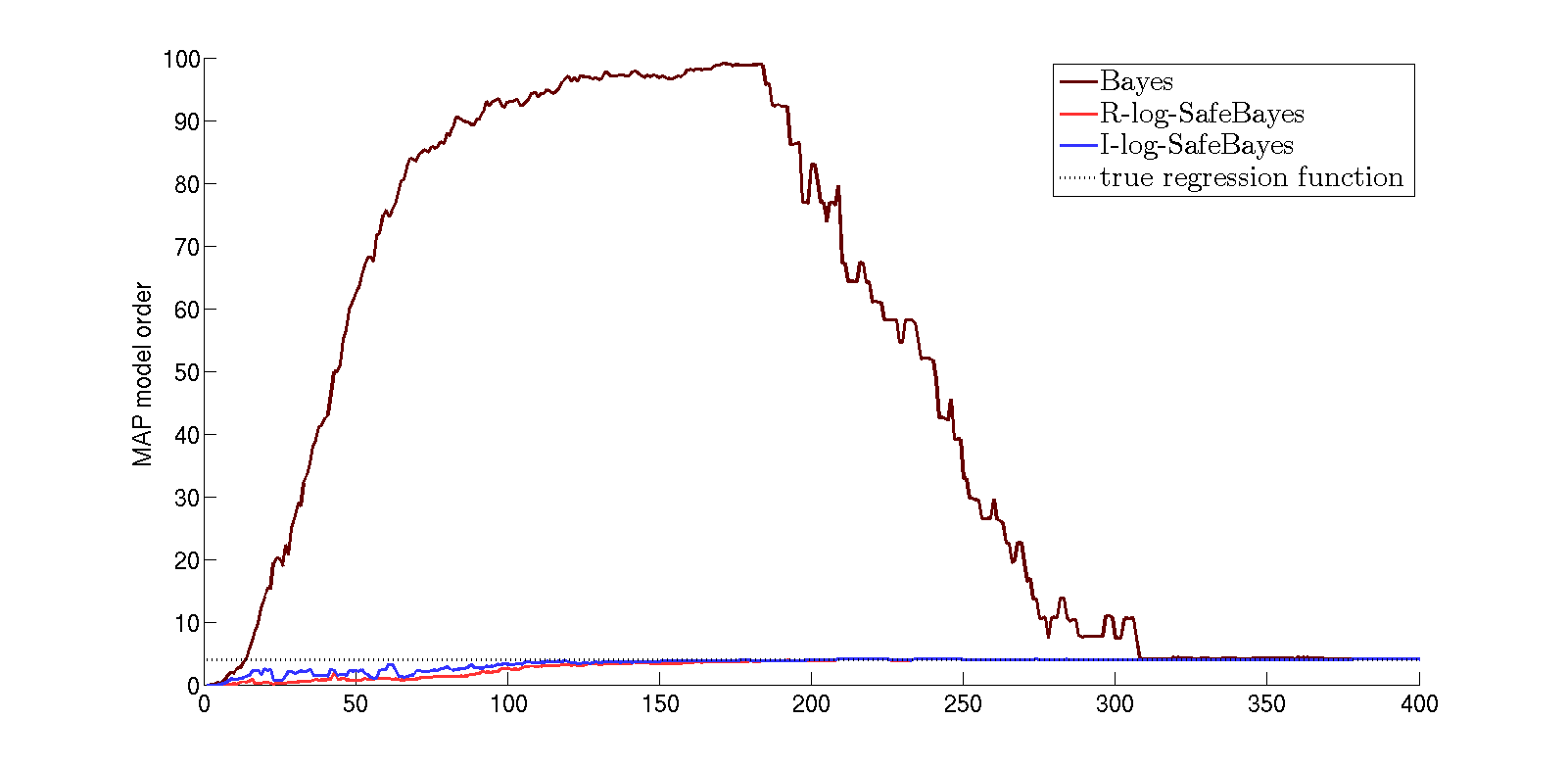} \\
\hspace*{0.2\textwidth}
\includegraphics[width=0.6\textwidth]{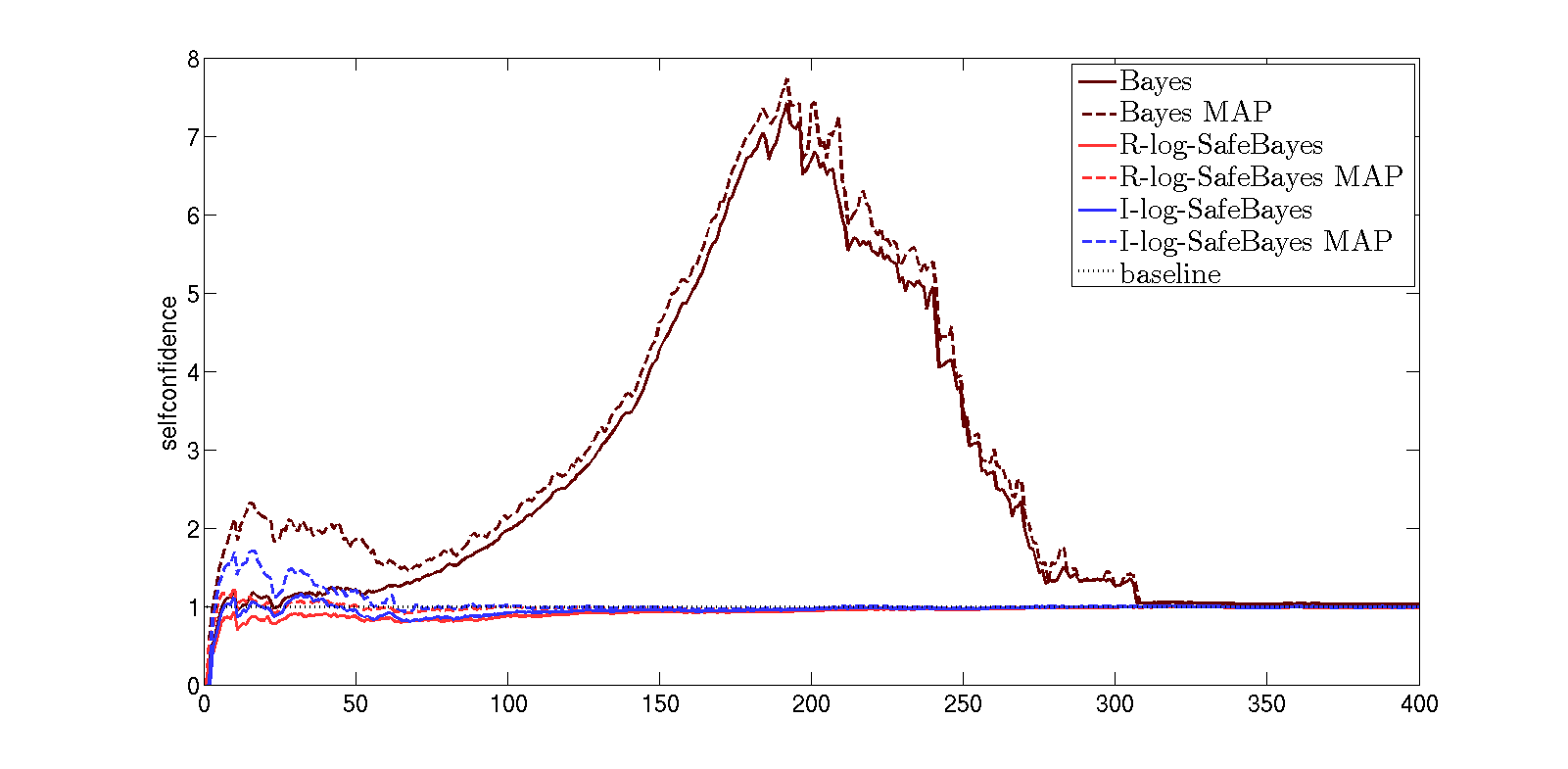} \\
\hspace*{0.2\textwidth}
\includegraphics[width=0.6\textwidth]{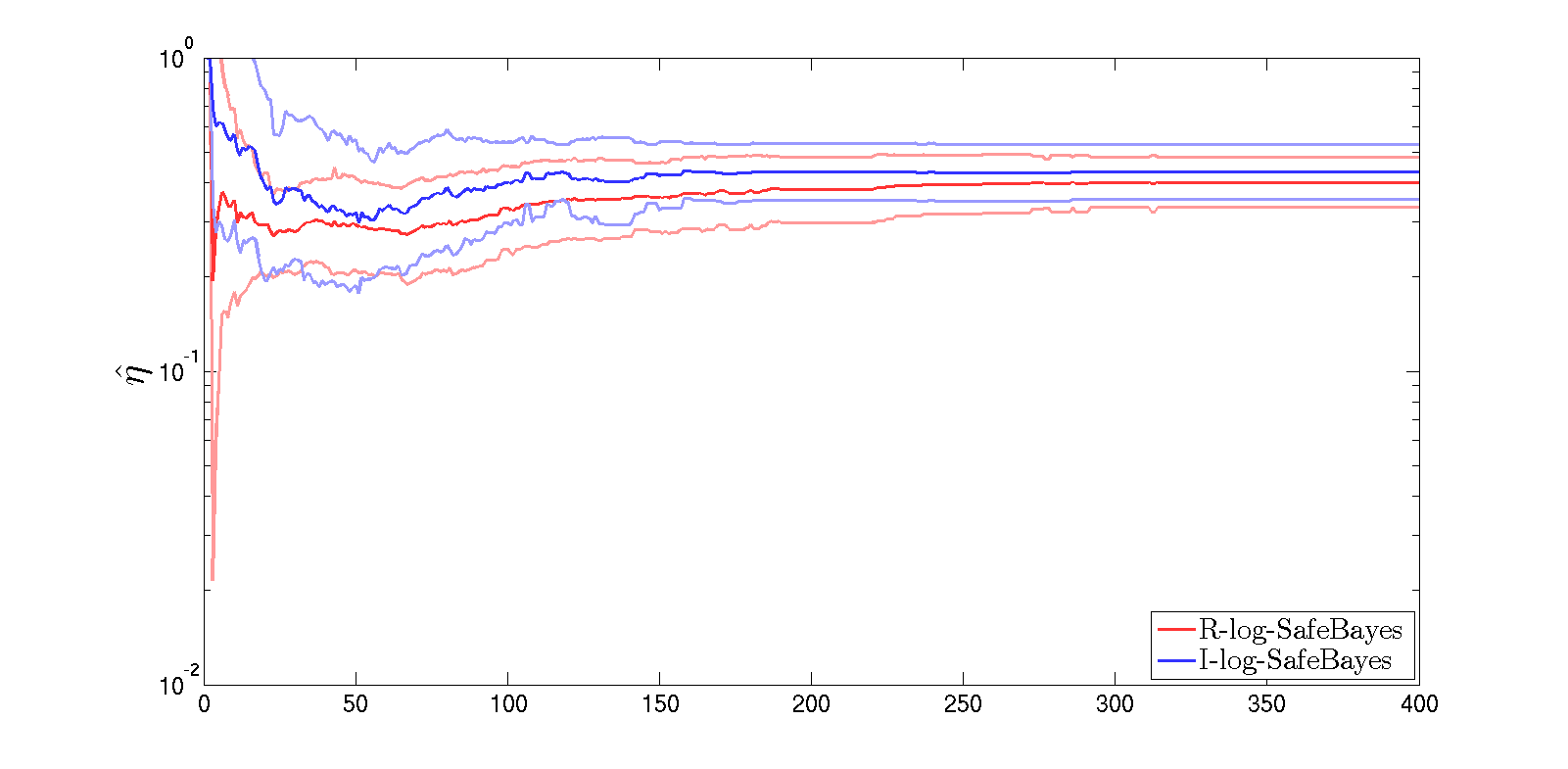}
\caption{\label{fig:mainexperimentc} Same four graphs as in
  Figure~\ref{fig:mainexperimenta}, for the wrong-model experiment
  with $\pmax = 100$.  }}
\end{figure}

\begin{figure}[h!]{\hspace*{0.2\textwidth}
\includegraphics[width=0.6\textwidth]{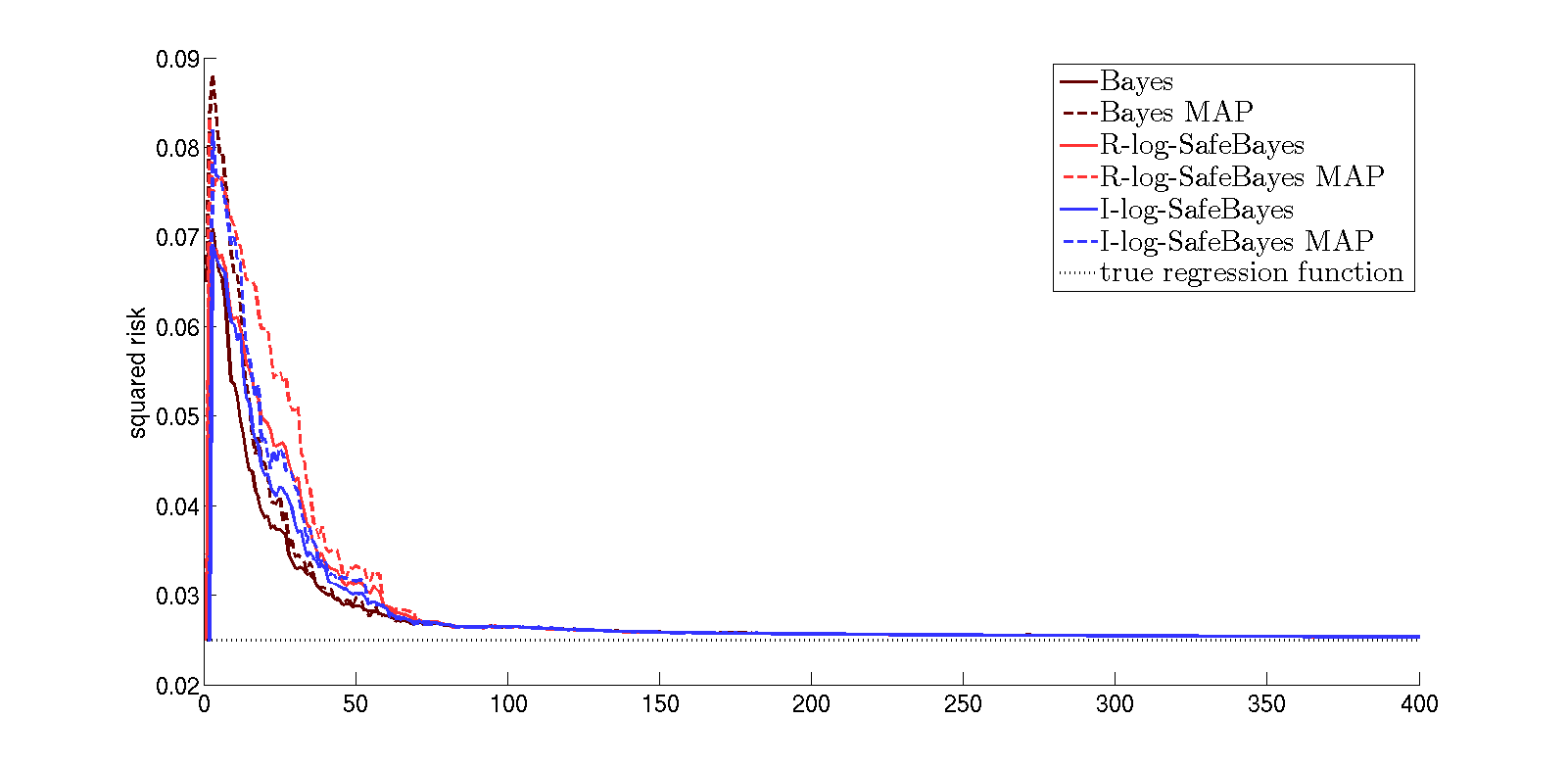} \\
\hspace*{0.2\textwidth}
\includegraphics[width=0.6\textwidth]{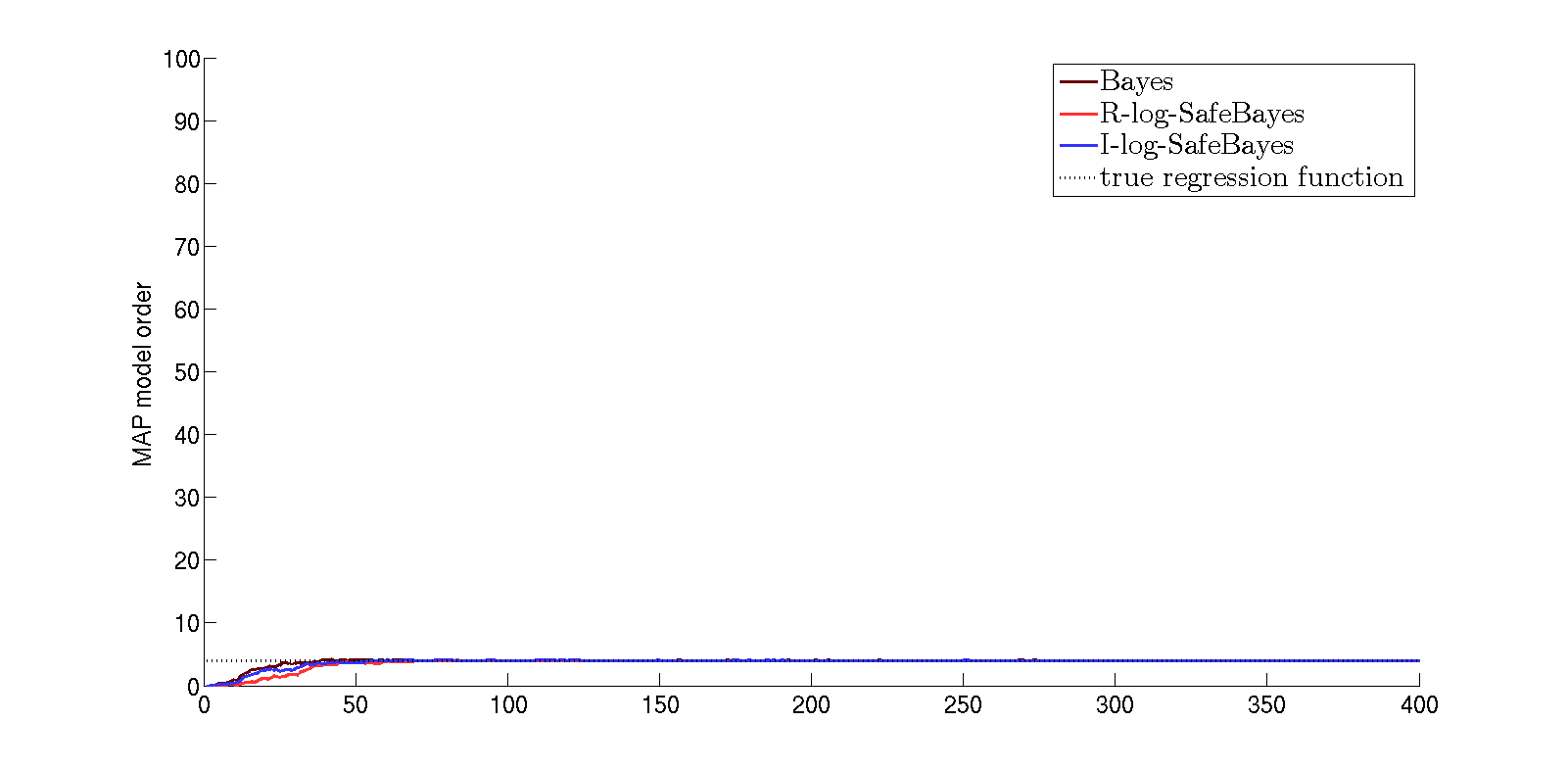} \\
\hspace*{0.2\textwidth}
\includegraphics[width=0.6\textwidth]{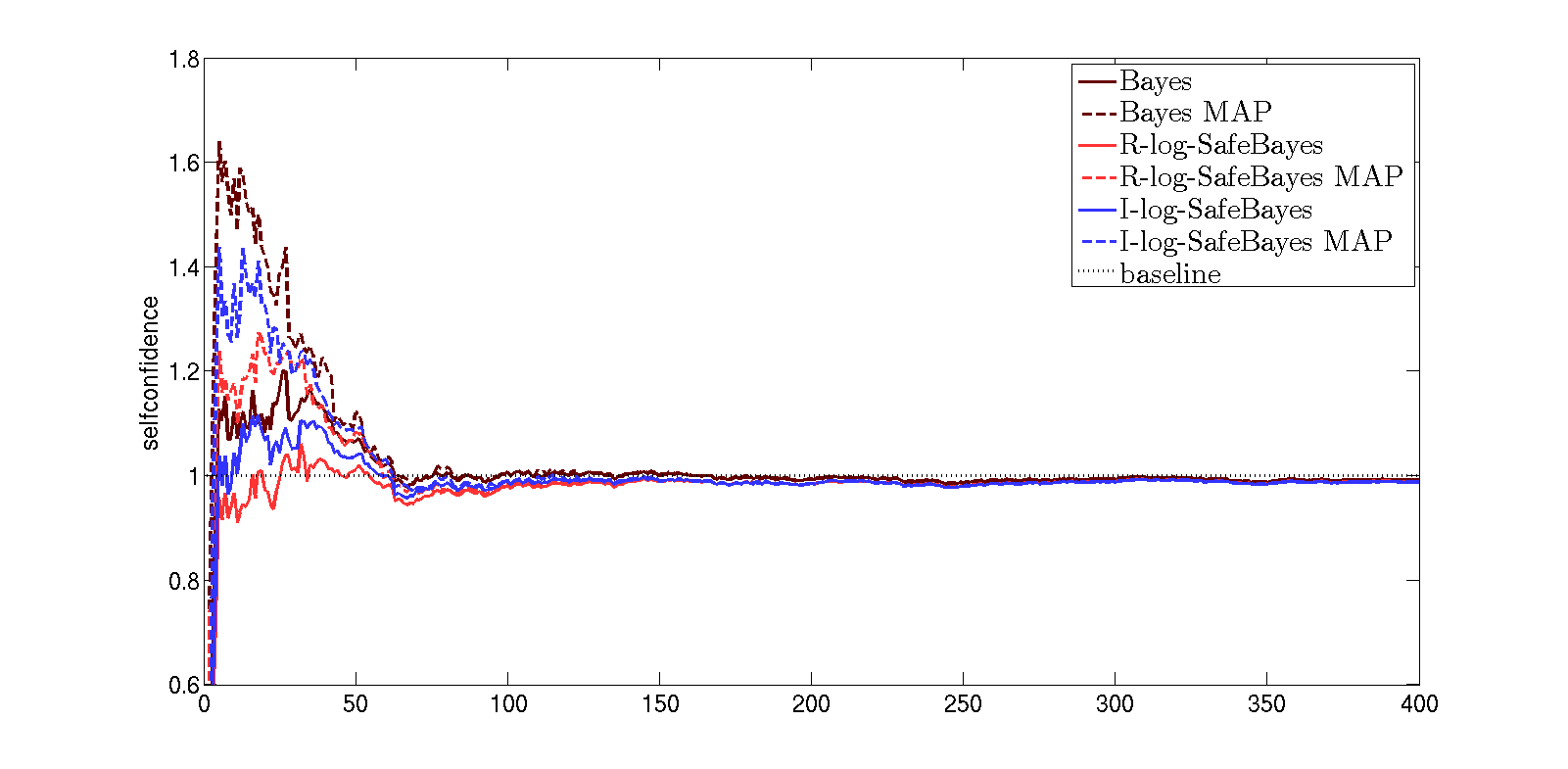} \\
\hspace*{0.2\textwidth}
\includegraphics[width=0.6\textwidth]{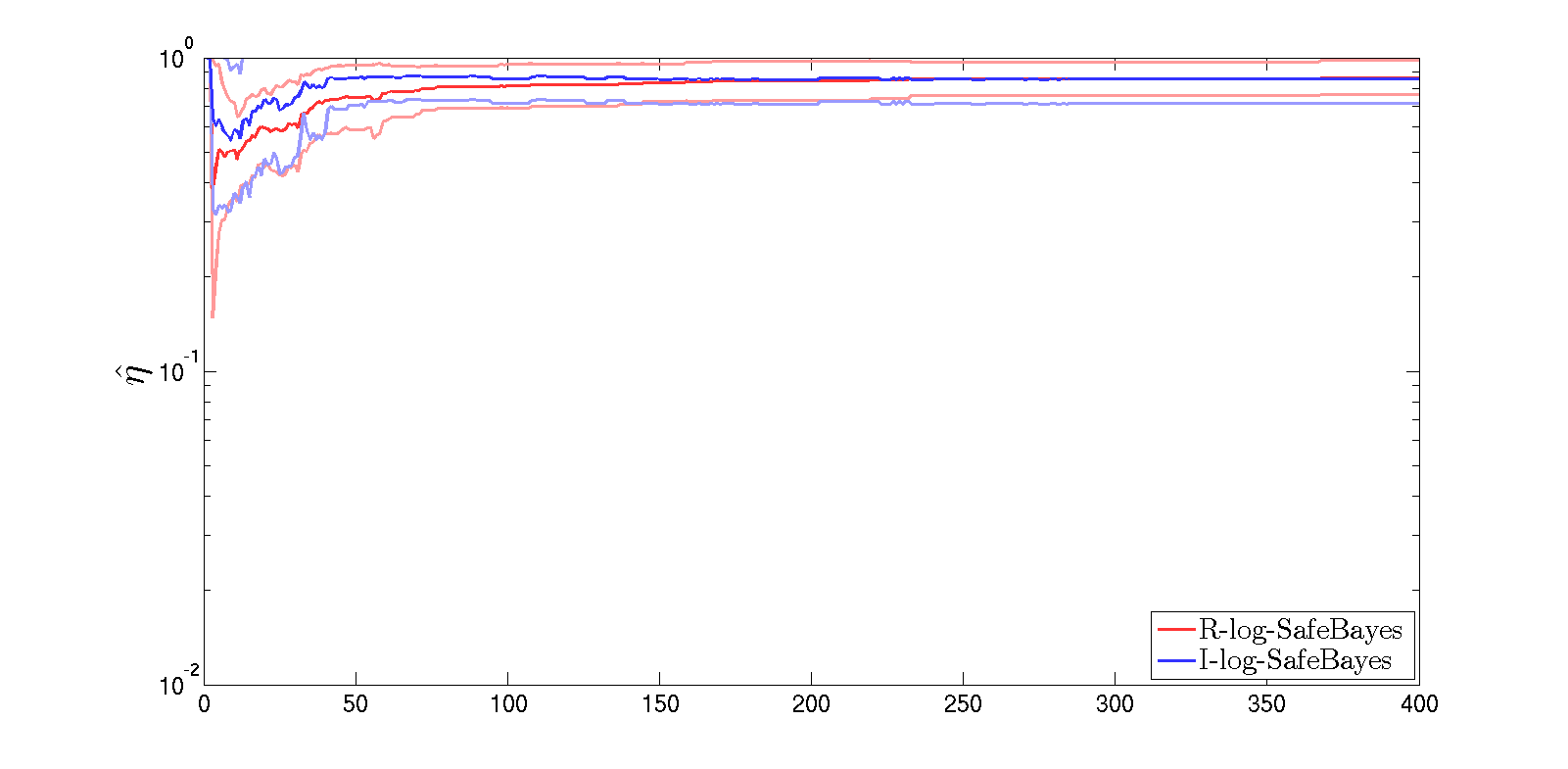}
\caption{\label{fig:mainexperimentd}
Same graphs as in Figure~\ref{fig:mainexperimenta} for the correct-model experiment with $\pmax = 100$.
}}
\end{figure}

We run the Safe Bayesian algorithm of Section~\ref{sec:safebayes} with
$z_i = (x_i,y_i)$ and $\ell_{\theta}(z_i) = - \log \dens_{\theta}(y_i
\mid x_i)$ is the (conditional) log-loss as described in that section.
As to the parameters of the algorithm (page~\pageref{alg:algorithm}),
in all experiments we set the step-size $\kappa_{\text{\sc step}}= 1/3$
and $\kappa_{\max} := 8$, i.e.\ we tried the following values of
$\eta$: $1, 2^{-1/3}, 2^{-2/3}, \ldots, 2^{-8}$. The result of the
wrong-model and correct-model experiment as described above with
$\pmax = 50$ and $\pmax = 100$, respectively, are given in
Figure~\ref{fig:mainexperimenta}--\ref{fig:mainexperimentd}.

\paragraph{Conclusion 1: Bayes performs well if model-correct, and
  dismally in model-incorrect experiment}
The four figures show that standard Bayes behaves excellently
in terms of all quality measures (square-risk, MAP model
identification and reliability) when the model is correct, and
dismally if the model is incorrect.
\paragraph{Conclusion 2: if (and only if) model incorrect, then the
  higher $\pmax$, the worse Bayes gets}

We see from Figure~\ref{fig:mainexperimentb}
and~\ref{fig:mainexperimentd} that standard Bayes behaves excellently
in terms of all quality measures (square-risk, MAP model
identification and reliability) when the model is correct,
both if
$\pmax=50$ and if $\pmax = 100$, the behavior at $\pmax=100$ being
essentially indistinguishable from the case with $\pmax = 50$. These
and other (unreported) experiments strongly suggests that, when the
data are sampled from a low-dimensional model, then, when the model is
correct, standard Bayes is unaffected (does not get confused) by
adding additional high-dimensional models to the model space. Indeed,
the same is suggested by various existing Bayesian consistency
theorems, such as those by \citet{Doob49,GhosalGV00,Zhang06density}.

At the same time, from Figure~\ref{fig:mainexperimenta}
and~\ref{fig:mainexperimentc} we infer that standard Bayes behaves
very badly in all three quality measures in our (admittedly very
`evilly chosen') model-wrong experiment. Eventually, at very large
sample sizes, Bayes recovers, but the main point here to notice is
that the $n$ at which a given level of recovery (measured in, say,
square-loss) takes place is much higher for the case $\pmax = 100$
(Figure~\ref{fig:mainexperimentc}) than for the case $\pmax = 50$
(Figure~\ref{fig:mainexperimenta}). This strongly suggests that, when
the model is incorrect but the best approximation lies in a
low-dimensional submodel, then standard Bayes gets confused by adding
additional high-dimensional models to the model space --- recovery
takes place at a sample size that increases with $\pmax$. Indeed, the
graphs strongly suggest that in the case that $\pmax = \infty$ (with
which we cannot experiment), Bayes will be inconsistent in the sense
that the risk of the posterior predictive will never ever reach the
risk attainable with the best submodel. \cite{GrunwaldL07} showed that
this can indeed happen with a simple, but much more unnatural
classification model; the present result indicates (but does not
prove) that it can happen with our standard model as well.

\paragraph{Conclusion 3: $R$-log-SafeBayes and $I$-log-SafeBayes generally
  perform well}
Comparing the four graphs for SafeBayes and $I$-log-SafeBayes, we see
that they behave quite well for {\em both\/} the model-correct and the
model-wrong experiments, being slightly worse than, though still
competitive to, standard Bayes when the model is correct and
incomparably better when the model is wrong. Indeed, in the
wrong-model experiments, about half of the data points are identical
and therefore do not provide very much information, so one would
expect that if a `good' method achieves a given level of square-risk
at sample size $n$ in the correct-model experiment, it achieves the
same level at about $2n$ in the incorrect-model experiment, and this
is indeed what happens. Also, we see from comparing
Figure~\ref{fig:mainexperimentc} and~\ref{fig:mainexperimentd} on the
one hand to Figure~\ref{fig:mainexperimenta}
and~\ref{fig:mainexperimentb} on the other that adding additional
high-dimensional models to the model space hardly affects the results
--- like standard Bayes when the model is correct, SafeBayes does not
get confused by the additional, larger model space. 

\paragraph{Secondary Conclusions}
We see that both types of SafeBayes converge quickly to the right
(pseudo-true) model order, which is pleasing since they were not
specifically designed to achieve this. Whether this is an artifact of
our setting or holds more generally would, of course, require further
experimentation.  We note that at small sample sizes, when both types
of SafeBayes still tend to select an overly simple model,
$I$-log-SafeBayes has significantly more variability in the model
chosen-on-average; it is not clear though whether this is `good' or
`bad'. We also note that the $\eta$'s chosen by both versions are very
similar for all but the smallest sample sizes, and are consistently
smaller than $1$. When instead of the full $\eta$-generalized
posteriors, the $\eta$-generalized posterior conditioned on the MAP
$\breve{\nc}_{\text{map}}$ is used, the behavior of all method consistently
deteriorates a little, but never by much.

For lack of space in the graphs, we did not show the Ces\`aro-versions
of Bayes, $R$-log-SafeBayes and $I$-log-SafeBayes (methods 3(a), 3(b), 3(c) in
Section~\ref{sec:statisticsshown}). Briefly, the curves look as
follows: Ces\`aro-Bayes performs significantly better than standard
Bayes in all three quality measures in the wrong-model experiments,
but is still far from competitive with the two (full-posterior)
SafeBayes versions. When Ces\`aroified, the SafeBayes methods become a
bit smoother but not necessarily better. Very similar behavior of
Ces\`aro (making bad methods significantly better but still not
competitive, and good methods smoother, sometimes a bit worse and
sometimes a bit better) has been explicitly depicted in the ridge
regression with varying $\sigma^2$ in Section~\ref{sec:ridge} below.
\ \pagebreak \ 
\ \pagebreak \ 

\subsection{Second Experiment: Ridge Regression, Varying $\sigma^2$} 
\label{sec:ridge}
\begin{figure}[ht!]{\hspace*{0.15\textwidth}
\includegraphics[width=0.7\textwidth]{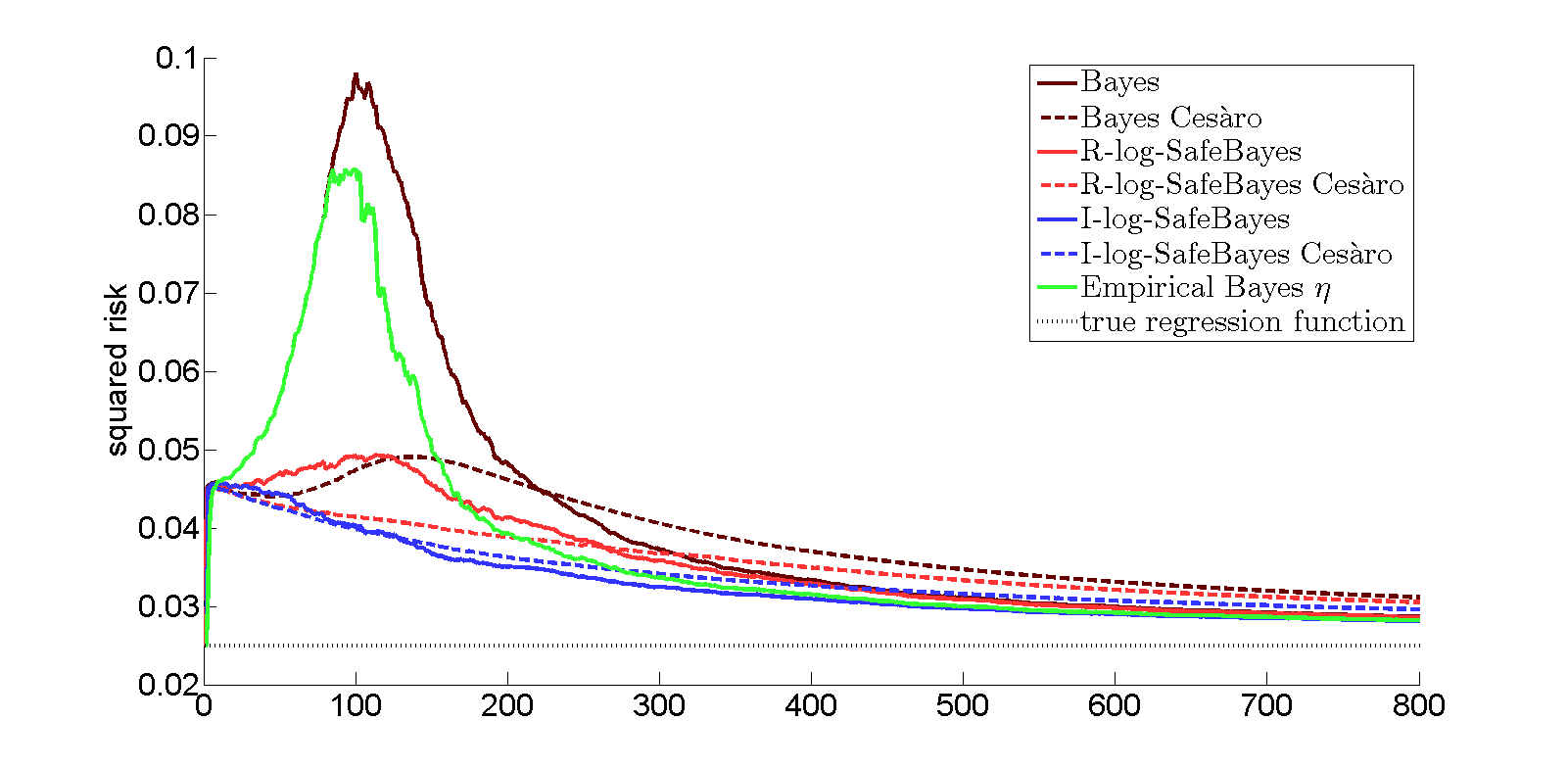} \\
\hspace*{0.15\textwidth}
\includegraphics[width=0.7\textwidth]{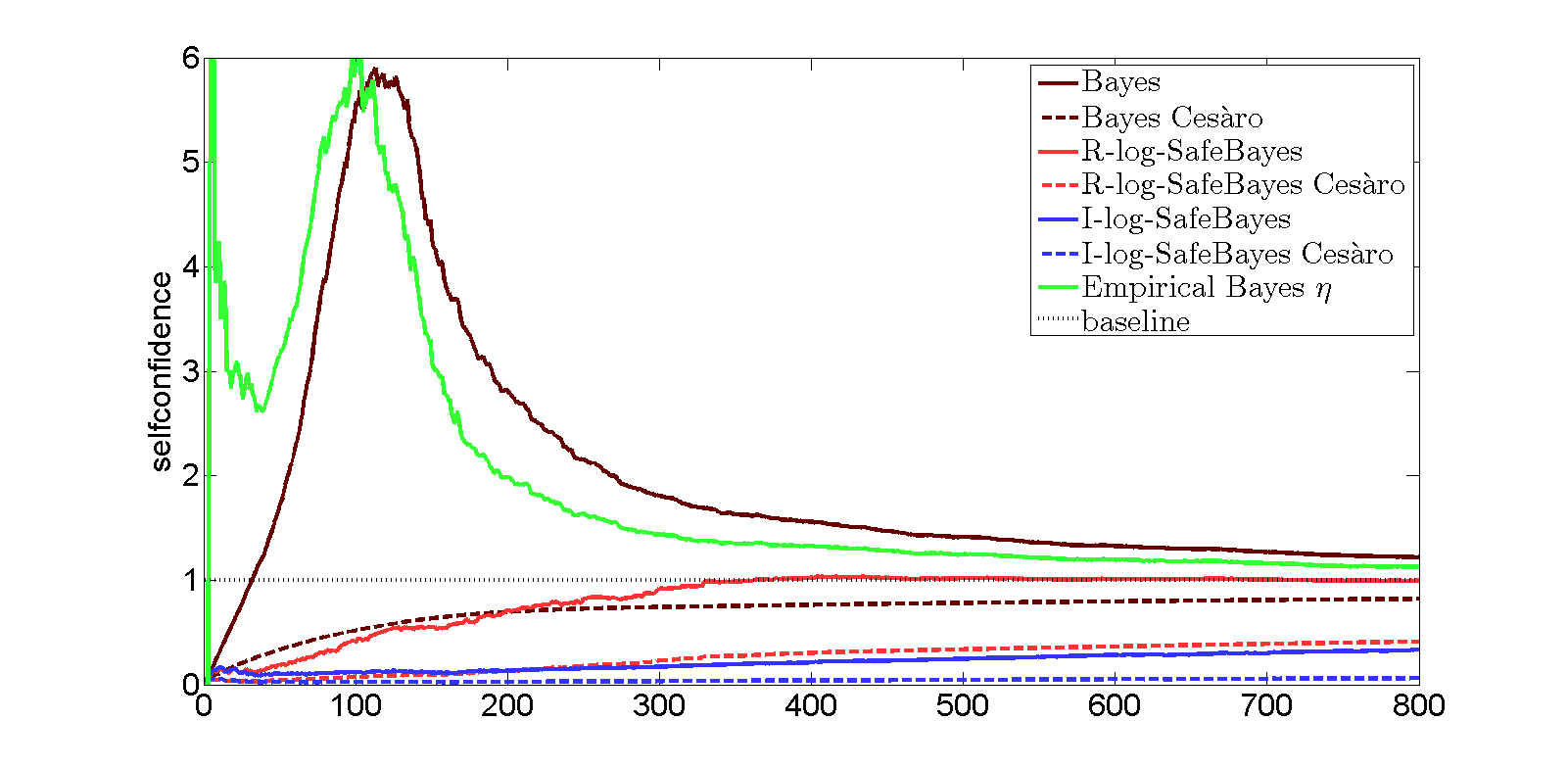}
\\
\hspace*{0.15\textwidth}
\includegraphics[width=0.7\textwidth]{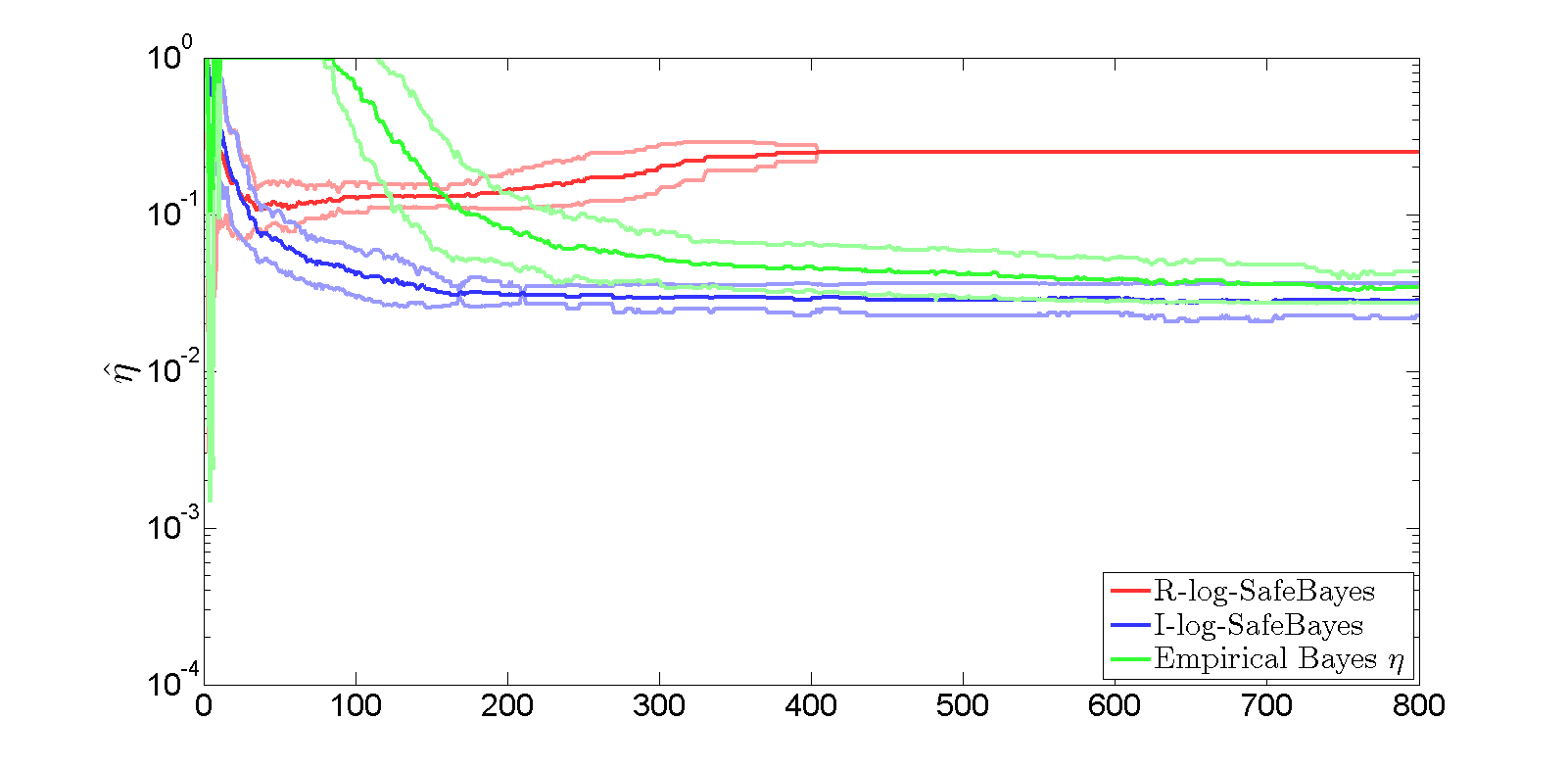}
\caption{\label{fig:ridgevarysigmaa} Bayesian Ridge Regression:
  Model-wrong experiment conditioned on $\nc := \pmax = 50$. The
  graphs (square-risk, overconfidence ratio and chosen $\eta$ as
  function of sample size) are as in
  Figure~\ref{fig:mainexperimenta}--\ref{fig:mainexperimentd}, except
  for the third graph there (MAP model order), which has no meaning here. The
  meaning of the curves is given in Section~\ref{sec:statisticsshown}
  except for {\em empirical Bayes}, explained in
  Section~\ref{sec:ridge}.  }}
\end{figure}

\begin{figure}[ht!]{\hspace*{0.15\textwidth}
\includegraphics[width=0.7\textwidth]{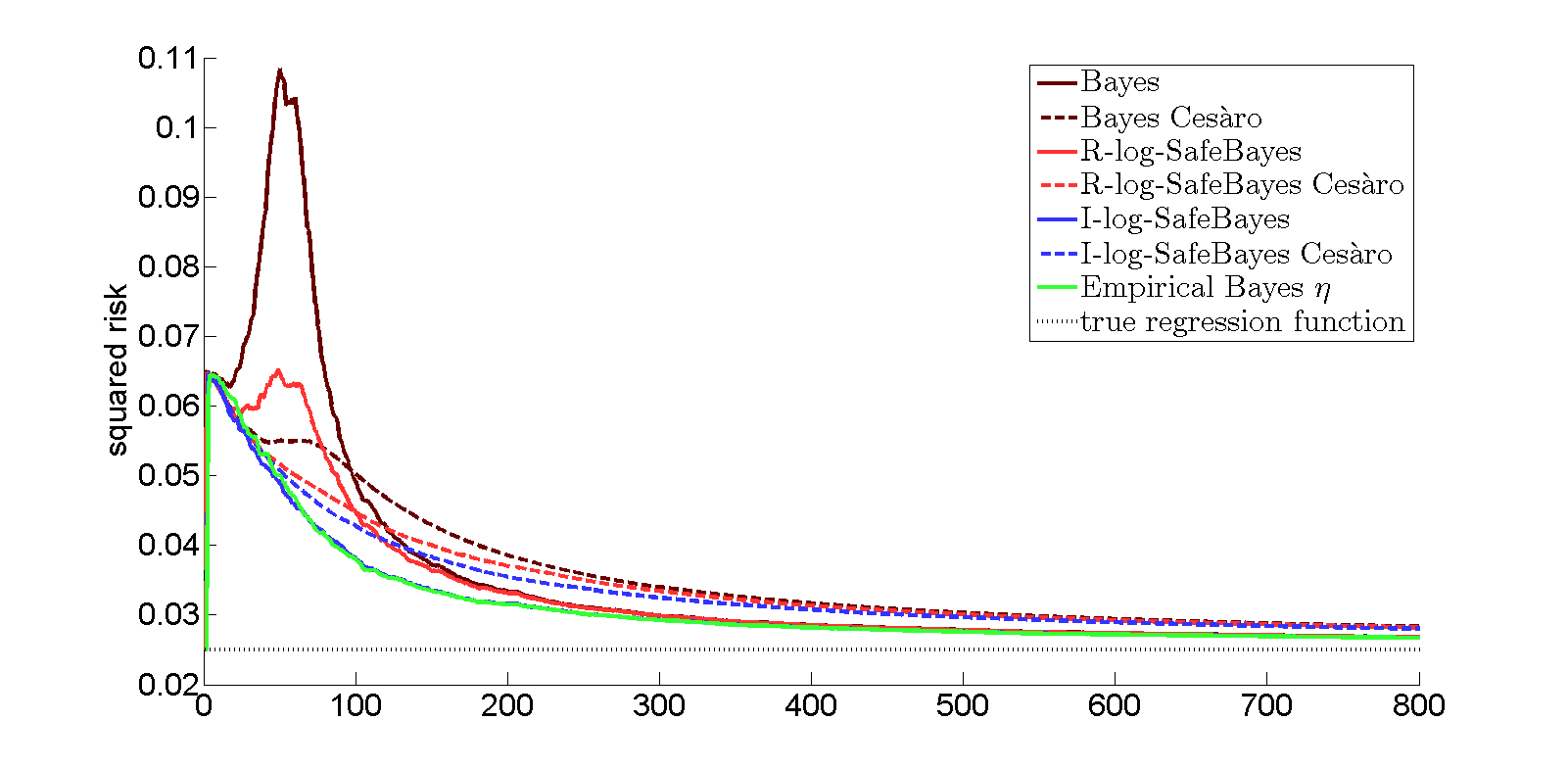} \\
\hspace*{0.15\textwidth}
\includegraphics[width=0.7\textwidth]{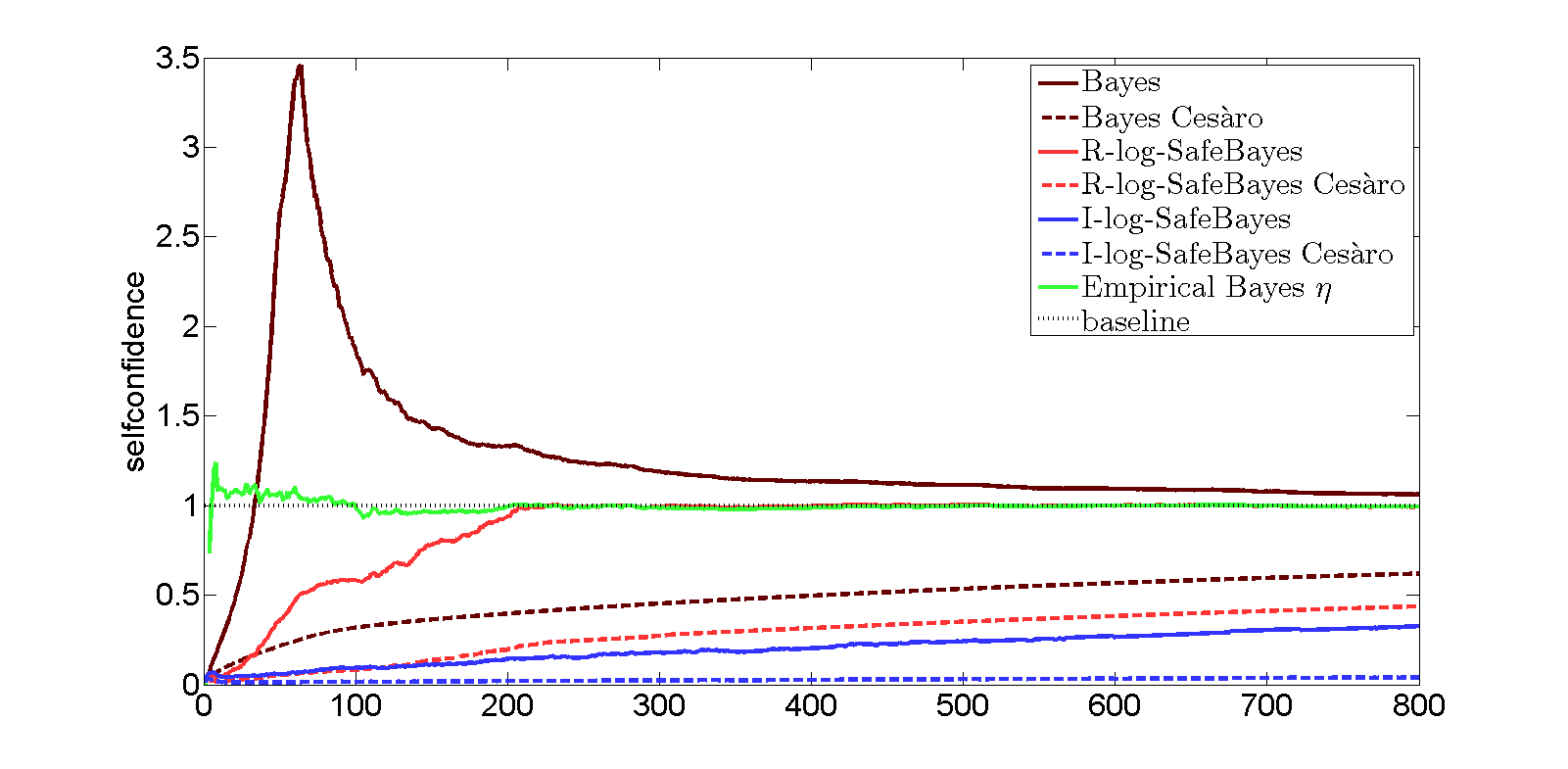}
\\
\hspace*{0.15\textwidth}
\includegraphics[width=0.7\textwidth]{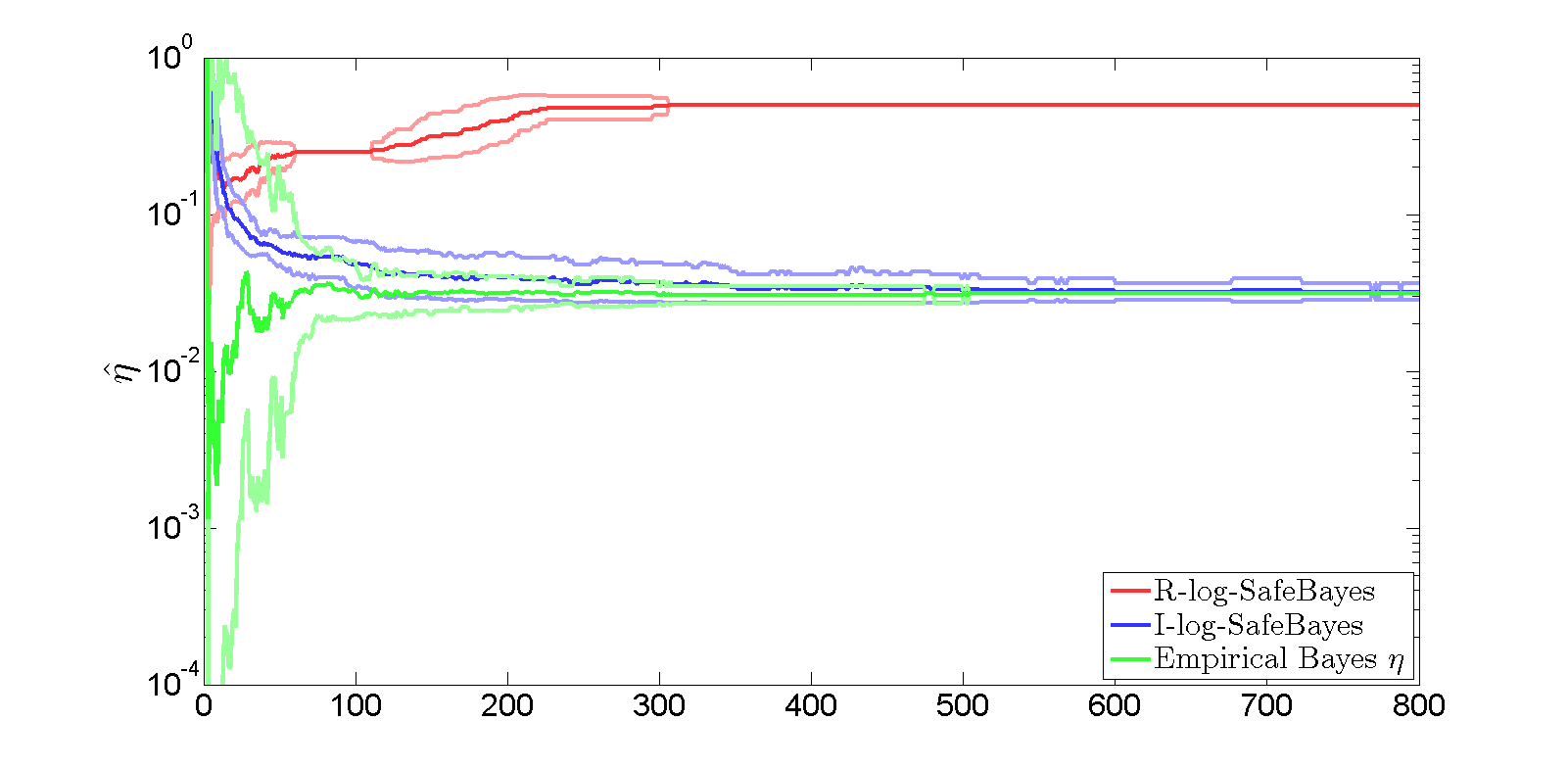}
\caption{\label{fig:ridgevarysigmab} Bayesian Ridge Regression: same
  graphs as in Figure~\ref{fig:ridgevarysigmaa}, but for the
  model-correct experiment conditioned on $\nc := \pmax = 50$. 
}}
\end{figure}

We repeat the model-wrong and model-correct experiment of
Figure~\ref{fig:mainexperimenta} and~\ref{fig:mainexperimentb}, with
just one major difference: all posteriors are conditioned on $\nc := \pmax
= 50$. Thus, we effectively consider just a fixed, high-dimensional
model, whereas the best approximation $\tilde\theta =
(50,\tilde\beta,\tilde\sigma^2)$ viewed as an element of $\cM_{\nc}$
is `sparse' in that it has only $\beta_1,\ldots, \beta_4$ not equal to
$0$. We note that the MAP model index graphs of
Figure~\ref{fig:mainexperimenta} and~\ref{fig:mainexperimentb} are
meaningless in this context (they would be equal to the constant 50)
so they are left out of the new Figure~\ref{fig:ridgevarysigmaa} and
~\ref{fig:ridgevarysigmab}. 

\paragraph{Instantiating Safe Bayes}
Since we noticed in preliminary experiments that some versions of
SafeBayes now have a tendency to select much smaller values of $\eta$ than
in the previous experiments, we now set $\kappa_{\max} = 16$ (large enough
so that in no experiment the optimal $\eta < 2^{-\kappa_{\max}}$); for
computational reasons we also increased the step size
and set $\kappa_{\text{\sc step}}= 1$.
\paragraph{Connection to Bayesian (B)ridge Regression}
From (\ref{eq:olvg}) we see that the posterior mean parameter
$\bar{\beta}_{i,\eta}$ is equal to the posterior MAP parameter and
depends on $\eta$ but not on $\sigma^2$, since $\sigma^2$ enters the
prior in the same way as the likelihood. Therefore, the square-loss
obtained when using the generalized posterior for prediction is always
given by $(y_i - x_i \bar{\beta}_{i,\eta})^2$ irrespective of whether
we use the posterior mean, or MAP, or the value of $\sigma^2$.
Interestingly, if we fix some $\lambda$ and perform standard (nongeneralized) Bayes with a
modified prior, proportional to the original prior raised to the power
$\lambda := \eta^{-1}$, then the prior becomes normal $N(\betao,
\sigma^2 \Sigma'_0)$ where $\Sigma'_0 = \eta \Sigma_0$ and the
standard posterior given $z^i$ is then (by (\ref{eq:olvg})) Gaussian
with mean
\begin{equation}\label{eq:invi}
\left( \left( \Sigma'_0\right)^{-1} + {\bf X}_n^T {\bf X} \right)^{-1} \cdot
\left( \left(\Sigma'_0\right)^{-1} \betao + {\bf X}_n^T y^n
\right) = \bar{\beta}_{i,\eta}.
\end{equation}
Thus we
see that in this special case, the (square-risk of the)
$\eta$-generalized Bayes posterior mean coincides with the (square-risk of) the standard Bayes posterior mean with prior $N(\betao,
\sigma^2 \eta \Sigma_0)$. But this means that the square-loss
obtained by $\eta$-generalized Bayes on a data sequence is 
precisely equal to the square-loss obtained by {\em Bayesian ridge
  regression\/} with penalty parameter $\lambda = \eta^{-1}$, as
defined, by, e.g., \cite{ParkC08} (to be precise, they call this
method Bayesian `Bridge' Regression with $q=2$; the choice of $q=1$ in
their formula gives their celebrated `Bayesian Lasso').  It is thus of
interest to see what happens if $\eta$ (equivalently, $\lambda$) is
determined by {\em empirical Bayes}, which is one of the methods
\cite{ParkC08} suggest. In addition to the graphs discussed earlier in
Section~\ref{sec:statisticsshown}, we thus also show the results for
$\eta$ set in this alternative way. Whereas this empirical-Bayesian
ridge regression is usually a very competitive method (indeed in our
model-correct experiment, Figure~\ref{fig:ridgevarysigmab}, it
performs best in al respects), we will see in
Figure~\ref{fig:ridgevarysigmaa} (the green line) that, just like
other versions of Bayes, it breaks down under our type of
misspecification.

We hasten to add that the correspondence between the
$\eta$-generalized posterior means and the standard posterior means
with prior raised to power $1/\eta$ only holds for the
$\bar{\beta}_{i,\eta}$ parameters. It does not hold for the
$\bar{\sigma}^2_{i,\eta}$ parameters, and thus, for fixed $\eta$, the
overconfidence of both methods may be quite different.

\paragraph{Conclusions for Model-Wrong Experiment} For most curves,
the overall picture of Figure~\ref{fig:ridgevarysigmaa} is comparable
to the corresponding model averaging experiment,
Figure~\ref{fig:mainexperimenta}: when the model is wrong, standard
Bayes shows dismal performance in terms of risk and reliability up to
a certain sample size and then very slowly recovers, whereas both
versions of SafeBayes perform quite well even for small sample sizes.
We do not show variations of the graph for $\nc = \pmax = 100$ (i.e.\ 
the analogue of Figure~\ref{fig:mainexperimentc}), since it relates to
Figure~\ref{fig:ridgevarysigmaa} in exactly the same way as
Figure~\ref{fig:mainexperimentc} relates to
Figure~\ref{fig:mainexperimenta}: with $\nc = 100$, bad square-risk
and reliability behavior of Bayes goes on for much longer (recovery
takes place at much larger sample size) and remains equally good as
for $\nc=50$ with the two versions of SafeBayes.

The results for the Ces\`aro-versions of our methods are exactly as
discussed at the end of Section~\ref{sec:modsel}.

We also see that, as we already indicated in the introduction,
choosing the learning rate by empirical Bayes (thus implementing one
version of Bayesian Bridge regression) behaves terribly. This complies
with our general theme that, to `save Bayes' in general
misspecification problems, the parameter $\eta$ cannot be chosen in a
standard Bayesian manner.

\paragraph{Conclusions for Model-Correct Experiment}
The model-correct experiment for ridge regression
(Figure~\ref{fig:ridgevarysigmab}) offers a surprise: we had expected
Bayes to perform best, and were surprised to find that the SafeBayeses
obtained smaller risk. Some followup experiments (not shown here),
with different true regression functions and different priors, shed
more light on the situation. Consider the setting in which the
coefficients of the true function are drawn randomly according to the
prior. In this setting standard Bayes performs at least as good in
expectation as any other method including SafeBayes (the Bayesian
posterior now represents exactly what an experimenter might ideally
know).  SafeBayes (still in this setting) usually chooses $\eta=1/2$
or $1/4$, and the difference in risks compared to Bayes is small. On
the other hand, if the true coefficients are drawn from a distribution
with substantially smaller variance than a priori expected by the
prior (a factor 1000 in the `correct'-model experiment of
Figure~\ref{fig:ridgevarysigmab}), then SafeBayes performs much better
than Bayes. Here Bayes can no longer necessarily be expected to have
the best performance (the model is correct, but the prior is
``wrong''), and it is possible that a slightly reduced learning rate
gives (significantly) better results. It seems that this situation,
where the variance of the true function is much smaller than its prior
expectation, is not exceptional: for example,
\cite{raftery1997bayesian} suggest choosing the variance of the prior
in such a way that a large region of parameter values receives
substantial prior mass. Following that suggestion in our experiments
already gives a variance that is large enough compared to the true
coefficients that SafeBayes performs better than Bayes even if the
model is correct.

\paragraph{A Joint Observation for Model-Wrong and Model-Correct Experiment}
Finally we note that we see an interesting difference between the two
SafeBayes versions here: $I$-log-SafeBayes seems better for risk, giving
a smooth decreasing curve in both experiments.  $R$-log-SafeBayes
inherits a trace of standard Bayes' bad behavior in both experiments,
with a nonmonotonicity in the learning curve.  On the other hand, in
terms of reliability, $R$-log-SafeBayes is consistently better than
$I$-log-SafeBayes (but note that the latter is underconfident, which is
arguably preferable over being overconfident, as Bayes is). All in al,
there is no clear winner between the two methods.

\subsection{Executive Summary: Joint Conclusions from Main and
  Additional Experiments}
\label{sec:summary}

\paragraph{Standard Bayes}
In almost all our experiments, Standard Bayesian inference fails in
its KL-associated prediction tasks (squared risk, reliability) when
the model is wrong. Adopting a different prior (such as the $g$-prior)
does not help, with two exceptions in model averaging: (a) when
Raftery's prior (Section~\ref{sec:raftery}) is used, then Bayes works
quite well, but there it fails dramatically again (in contrast to
SafeBayes) once the percentage of easy points is increased; (b) when
it is run with a fixed variance that is significantly larger than the
`best' (pseudo-true) variance $\tilde{\sigma}^2$. Moreover, in the
ridge regression experiment with fixed $\sigma^2$, we find that
standard Bayes can even perform much worse than SafeBayes when the
model is correct --- so all in all we tentatively conclude that SafeBayes is safer
to use for linear regression. 

\paragraph{Safe Bayes} $R$-square-SafeBayes is not competitive with
the other SafeBayes methods and can even get worse than Bayes
sometimes; this is due to an unwanted dependence on the specified
scale $\sigma^2$ as explained in Section~\ref{sec:priorvariations}.
The other three SafeBayes methods behave reasonably well in all our
experiments, and there is no clear winner among them.  $I$-square-
SafeBayes usually behaves excellently for the square-risk but cannot
directly be used to assess its own performance.  $I$-log-SafeBayes
usually behaves excellently in terms of square-risk as well but is
underconfident about its own performance (which is perhaps acceptable,
overconfidence being a lot more dangerous).  $R$-log-SafeBayes is
usually good in terms of square-risk though not as good as
$I$-log-SafeBayes, yet it is highly reliable. However, in
Appendix~\ref{sec:etalater}, we describe an initial idea for
discounting the importance of the first few outcomes and explain why
this might improve performance. When combined with this discounting
idea, $R$-log-SafeBayes may actually always be competitive with the
other two methods in terms of square-risk as well.
\paragraph{Learning $\eta$ in Bayes- or Likelihood Way Fails}
Despite its intuitive appeal, fitting $\eta$ to the data by e.g.\ 
empirical Bayes fails both in the model-wrong ridge experiment with a
prior in $\sigma^2$, where it amounts to Bayesian ridge regression
(Figure~\ref{fig:ridgevarysigmaa}) and in the model-wrong
fixed-variance ridge experiment (where it amounts to a method for
learning the variance, see Section~\ref{sec:ridgefixedsigma}).
\paragraph{Robustness of Experiments}
It does not matter whether the $X_{i1}, X_{i2}, \ldots$ are
independent Gaussian, uniform or represent polynomial basis functions:
all phenomena reported here persist for all choices. 
If the `easy' points are not precisely $(0,0)$, but have themselves a
small variance in both dimensions, then all phenomena reported here
persist, but on a smaller scale. 
\paragraph{Centering}
We repeated several of our experiments with centered data, i.e.\ 
preprocessed data so that the empirical average of the $Y_i$ is
exactly 0 \cite{raftery1997bayesian,HastieTF01}. In none of our
experiments did this affect any results. While this is not further
mentioned in the appendix, there we also looked at the case where the
true regression function has an intercept far from $0$, and data are
{\em not\/} centered. This hardly affected the SafeBayes methods.

\paragraph{Other Methods}
We also repeated the wrong-model experiment for other methods of model
selection: AIC, BIC, and various forms of cross-validation. Briefly,
we found that all these have severe problems with our data as well.
Whereas in these experiments, cross-validation was used to identify a
model index $p$ and $\eta$ played no role, in our final experiment we
used leave-one-out cross-validation again to learn $\eta$ itself. With
the squared error loss it worked fine, which is not too surprising
given its close similarity to $I$-square-SafeBayes.  However, when we
tried it with log-loss (as a likelihoodist or information-theorist
might be tempted to do), it behaved terribly.

\pagebreak
\section{Bayes' Behavior Explained}
\label{sec:explanation}
In this section we explain how anomalous behavior of the Bayesian
posterior may arise, taking a frequentist perspective.
Section~\ref{sec:variance} is merely provided to give some initial
intuition and may be skipped.

\subsection{Explanation I: Variance Issues}
\label{sec:variance}
\begin{example}{\bf [Bernoulli]}\label{ex:bernoulli}
Consider the
  following very simple scenario: our `model' consists of two
  Bernoulli distributions, $\cM = \{P_{\theta} \mid \theta \in
  \{0.2,0.8\} \}$, with $P_{\theta}$ expressing that $Y_1, Y_2, \ldots
  \sim \text{i.i.d.~{\sc Ber}}(\theta)$. We perform Bayesian inference
  based on a uniform prior on $\cM$. Suppose first that the data are,
  in fact, sampled i.i.d.~from $P_{\theta^*}$, where $\theta^*$ is the
  `true' parameter. The model is misspecified, in particular we will
  take a $\theta^* \not \in \{0.2,0.8\}$. The log-likelihood ratio
  between the two distributions for data $Y^n$ with $n_1$ ones and
  $n_0 = n- n_1$ zeroes, measured for convenience in bits (base 2), is
  given by
\begin{equation}\label{eq:BernoulliLL}
L = \log_2 \frac{\dens_{0.8}(Y^n)}{\dens_{0.2}(Y^n)} = 
\log_2 \frac{(0.8)^{n_1} (0.2)^{n_0}}{(0.2)^{n_1} (0.8)^{n_0}} = 2 (n_1 - n_0).
\end{equation}
With uniform priors, the posterior will prefer $\theta = 0.2$ as soon
as $L < 0$. 

First suppose $\theta^* = 1/2$. Then both distributions in $\cM$ are
equally far from $\theta^*$ in terms of KL divergence (or any other commonly used
measure). By the central limit theorem, however, we expect that the
probability that $|L| > \sqrt{n}/2$ is larger than a constant for all
large $n$; in this particular case we numerically find that, for all
$n$, it is larger than $0.32$.

This implies, that, at each $n$, with `true' probability at least
$0.32$, $\min_{\theta \in \{0.2, 0.8 \}} \pi(\theta \mid Y^n) \approx
2^{- \sqrt{n}/2}$. Thus, there is a nonnegligible `true' probability
that the posterior on one of the two distributions is negligibly
small, and a naive Bayesian who adopted such a model would be strongly
convinced that the other distribution would be better even though both
distributions are equally bad. While this already indicates that
strange things may happen under misspecification, we are of course
more interested in the situation in which $\theta^* \neq 1/2$, so that
one of the two distributions in $\cM$ is truly `better'. Now, if the
`true' parameter $\theta^*$ is within $O(1/\sqrt{n})$ of $1/2$, then,
by the central limit theorem, the probability that $L < 0$ is
nonnegligible.  For example, if $\theta^*$ is exactly $1/2+
1/\sqrt{n}$, then this probability is larger than 0.16 for all $n$.
Thus, for values of $\theta^*$ this close to $1/2$, there is no way we
can even expect Bayes to learn the `best' value. For fixed
(independent of $n$), larger values of $\theta^*$, like $0.6$, the
posterior will concentrate at $0.8$ at an exponential rate, but the
sample size at which concentration starts is considerably larger than
the sample sized needed when the true parameter is, in fact $0.8$. For
example, at $n=50$, $P_{0.6}(L < 0) \approx 0.1, P_{0.8}(L < 0)
\approx 2 \cdot 10^{-5}$; both probabilities go to $0$ exponentially
fast but their ratio increases exponentially with $n$.  So, under a
fixed $\theta^*$, with increasing $n$, Bayes may take longer to
concentrate on the best $\tilde{\theta}$ if $\tilde\theta \neq
\theta^*$ (misspecification) than if $\tilde\theta = \theta^*$, but it
eventually `recovers' (this was seen in the ridge experiments of Section~\ref{sec:ridge}). Now, for larger models, the consequence of slower concentration
of the log-likelihood ratio $L$ is that the probability that {\em
  some\/} `bad' $P_{\theta}$ happens to `win' is substantially larger
than with a correct model.  \cite{GrunwaldL07} showed that, in a
classification context with an infinite-dimensional model, there are
so many of such `bad' $P_{\theta}$ that Bayes does not recover any
more, and the posterior keeps putting most of its mass on a bad model
for ever (although the particular bad model on which it puts its mass,
keeps changing). In this paper we empirically showed the same in a
regression problem.
\end{example}
Now one might conjecture that the issues above are caused by the fact
that the model $\cM$ is `disconnected'.  In the Bernoulli example
above, the problem indeed goes away if instead of the model $\cM$, we
adopt its `closure' $\cM' = \{ P_{\theta} \mid \theta \in [0.2,0.8]
\}$. However, high-dimensional regression problems exhibit the same
phenomenon, even if their parameter spaces are connected. It turns out
that in general, to get concentration at the same rates as if the model
were correct, the model must be {\em convex}, i.e.\ closed under
taking any finite mixture of the  densities, which is a much stronger
requirement than mere connectedness. 
For standard Gaussian regression problems with $Y \mid X \sim N(0,
\sigma^2)$, this would mean that we would have to adopt a model in
which $Y \mid X$ can be any Gaussian mixture with arbitrarily many
components --- which is clearly not practical (note that `convex' refers to the
densities, not the regression functions \citep[Section 6.3.5]{GrunwaldL07}).

\subsection{Explanation II: Good vs.\ Bad Misspecification}\label{sec:badmisspecification} 
\cite{Barron98} showed that sequential Bayesian prediction under a
logarithmic score function shows excellent behavior in a cumulative
risk sense; for a related result see \cite[Lemma
4]{barron1999consistency}. Although \cite{Barron98} focuses on the
well-specified case, this assumption is not required for the proof and
the result still holds even if the model $\cM$ is completely wrong.
For a precise description and proof of this result emphasizing that it
holds under misspecification, see \cite[Section 15.2]{Grunwald07}. At
first sight, this leads to a paradox, as we now explain.
\paragraph{A Paradox?}
Let $\tilde{\theta}$ index the KL-optimal distribution in $\Theta$ as
in Section~\ref{sec:optimal}. The result of \cite{Barron98}
essentially says that, for arbitrary models $\Theta$, for all $n$,
\begin{equation}\label{eq:barron}
\Exp_{Z^n \sim P^*} \left[ \sum_{i=1}^n \logrisk(\Pi \mid
  Z^{i-1})- \logrisk(\tilde{\theta})
\right]
\leq  \red_n,
\end{equation}
where $\logrisk(W)$, for a distribution $W$ on $\Theta$, is defined as
the log-risk obtained when predicting by the $W$-mixture of
$\dens_{\theta}$, i.e. 
\begin{equation}\label{eq:logriskW}
\logrisk(W) = \Exp_{X,Y
  \sim P^*} [ - \log \Exp_{\theta \sim W} \dens_{\theta}(Y \mid X)].
\end{equation}
In (\ref{eq:barron}), this coincides with log-risk of the Bayes
predictive density $\pbayes(\cdot \mid Z^{i-1})$, as defined by
(\ref{eq:bayespred}). Here, as in the remainder of this
section, we look at the standard Bayes predictive density, i.e.\ $\eta
= 1$. $\red_n$ is the so-called {\em relative expected stochastic
  complexity\/} or {\em redundancy\/} \citep{Grunwald07}, which
depends on the prior and for `reasonable' priors is typically {\em
  small}. The result thus means that, when sequentially predicting
using the standard predictive distribution under a log-scoring rule,
one does not lose much compared to when predicting with the log-risk
optimal $\tilde\theta$.

When $\cM$ is a union of a finite or countably infinite
number of parametric exponential families and $\tilde{\nc}< \infty$ is
well-defined, then, under some further regularity conditions, which
hold in our regression example, \cite{Grunwald07}, the redundancy is, up to $O(1)$, equal
to the BIC term $(\tilde{k}/2) \log n$, where $\tilde{k}$ is the
dimensionality of the smallest model containing $\tilde{\theta}$.  In
the regression case, $\cM_{\tilde{\nc}}$ has $\tilde{\nc} +2$
parameters $(\beta_0, \ldots, \beta_{\nc},\sigma^2)$, so in the two
experiments of Section~\ref{sec:mainexperiments}, $\tilde{k} = 6$.
Thus, in our regression example, 
when sequentially predicting with the standard Bayes predictive
$\pbayes(\cdot \mid Z^{i-1})$, the cumulative log-risk is at most $n
\cdot \logrisk(\tilde\theta)$ which is linear in $n$, plus a
logarithmic term that becomes comparatively negligible as $n$
increases. This is confirmed by Figure~\ref{fig:hypercompressiona}
below.  Now, for each individual $\theta = (p,\beta,\sigma^2)$ we know
from Section~\ref{sec:optimalb} that, if its log-risk is close to that
of $\tilde\theta$, then its square-risk must also be close to that of
$\tilde\theta$; and $\tilde\theta$ itself has the smallest square-risk
among all $\theta \in \Theta$. Hence, one would expect the reasoning
for log-risk to transfer to square-risk: it seems that when
sequentially predicting with the standard Bayes predictive
$\pbayes(\cdot \mid Z^{i-1})$, the cumulative square-risk should at
most be $n$ times the instantaneous square-risk of $\tilde\theta$ plus
a term that hardly grows with $n$; in other words, the cumulative
square-risk from time $1$ to $n$, averaged over time by dividing by
$n$, should rapidly converge to the constant instantaneous risk of
$\tilde\theta$.  Yet the experiments of
Section~\ref{sec:mainexperiments} clearly show that this is {\em
  not\/} the case: Figure~\ref{fig:mainexperimenta} shows that, until
$n=100$, it is about 3  times as large.

\begin{figure}{\hspace*{0.25\textwidth}
\includegraphics[width=0.49\textwidth]{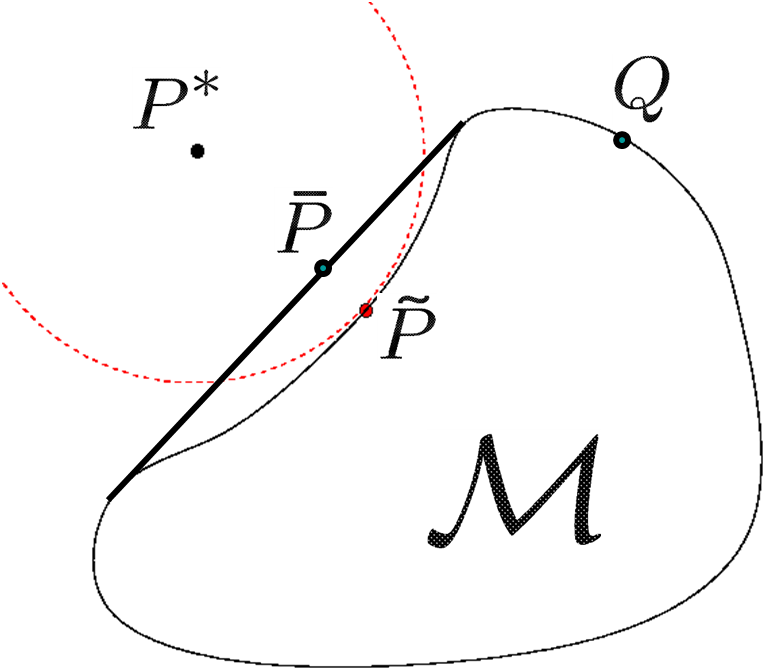}}
\caption{Benign vs.\ Bad Misspecification: $\tilde{P} = \arg
  \min_{P \in\cM} D(P^* \| P)$ is the distribution in model $\cM$ that
  minimizes KL divergence to the `true' $P^*$, but, since the model is
  nonconvex, the Bayes predictive distribution $\bar{P}$ may happen to
  be very different from any $P \in \cM$. When this
  happens, we can have `bad misspecification' and then it  may be necessary to decrease the learning rate (in this
  simplistic drawing $\bar{P}$ is a mixture of just two distributions;
  in our regression example it mixes infinitely many). Yet if
  $P^*$ were such that $\inf_{P \in \cM} D(P^* \| P)$ does not
  decrease if the infimum is taken over the convex hull of $\cM$
  (e.g.\ 
  if $Q$ rather than $\tilde{P}$ reached the minimum), then any
  learning rate $\eta < 1$ is fine (`benign' misspecification).
In the picture, we even have $D(P^* \| \bar{P}) < D(P^* \| \tilde{P})$; in this case we can get  hypercompression.
\label{fig:badmix}}
\end{figure}

This `paradox' is resolved once we realize that the Bayesian
predictive density $\pbayes(\cdot \mid ^{i-1})$ is a {\em mixture\/}
of various $\dens_{\theta}$, and not necessarily similar to
$\dens_{\theta}$ for any individual $\theta$ --- the link between
log-risk and square-risk (\ref{eq:THEidentity}) only holds for
individual $\theta = (\nc,\beta,\sigma^2)$, not for mixtures of them.
Indeed, if at each point in time $i$, $\pbayes(\cdot \mid Z^{i})$
would be very similar (in terms of e.g.\ Hellinger distance) to some
particular $\dens_{\theta_i}$ with $\theta_i \in \Theta$, then there
would really be a contradiction. Thus, the discrepancy between the
good log-risk and bad square-risk results in fact {\em implies\/} that
at a substantial fraction of sample sizes $i$, $\pbayes(\cdot \mid
Z^{i})$ must be substantially different from {\em every\/} $\theta \in
\Theta$. In other words, {\em the posterior is not concentrated at
  such $i$}.  A cartoon picture of this situation is given in
Figure~\ref{fig:badmix}: the Bayes predictive achieves small log-risk
because it mixes together several distributions into a single
predictive distribution which is very different from any particular
single $\dens_{\theta} \in \cM$. By Barron's bound, (\ref{eq:barron}),
the resulting $\pbayes(\cdot \mid Z^i)$ must, averaged over $i$, have
at most a risk almost as small as the risk of $\tilde{\theta}$.  We
can thus, at least informally, distinguish between ``benign'' and
``bad'' misspecification. Bad misspecification occurs if there is a
nonnegligible probability that for a range of sample sizes, the
predictive distribution is substantially different from any of the
distributions in $\cM$.  As Figure~\ref{fig:badmix} suggests, `bad'
misspecification cannot occur for convex models $\cM$ --- and indeed,
the results by \cite{Li99} suggest that for such models consistency
holds under weak conditions for any $\eta < 1$, even under misspecification.
\subsection{Hypercompression}\label{sec:hypercompression}
The picture suggests that, if, as in
our regression model, the model is nonconvex (i.e.\ the set of
densities $\{ \dens_{\theta} \mid \theta \in \Theta \}$ is not closed
under taking mixtures), then $\pbayes(\cdot \mid Z^i)$ might in fact
be significantly {\em better\/} in terms of log-risk than the best
$\tilde{\theta}$, and its individual constituents might even all be
substantially worse than $\tilde{\theta}$.  If this were indeed the
case then, with high $P^*$-probability, we would also get the
analogous result for an actual sample (and not just in expectation):
the cumulative log-risk obtained by the Bayes predictive should be
significantly smaller than the cumulative log-risk achieved with the
optimal $\tilde{\dens}$.  Figure~\ref{fig:hypercompressiona} below
shows that this indeed happens with our data, until $n \approx 100$.
\paragraph{The No-Hypercompression Inequality}
In fact, Figure~\ref{fig:hypercompressiona} shows a phenomenon that is
virtually impossible if the Bayesian's model and prior are `correct'
in the sense that data $Z^n$ would behave like a typical sample from
them:
it
easily follows from Markov's inequality (for details see \cite[Chapter 3]{Grunwald07})
that, letting $\Pi$ denote the Bayesian's joint distribution on $\Theta
\times \cZ^n$, 
for each $K \geq 0$,
$$
\Pi \left\{ (\theta,Z^n)\ :\ 
\sum_{i=1}^n \left(
- \log \pbayes(Y_i \mid X_i, Z^{i-1}) \right) \leq \sum_{i=1}^n \left(
- \log \dens_{\theta}(Y_i \mid X_i, Z^{i-1} )\right) -K \right\} \leq e^{-K},
$$
i.e.\ the probability that the Bayes predictive $\pbayes$ cumulatively outperforms
$\dens_{\theta}$, with $\theta$ drawn from the prior, by $K$ log-loss units is exponentially
small in $K$. Figure~\ref{fig:hypercompressiona} below thus shows that
at sample size $n \approx 90$, an a-priori formulated event has
happened of probability less than $e^{-30}$, clearly indicating that
something about our model or prior is quite wrong.

Since the difference in cumulative log-loss between $\pbayes$ and
$\dens_{\theta}$ can be interpreted as the amount of bits saved when coding the
data with a code that would be optimal under $\pbayes$ rather than
$\dens_{\theta}$, this result has been called the {\em no hyper-compression
  inequality\/} by \cite{Grunwald07}. The figure shows that for our
data, we have substantial hypercompression.

\paragraph{The Safe Bayes Error Measure}
As seen from (\ref{eq:kibbeling}), 
SafeBayes measures the performance of $\eta$-generalized Bayes not by
the cumulative log-loss, as standard Bayes does, but
instead by the cumulative posterior-expected error when predicting by
drawing from the posterior. One way to interpret this alternative error measure is that, at least in expectation, we cannot get hypercompression. Defining (compare to (\ref{eq:logriskW})!)
\begin{equation}\label{eq:logriskR}
\Rlogrisk(W) = \Exp_{X,Y
  \sim P^*} E_{\theta \sim W} [ - \log  \dens_{\theta}(Y \mid X)],
\end{equation}
we get by Fubini's theorem, 
\begin{equation}\label{eq:coldfinger}
\Rlogrisk(W) - \logrisk(\tilde\theta) = E_{\theta \sim W} \Exp_{X,Y
  \sim P^*}  [[ - \log  \dens_{\theta}(Y \mid X)] - [ - \log  \dens_{\tilde\theta}(Y \mid X)]
] \geq 0,
\end{equation}
where the inequality follows by definition of $\tilde\theta$ being
log-risk optimal among $\Theta$.  There is thus a crucial difference
between $\Rlogrisk$ and $\logrisk$ --- for the latter we just argued
that, under misspecification, $\logrisk(W) - \logrisk(\tilde\theta)
\leq 0$ is very well possible.  Thus, in contrast to predicting with
the mixture density $\Exp_{\theta \sim W} f_{\theta}$, prediction
by randomization (first sampling $\theta \sim W$ and then predicting
with the sampled $f_{\theta}$) cannot `exploit' the fact that mixture
densities might have smaller log-risk than their components. Thus, if
the difference (\ref{eq:coldfinger}) is small, then $W$ must put most
of its mass on distributions $\theta \in \Theta$ that have small
log-risk themselves. For {\em individual\/} $\theta$, we know that
small log-risk implies small square- risk. This implies that if
(\ref{eq:coldfinger}) is small, then the (standard) posterior is
concentrated on distributions with small $R$-square-risk.

\paragraph{Experimental Demonstration of Hypercompression for Standard
  Bayes}
Figure~\ref{fig:hypercompressiona} and
Figure~\ref{fig:hypercompressionb} show the predictive capabilities of
Standard Bayes in the wrong model example in terms of {\em
  cumulative\/} and {\em instantaneous log-loss\/} on a simulated
sample. 
The graphs clearly demonstrate hypercompression: the Bayes predictive
cumulatively performs {\em better\/} than the best single model/the
best distribution in the model space, until at about $n \approx 100$
there is a phase transition. At individual points, we see that it
sometimes performs a little worse, and sometimes (namely at the `easy'
$(0,0)$ points) much better than the best distribution. We also see
that, as anticipated above, randomized and in-model Bayesian
prediction do {\em not\/} show hypercompression and in fact perform
terribly on the log-loss until the phase transition at $n=100$, when
they becomes almost as good as standard Bayes. We see that for $\eta =
1$, they perform much worse. An important conclusion is that {\em if
  we are only interested in log-loss prediction, it is clear that we
  just want to use Bayes rather than randomized or in-model prediction
  with different $\eta$}.
\begin{quote}
As an aside, we note that the first few outcomes have a dramatic effect on
cumulative $R$-and $I$-log-loss (it disappears from
Figure~\ref{fig:hypercompressionb}); this may be due to the fact that
our densities --- other than those considered by \cite{Grunwald12} ---
have unbounded support so that there is no $v$ such that
Theorem~\ref{thm:concentration} below holds. This observation inspired
the idea described in Appendix~\ref{sec:etalater} about ignoring the
first few outcomes when determining the optimal $\eta$.
Also, we emphasize that the
hypercompression phenomenon takes places more generally, not just in
our regression setup --- for example, the classification inconsistency
noted by \cite{GrunwaldL07} can be understood in terms of
hypercompression as well.
\end{quote}
\begin{figure}{\hspace*{0.125\textwidth}
\includegraphics[width=0.75\textwidth]{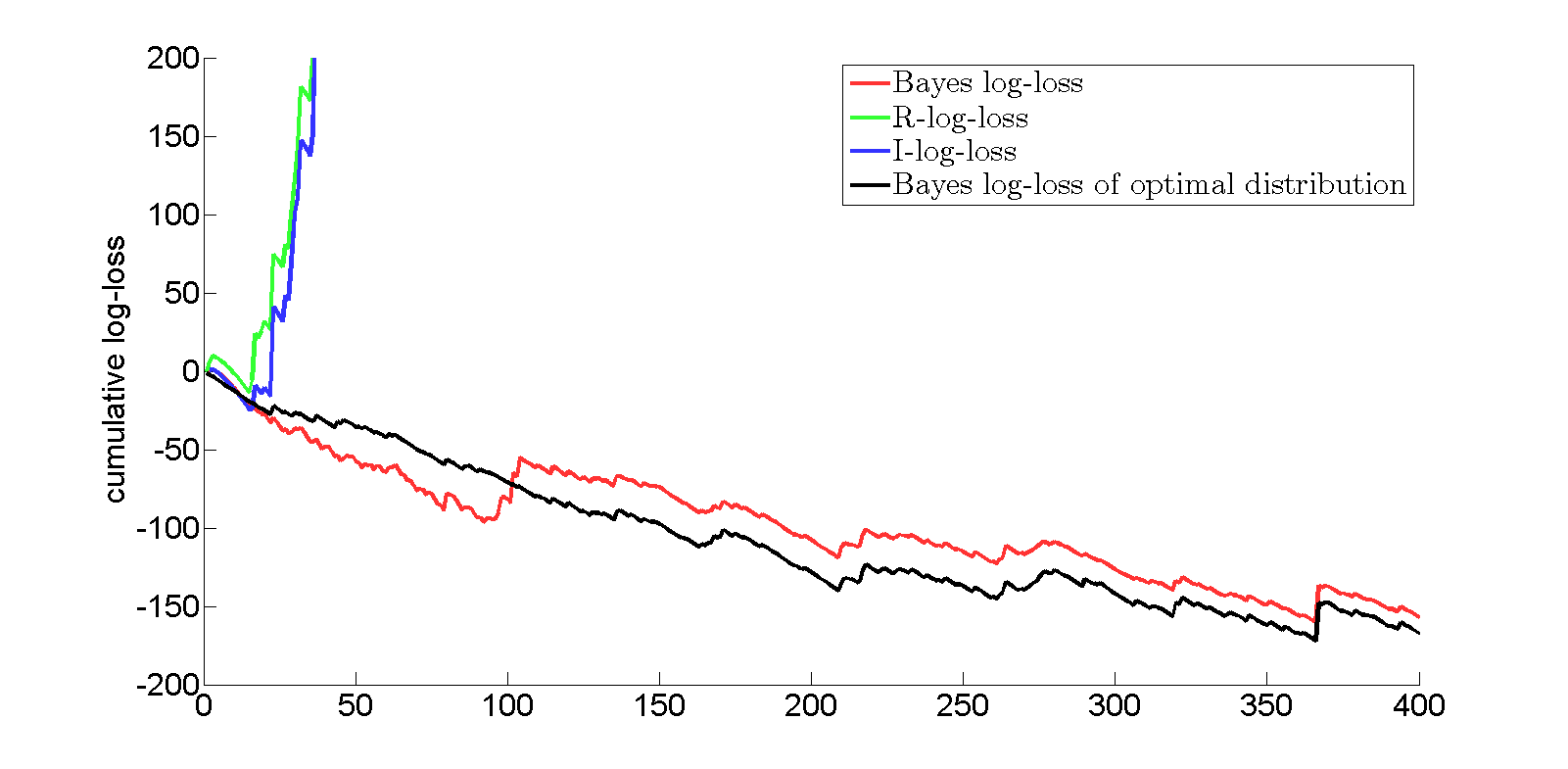}}
\caption{Cumulative standard, $R$-, and $I$-log-Loss as defined in (\ref{eq:kibbeling}) and (\ref{eq:jibbeling}) respectively of standard Bayesian prediction ($\eta =
  1$) on a single
  run for the model-averaging experiment of Figure~\ref{fig:mainexperimenta}.  
  We clearly see that standard Bayes achieves {\em hypercompression},
  being better than the best single model. And, as predicted by
  theory, randomized Bayes is never better than standard
  Bayes, whose curve has negative slope because the densities of outcomes are $> 1$ on average.\label{fig:hypercompressiona} }
\end{figure}
\begin{figure}{\hspace*{0.12\textwidth}
\includegraphics[width=0.75\textwidth]{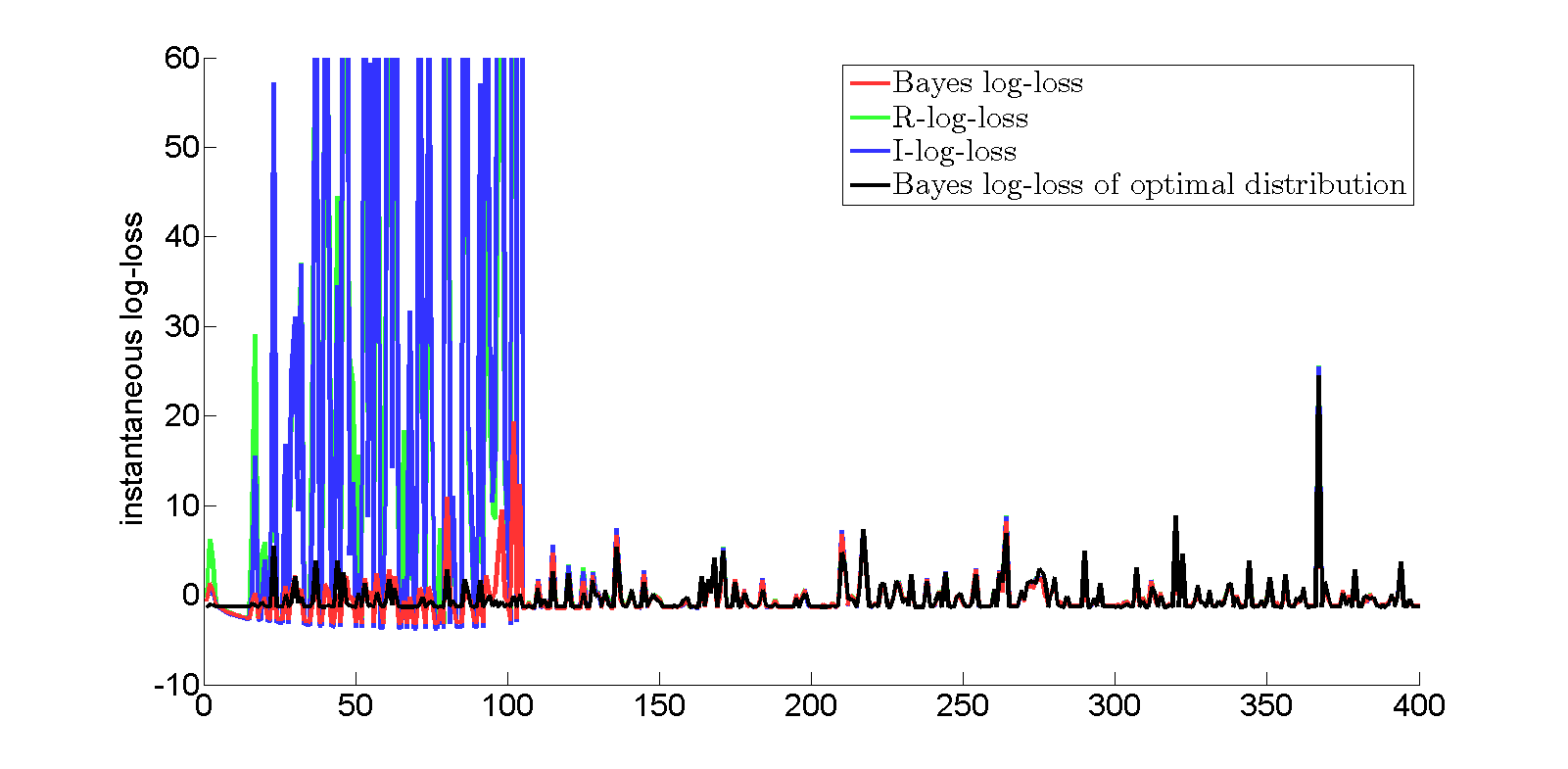}}
\caption{Instantaneous standard, $R$- and $I$-log-Loss of standard Bayesian prediction for the run depicted in Figure~\ref{fig:hypercompressiona}.
\label{fig:hypercompressionb} }
\end{figure}

\begin{figure}
\hspace*{0.15\textwidth}
\includegraphics[width=0.7\textwidth]{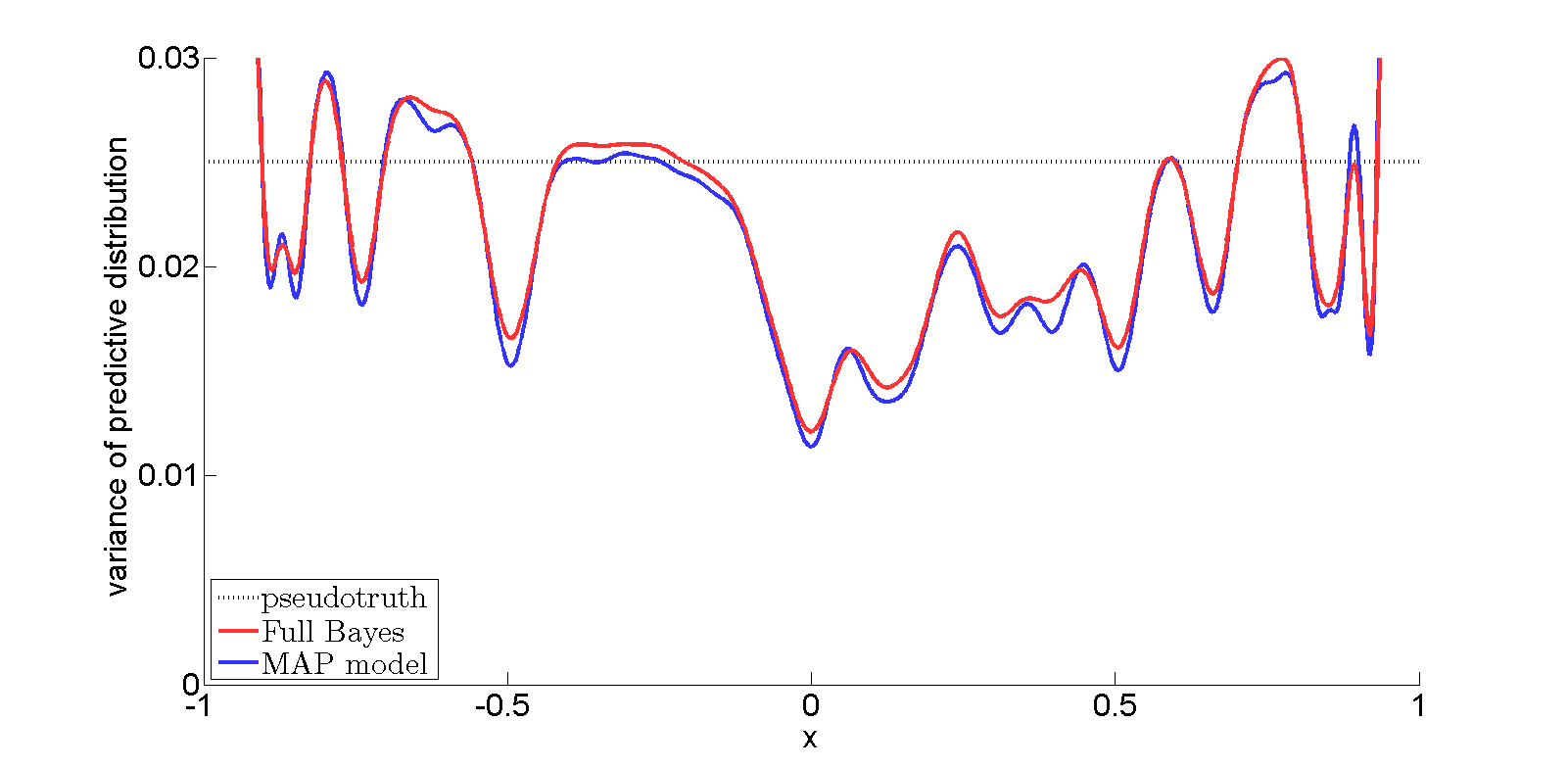}
\caption{Variance of standard Bayes predictive distribution
  conditioned on a new input $S$ as a function of $S$ after 50
  examples for the model-wrong experiment
  (Figure~\ref{fig:mainexperimenta}), shown both for the predictive
  distribution based on the full, model-averaging posterior and for
  the posterior conditioned on the MAP model
  $\M_{\breve{\nc}_{\text{map}}}$. For both posteriors, the posterior
  mean of $Y$ is incorrect for $x \neq 0$, yet $\pbayes(Y \mid Z^{50},X)$
  still achieves small risk because of its small variance at $X=0$.
  \label{fig:polynomialpredvar}}
\end{figure}

\paragraph{\bf How Hypercompression arises in Regression}
Figure~\ref{fig:polynomialpredvar} gives some clues as to how
hypercompression is achieved: it shows the variance of the predictive
distribution $\pbayes(\cdot \mid Z^{50})$ as a function of $S \in
[-1,1]$ for the polynomial example of Figure~\ref{fig:polynomial} in
the introduction, at sample size $n=50$, where hypercompression takes
place. Figure~\ref{fig:polynomial} gave the posterior mean (regression
function) at $n = 100$; the function at $n=50$ looks similar,
correctly having mean $0$ at $S=0$ but, incorrectly, mean far from $0$
at most other $S$. The predictive distribution conditioned on the MAP
model $\M_{\breve{\nc}_{\text{map}(Z^{50})}}$ is a t-distribution with
approximately $\breve{\nc}_{\text{map}(Z^{50})} \approx 50$ degrees of
freedom, which means that it is approximately normal.
Figure~\ref{fig:polynomialpredvar} shows that its variance is {\em
  much\/} smaller than the variance $\tilde{\sigma}^2$ at $S=0$; as a
result, its log-risk conditional on $U=0$ is smaller than that of
$\tilde\theta = (\tilde{\nc},\tilde\beta,\tilde{\sigma}^2)$ by some
large amount $A$. Conditioned at $S \neq 0$, its conditional mean is
off by some amount, and its variance is, on average, slightly (but not
much) smaller than $\tilde{\sigma}^2$, making its conditional log-risk
given $U \neq 0$ larger than that of $\tilde\theta$ by an amount $A'$
where, it turns out, $A'$ is smaller than $A$. Both events $S = 0$ and
$S \neq 0$ happen with probability $1/2$, so that the final,
unconditional log-risk of $\pbayes(\cdot \mid Z^{50})$ is smaller than
that of $\tilde\theta$.

Summarizing, hypercompression occurs because the variance of the
predictive distribution conditioned on past data and a new $X$ is much
smaller than $\tilde{\sigma}^2$ at $X=0$. This suggests that, if
instead of a prior on $\sigma^2$ we use models $\cM_p$ with a
fixed $\sigma^2$, we can only get hypercompression (and correspondingly
bad square-risk behaviour) if $\sigma^2 \ll \tilde{\sigma}^2$, because
the predictive variance based on linear models $\cM_p$ with fixed
variance $\sigma^2$ given $X=x$ is, for all $x$, lower bounded by
$\sigma^2$. Our experiments in Appendix~\ref{sec:mainfixedsigma} confirm that this is indeed what happens. 

\subsection{Explanation III: The Mixability Gap \& The Bayesian
  Belief in Concentration}
\label{sec:mixability}
As we indicated at the end of Section~\ref{sec:badmisspecification}, bad misspecification
can occur only if the standard ($\eta = 1$) posterior is {\em
  nonconcentrated}\footnote{Things would simplify if we could say `bad
  misspecification can occur if and only if there is hypercompression', but we do
  not know whether that is the case, see Section~\ref{sec:openproblem}.}.  Intriguingly, by formalizing
`concentration' in the appropriate way, we will now show, under some
conditions on the prior, that a {\em Bayesian a priori always believes
  that the posterior will concentrate very fast}. Thus, if we observe
data $Z^n$, and for many $n' \le n$, the posterior based on $Z^{n'}$
is not concentrated, then we can view this as an indication of bad
misspecification. In the next subsection we will see that  SafeBayes
selects a $\hat\eta \ll 1$ iff we have such nonconcentration at $\eta
= 1$. Thus, SafeBayes can partially be understood as a prior
predictive check, i.e.\ a test whether the assumptions implied by the
prior actually hold on the data \citep{box1980sampling}. 

\paragraph{The Mixability Gap}
We express posterior nonconcentration in terms of the {\em mixability
  gap\/} \citep{Grunwald12,RooijEGK14}. In this section we only  consider the  special case of $\eta = 1$ (standard Bayes), for which the mixability gap
$\delta_{i}$ is defined as the difference between $1$-$R$-log-loss
(\ref{eq:kibbeling}) and standard log-loss as obtained by predicting
with the posterior predictive, at sample size $i$:
\begin{align}\label{eq:mixgap}
\delta_{i} & :=  \Exp_{\theta \sim \Pi \mid z^{i-1} } 
\left[ - \log \dens(y_i \mid x_i,\theta) \right] 
- \left( - \log \Exp_{\theta \sim \Pi \mid z^{i-1} }  [\dens(y_i \mid x_i,\theta)]
\right) \nonumber \\  &\; = 
\Exp_{\theta \sim \Pi \mid z^{i-1} } 
\left[- \log \dens_{\theta}(y_i \mid x_i )\right] 
- \left( - \log \pbayes(y_i \mid x_i,  z^{i-1}) \right),
\end{align}
Straightforward application of Jensen's inequality as in
(\ref{eq:jensen}) gives that $\delta_{i} \geq 0$. $\delta_{i}$, which depends on $z_1, \ldots,
z_{i}$, is a measure of the posterior's
concentratedness at sample size $i$ when used to predict $y_i$ given
$x_i$: it is small if  $\dens_{\theta}(y_i \mid
x_i)$ does not vary much among the $\theta$ that have substantial
$\eta$-posterior mass; by strict convexity of $- \log$, it is $0$ iff
there exists a set $\Theta_0$ with $\Pi(\Theta_0 \mid Z^{i-1})
=1$ such that for all $\theta, \theta' \in \Theta_0$, $f_{\theta}(y_i
\mid x_i) = f_{\theta'}(y_i \mid x_i)$.

We set the {\em cumulative mixability gap\/} to be $\Delta_{n}
:= \sum_{i=1}^n \delta_{i}$. 
\paragraph{The Bayesian Belief in Posterior Concentration}

As a theoretical contribution of this paper, we now show that, under
some conditions on model and prior, if the data are as expected by the
model and prior, then the expected mixability gap goes to $0$ as
$O((\log n)/n)$, and hence a Bayesian automatically a priori believes
that the posterior will concentrate fast. For simplicity we restrict
ourselves to a model $\cM = \{ P_{\theta} : \theta \in \Theta\}$ where
$\Theta$ is countable, and we let all $\theta \in \Theta$ represent a
conditional distribution for $Y$ given $X$, extended to $n$ outcomes
by independence. We let $\pi$ be a probability mass on $\Theta$, and
define the joint Bayesian distribution $\Pi$ on $\Theta \times \cY^n \mid \cX^n$
in the usual way, so that for measurable $\cA \subset \cY^n$, $\Pi(
(\theta^*, \cA) \mid X^n = x^n) = \pi(\theta^*) \cdot P_{\theta^*}(\cA
\mid X^n = x^n)$. The random variable $\theta^*$ refers to the
$\theta$ chosen according to density $\pi$. We will look at the
Bayesian probability distribution of the $\theta^*$-expected
mixability gap, $\bar{\delta}_n := \Exp_{\theta^*}[\delta_n]$.
\begin{theorem}\label{thm:concentration}
Consider a countable model with prior $\Pi$ as above.   
Suppose that the density ratios in $\Theta$ are uniformly bounded,
  i.e.\ there is a $v > 1$ such that for all $x,y \in \cX \times \cY$,
  all $\theta, \theta' \in \Theta$, $\dens_{\theta}(y \mid
  x)/\dens_{\theta'}(y \mid x) \leq v$. Suppose that for some $\eta <
  1$ we have $\sum_{\theta} \pi(\theta)^{\eta} < \infty$. Then for every $a > 0$ there are constants $C_0$ and $C_1$ such that, for all $n$,
\begin{equation}
\Pi\left(\  
\bar{\delta}_n
\geq C_0 \cdot \frac{\log n}{n} \ \right)\leq C_1 \cdot \frac{1}{n^a}.
\end{equation}
Moreover, for any $0 < a' \leq 1$, there exist $C_2$ and $C_3$ such that 
\begin{equation}\label{eq:heathrow}
\Pi\left(\  
{\Delta}_n
\geq C_2  \cdot n^{a'} \ \right) \leq C_3 \cdot \frac{(\log n)^2}{n^{a'}},
\end{equation}
i.e.\ the Bayesian believes that the mixability gap will be small on
average and that the cumulative mixability gap will be small with high probability.  
\end{theorem}
Thus, with high probability, $\Delta_{n}$ grows only
polylogarithmically, even though it is the difference of two
quantities that are typically linear in $n$. This means that observing
a large value of $\Delta_n$ strongly indicates misspecification.
\begin{ownquote}
  We hasten to add that the regularity conditions for
  Theorem~\ref{thm:concentration} do {\em not\/} hold in the
  regression problem we study in this paper; the theorem is merely
  meant to show that $\Delta_{n}$ is believed to be small in idealized
  circumstances that have been simplified so as to make mathematical
  analysis easier. Note however, that the regularity conditions do not
  constrain $\Theta$ in the most important respect: by allowing
  countably infinite $\Theta$, we can approximate nonparametric models
  arbitrarily well by suitable covers \citep{BarronC91}. In particular
  we do allow sets $\Theta$ for which maximum likelihood methods would
  lead to disastrous overfitting at all sample sizes. Also the
  condition that $\sum \pi(\theta)^{\eta} < \infty$ is standard in
  proving Bayesian and MDL convergence theorems
  \citep{BarronC91,Zhang06density}. In fact, since the constants $C_0$ and $C_1$
  scale logarithmically in $v$, we expect that
  Theorem~\ref{thm:concentration} can be extended to the regression
  setting we are dealing with here as long as all distributions in the
  model have exponentially small tails, using methods similar to those
  in \cite{grunwald2016fast}.
\end{ownquote}

\begin{example}{\bf [Cumulative Nonconcentration can (and will) go
    together with  Momentary Concentration: 
Example~\ref{ex:bernoulli}, Bernoulli, Cont.]}\label{ex:bernoullib}
Consider the first instance of the Bernoulli
Example~\ref{ex:bernoulli} again, where we again look at what happens
if both distributions are equally bad: $\cM = \{ P_{0.2}, P_{0.8} \}$,
whereas $Y_1, Y_2, \ldots$ are i.i.d.~$\sim P_{\theta^*}$ with
$\theta^* = 1/2$.  As we showed in that example, at any given $n$,
with $P_{\theta^*}$-probability at least $0.32$, $\min_{\theta \in
  \{0.2, 0.8 \}} \pi(\theta \mid Y^n) \approx 2^{- \sqrt{n}/2}$: the
posterior puts almost all mass on one $\theta$. Lemma 6 of
\cite{ErvenGKR11} shows that in such cases $\delta_{n}$ is small; in
this particular case, $\delta_{n} \leq 2 (e-2) \min_{\theta \in \{0.2,
  0.8 \}} \pi(\theta \mid Y^n) \approx 1.42 \cdot 2^{-
  \sqrt{n}/2}$. Thus, the posterior {\em looks\/} exceedingly
concentrated at time $n$, with nonnegligible probability (this
unwarranted confidence is a simplified version of what was called the
{\em fair balance paradox\/} by \cite{Yang07paradox}, who conjectured it is
the underlying reason for the problem of `overconfident posteriors' in
Bayesian phylogenetic tree inference).  However, Safe Bayes detects
misspecification by looking at {\em cumulative\/} concentration,
i.e.\ the sum of the $\delta$'s: $L$ as in (\ref{eq:BernoulliLL}) can
be interpreted as a random walk on ${\mathbb Z}$ starting at the
origin, with equal probabilities to move to the left and to the
right. By the central limit theorem, the random walk crosses the
origin at time $n$ with probability about $1/\sqrt{n \pi /2} =
\tilde{O}(n^{-1/2})$, so that we may conjecture that, with high
probability, it crosses the origin $\tilde{O}(n \cdot n^{- 1/2}) =
\tilde{O}(n^{1/2})$ times. Each time it crosses the origin, the
posterior is uniform and hence as nonconcentrated as it can be, and
$\Delta_{n}$ is increased by at least a fixed constant. One would
therefore expect (under the `true' $\theta^*$) that $\Delta_{n}$ is of
order $\sqrt{n}$, which by Theorem~\ref{thm:concentration} is much
larger than a Bayesian a priori expect it to be --- the model fails
the `prior predictive check'.\footnote{This heuristic argument can
  actually be formalized: if data are i.i.d.~Bernoulli$(1/2)$, then
  the expected regret for every absolute loss predictor is of order
  $\tilde{O}(n^{1/2})$ \citep{CesaBianchiL06}, which implies, via the
  connections between regret and $\Delta_{n}$ given by
  \cite{RooijEGK14}, that $\Delta_{n}$ must also be of order
  $n^{1/2}$; we omit further details. }\end{example}

\subsection{How Safe Bayes Works}
\label{sec:howitworks}
In its simplest form, the in-model fixed variance case, SafeBayes
finds the $\hat\eta$ that minimizes cumulative square-loss on the
sample and thus can simply be understood as a pragmatic attempt to
find a $\hat{\eta}$ that achieves small risk. However, the other
versions of SafeBayes do not have such an easy interpretation. To explain them further, we need to generalize the notion of mixability gap in
terms of the {\em $\eta'$-flattened\/} $\eta$-generalized Bayesian
predictive density. The latter is defined, for $\eta, \eta' \leq 1$,
as:
\begin{equation}\label{eq:flattened}
\pbayes(y_i \mid x_i, z^{i-1}, \langle\eta'\rangle ; \eta) 
:= \left( \Exp_{\theta \sim \Pi \mid z^{i-1}, \eta} \left[
    \dens_{\theta}^{\eta'}(y_i \mid x_i) \right] \right)^{1/\eta'}.
\end{equation}
By Jensen's inequality, for any $\eta' \leq 1$, any $(x_i,y_i)$, we
have $\pbayes(y_i \mid x_i, z^{i-1}, \langle\eta'\rangle; \eta) \leq
\pbayes(y_i \mid x_i,z^{i-1} ,\eta)$. Indeed, intentionally, $\pbayes
(\cdot \mid \langle\eta'\rangle; \eta)$ is a `defective' density in
the sense that $\int_{\reals} \pbayes(y \mid x_i, z^{i-1},
\langle\eta'\rangle; \eta) d y < 1$.  The
log-loss achieved by $\eta$-generalized, $\eta'$-flattened Bayesian
prediction is called {\em $(\eta,\eta')$-mix-loss\/} from now on,
following terminology from \cite{RooijEGK14}. For $0 < \eta \leq
\eta'\leq 1$, the {\em mixability gap\/} $\delta_{i,\eta,\eta'}$ is
defined as the difference between the $\eta$-$R$-log-loss and the
$\eta'$-mix-loss:
\begin{equation}\label{eq:mixgapb}
\delta_{i, \eta, \eta'} := \Exp_{\theta \sim \Pi \mid Z^{i-1}, \eta}
\left[ - \log \dens_{\theta}(Y_i \mid X_i) \right] - \left( - \log
  \pbayes(Y_i \mid X_i, Z^{i-1}; \langle \eta'\rangle; \eta)
\right).
\end{equation} We once again define a cumulative version
$\Delta_{n,\eta,\eta'} = \sum_{i=1}^n \delta_{i,\eta,\eta'}$, and note
that the definitions are compatible with the special cases $\delta_{i}
:= \delta_{i,1,1}$ and $\Delta_{n} := \Delta_{n,1,1}$ defined in the
previous subsection. Now we can rewrite the cumulative $R$-log-loss
achieved by Bayes with the $\eta$-generalized posterior as
\begin{equation}\label{eq:mixabilitygap}
\sum_{i=1}^n \Exp_{\theta \sim \Pi \mid z^{i-1}, \eta} \left[ - \log \dens_{\theta}(y_i \mid x_i) \right] 
= \Delta_{n,\eta,\eta'} + \text{CML}_{n,\eta,\eta'},
\end{equation}
where 
$$\text{CML}_{n,\eta,\eta'} = \left( \sum_{i=1}^n - \log \pbayes(y_i
  \mid x_i, z^{i-1}, \langle \eta'\rangle ; \eta) \right)
$$ is the cumulative $(\eta,\eta')$-mix-loss. 
(\ref{eq:mixabilitygap}) holds for all $0 < \eta \leq \eta' \leq 1$.
Consider first $\eta'= 1$. As was seen, if $\Delta_{n,1,1}$ is large,
then this indicates potential bad misspecification.  But
(\ref{eq:mixabilitygap}) still holds for smaller $\eta'< 1$; by
Jensen's inequality, for any fixed $\eta$, decreasing $\eta'$ will
make $\Delta_{n,\eta,\eta'}$ smaller as well.  Indeed, for any fixed $P^*$, defining 
$$
\bar{\delta}_{\eta'} := 
\sup_{W} \; \Exp_{X,Y \sim P^*}  \left[ \;   \Exp_{\theta \sim W}[ - \log \dens_{\theta}(Y \mid X)] 
-  \left( - \frac{1}{\eta'}\log \Exp_{\theta \sim W}[\dens_{\theta}(Y \mid X)^{\eta'}]
\right) \right],
$$
where the supremum is over {\em all\/} distributions on $\Theta$, we have
$$
\lim_{\eta'\downarrow 0} \bar{\delta}_{\eta'} = 0,
$$
so we have an upper bound on the expectation of
$\Delta_{n,\eta,\eta'}$ independent of the actual data that, for small
enough $\eta'$, will become negligibly small. But the left-hand side of
(\ref{eq:mixabilitygap}) does not depend on $\eta'$, so if, by
decreasing $\eta'$, we decrease $\Delta_{n,\eta,\eta'}$,
$\text{CML}_{n,\eta,\eta'}$ must increase by the same amount --- so as
yet we have gained nothing.  Indeed, not surprisingly, Barron's bound
does not hold any more for $\text{CML}_{n,\eta,\eta'}$ with $\eta = 1$
and $\eta'< 1$ (and in general, it does not hold for $\eta,\eta'$
whenever $\eta'< \eta$). {\em But}, it turns out, a version of
Barron's bound still holds for $\text{CML}_{n,\eta',\eta'}$, for all
$\eta'> 0$: the cumulative log-risk of $\eta'$-flattened,
$\eta'$-generalized Bayes is still guaranteed to be within a small
$\text{\sc red}_n$ of the cumulative log-risk of $\tilde\theta$,
although $\text{\sc red}_n$ does monotonically increase as $\eta'$
gets smaller --- simply because the prior becomes more important
relative to the data (standard results in learning theory show that
$\text{CML}_{n,\eta,\eta}$ is monotonically decreasing in $\eta$, and
can be upper bounded as $O(1/\eta)$; see e.g.\ \citep[Lemma
1]{RooijEGK14}. Thus, it makes sense to consider the special case
$\eta'= \eta$, and think of SafeBayes as finding the $\eta$ minimizing
\begin{equation}\label{eq:mixabilitygapb}
\sum_{i=1}^n \Exp_{\theta \sim \Pi \mid z^{i-1}, \eta} \left[ - \log \dens_{\theta}(y_i \mid x_i) \right] 
= \Delta_{n,\eta,\eta} + \text{CML}_{n,\eta,\eta},
\end{equation}
since we have clear interpretations of both terms: the second 
indicates, by Barron's bound, how much worse the $\eta$-generalized
posterior predicts in terms of log-loss compared to the optimal
$\tilde\theta$; the first indicates how much is additionally lost if
one is forced to predict by distributions inside the model. The second
term decreases in $\eta$, the first has an upper bound which increases
in $\eta$. SafeBayes can now be understood as trying to minimize both
terms at the same time.

Now broadly speaking, the central convergence result of
\cite{Grunwald12} states that $\Delta_{n,\eta,\eta}$ will be
`sufficiently small' for all $\eta < 1$, and under some further
conditions even for $\eta = 1$, if the model is correct or convex; and it will
also be `sufficiently small' if the model is incorrect, as long as
$\eta$ is smaller than some `critical' value $\etacrit$ (which may
depend on $n$ though).  Here `sufficiently small' means that it is not
the dominating term in (\ref{eq:mixabilitygapb}). Intuitively, we
would like the $\hat\eta$ determined by SafeBayes to be the largest
$\eta$ that is smaller than $\etacrit$.  \cite{Grunwald12} shows that
Safe Bayes indeed finds such an $\eta$, and that prediction based on
the generalized posterior with this $\eta$ achieves good frequentist
convergence rates.

\begin{figure}[h]{\hspace*{0.125\textwidth}
\includegraphics[width=0.75\textwidth]{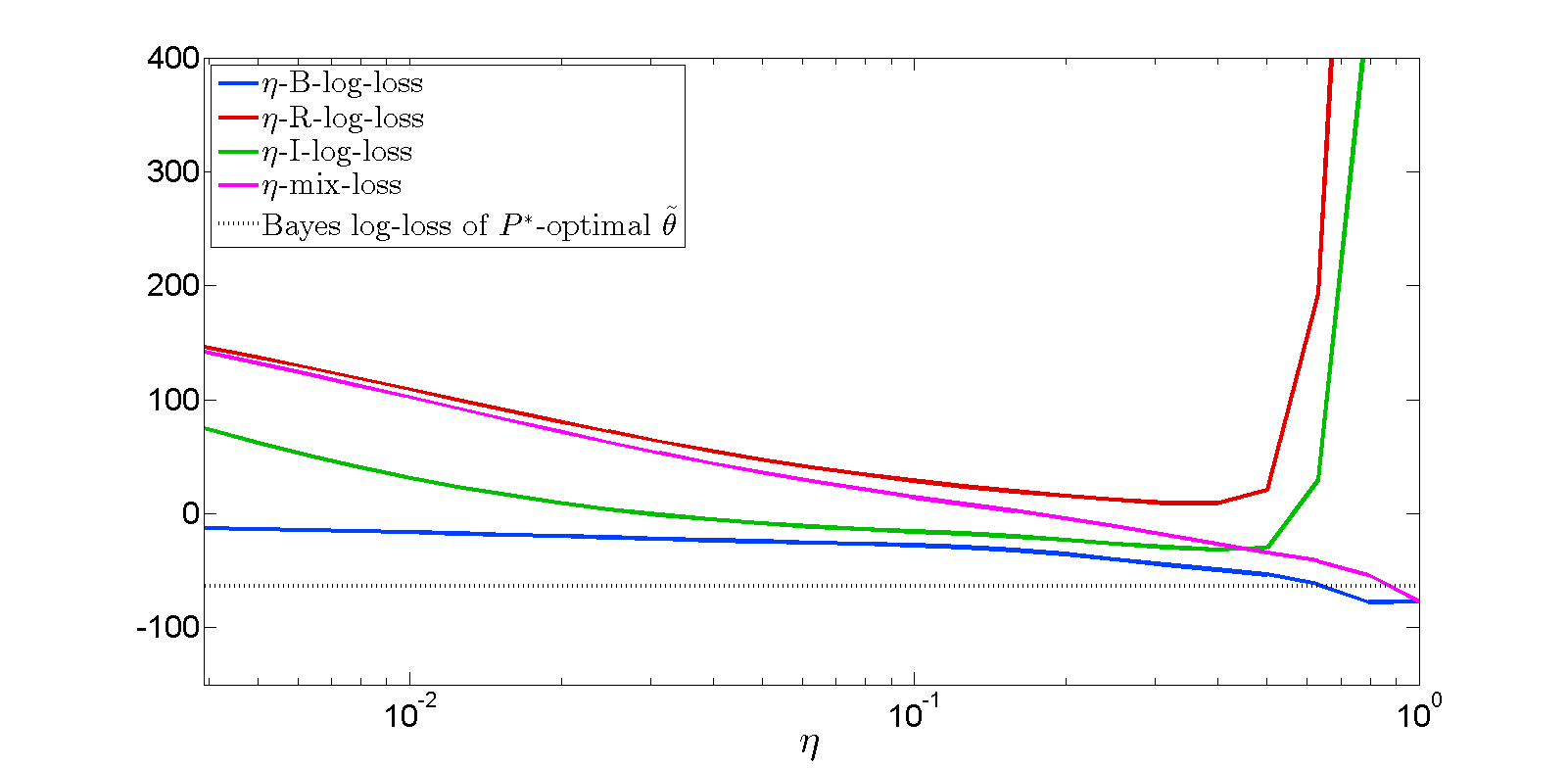}}
\caption{\label{fig:etafunction} Cumulative losses up to sample 100
  (where the posterior has not converged yet) as a function of $\eta$,
  averaged over 30 runs, for the experiment of Figure~\ref{fig:mainexperimenta}.
  $\eta$-B-log-loss is the cumulative log-loss achieved by standard
  Bayes with the $\eta$-generalized posterior.}
\end{figure}

\paragraph{Experimental Illustration:}
Consider the main wrong-model experiment of
Section~\ref{sec:mainexperiments}.
Figure~\ref{fig:etafunction} shows, as a function of $\eta$, in red, the cumulative $\eta$-$R$-log-loss
achieved by Safe Bayes, averaged over 30 runs of
Experiment 1 (Bayesian model averaging with incorrect model) of
Figure~\ref{fig:mainexperimenta}. In each individual run, Safe Bayes picks the
$\hat{\eta}$ minimizing this quantity; we thus get that on most runs,
$\hat\eta$ is close to $0.4$. In contrast to $\eta$-$R$-log-loss, and
as predicted by theory, the $\eta$-mix-loss (in purple) decreases
monotonically and coincides with the standard Bayesian log-loss at
$\eta =1$ and with the $\eta$-$R$-log-loss as $\eta \downarrow 0$. We also see hypercompression again: near $\eta = 1$, both
the Bayesian log-loss and the mix-loss are smaller than the log-loss
achieved by the best $\tilde\theta$ in the model.  At $\eta = 0.5$,
there is a sudden sharp rise in $\Delta_{n,\eta,\eta}$ (the difference
between the red and purple curves). We can think of Safe Bayes as
trying to identify this `critical' $\etacrit$.

\begin{ownquote}
  Theorem~\ref{thm:concentration} shows that, if both model and prior
  are well-specified, then the Bayesian posterior cumulatively
  concentrates in a very strong sense. More generally, if the model is
  correct but also if there is `benign' misspecification, then, under
  some conditions on the prior, by the results of \cite{Grunwald12},
  the Bayesian posterior eventually cumulatively concentrates at
  $\eta = 1$. One might thus be tempted to interpret $\etacrit$ (the
  learning rate which SafeBayes tries to learn) as `largest learning
  rate at which the posterior cumulatively concentrates'. However,
  this interpretation works only if $\etacrit = 1$. If $\etacrit< 1$,
  we can only show that, for every $\eta < \etacrit$,
  $\Delta_{n,\eta,\eta}$ is small; true cumulative concentration would
  instead mean that $\Delta_{n,\eta,1}$ is small for such $\eta$ (note
  we must have $\Delta_{n,\eta,\eta} \leq \Delta_{n,\eta,1}$ by
  Jensen). The figure shows that $\Delta_{n,\eta,1}$ (the difference
  between the red and blue curve) may indeed be large even at small
  $\eta$.  A better interpretation is that, for every fixed $\eta$, with
  decreasing $\eta'$, the geometry of the $(\eta,\eta')$-mix-loss
  changes, so that the loss difference between the mix loss and the
  $R$-log-loss obtained by randomization gets smaller. By then further using
  the generalized posterior for the same $\eta'$, we guarantee that a
  version of Barron's bound  holds for
  the $(\eta',\eta')$-mix-loss. 
\end{ownquote}
\begin{ownquote}
{\bf Replacing $R$- by $I$-loss}
Although the proofs of \cite{Grunwald12} are optimized for
$R$-SafeBayes, the same story as above can be told for any fixed
transformation from the posterior to a possibly randomized prediction,
i.e.\ anything of the form (\ref{eq:haring}); in particular for the
most extreme transformation where we replace the posterior predictive
by the distribution indexed by the posterior mean parameters so that
instead of $R$-SafeBayes we end up with $I$-SafeBayes.  In fact, the
importance of the distinction between `in-model' and `out-model'
prediction under model misspecification has been emphasized before
\citep{Grunwald07,BarronH98,kotlowski2010following}. In general,
although we do not know how to exploit this intuition to strengthen
the convergence proofs of \cite{Grunwald12}, it seems more natural to
replace the randomized predictions by deterministic, in-model
predictions.
\end{ownquote}
\section{Discussion, Open Problems and Conclusion}
\label{sec:discussion}
\begin{quote}{\small 
``If a subjective distribution $\Pi$ attaches probability zero to a
non-ignorable event, and if this event happens, then $\Pi$ must be
treated with suspicion, and {\em modified\/} or replaced'' (emphasis
added) \\
--- A.P. Dawid (\citeyear{Dawid82}).}
\end{quote}
\begin{quote}{\small 
``Some models are obviously wrong, yet evidently useful'' \\
---  (very freely
paraphrasing \cite{Box79}).}
\end{quote}
We already discussed the theoretical significance of the inconsistency
result in the introduction. Extensive further discussion on Bayesian
inference under misspecification is given by \cite{walker2013bayesian}
and \cite{GrunwaldL07}. For us, it remains to discuss the place of
both the inconsistency result and our solution in Bayesian
methodology.

Following the well-known Bayesian statisticians
\cite{box1980sampling}, \cite{Good83}, \cite{Dawid82,Dawid04} and
\cite{Gelman04} (see also \cite{GelmanS12}), we take the stance that
model checking is a crucial part of successful Bayesian practice. When
there is a large discrepancy between a model's predictions and actual
observations, it is not merely sufficient to keep gathering data and
update one's posterior: something more radical is needed. In many such
cases, the right thing to do is to go back to the drawing board and
try to devise a more realistic model. However, we think this story is
incomplete: in machine learning and pattern recognition, one often
encounters situations in which the model employed is {\em obviously\/}
wrong in some respects, yet there is a model instantiation (parameter
vector) that is {\em pretty adequate\/} for the specific prediction
task one is interested in. Examples of such
obviously-wrong-yet-pretty-adequate models are, like in this paper,
assuming homoskedasticity in linear regression when the goal is to
approximate the true regression function and the true noise is
heteroskedastic\footnote{As long as, as in this paper, the tails of
  the conditional distribution of $Y$ given $X=x$ are sub-Gaussian,
  for each $x$; if they are not, there may be real outliers and then
  one cannot say that the model is `pretty adequate' any more.}, but
also the use of $N$-grams in language modeling (is the probability of
a word given the previous three words really independent of everything
that was said earlier?), logistic regression in e.g.\ spam filtering,
and every single successful data compression method that we know of
(see {\em Bayes and Gzip\/} \cite[Chapter 17, page 537]{Grunwald07}).
The difference with the more standard statistical (be it Bayesian or
frequentist) mode of reasoning is eloquently described in Breiman's
(\citeyear{Breiman01}) {\em the two cultures}\footnote{The `two
  cultures' does {\em not\/} refer to the Bayesian-frequentist divide,
  but to the modeling vs.  prediction-divide. We certainly do not take
  the extreme view that statisticians should {\em only\/} be
  interested in prediction tasks such as classification and
  square-error prediction rather than density estimation and testing;
  our point is merely that in some cases, the goal of inference is
  clearly defined (it could be classification, but it could also be
  determination whether some random variables are (conditionally)
  (in)dependent), whereas part of our model is unavoidably
  misspecified; and in such cases, one may want to use a generalized
  form of Bayesian inference. }. Bayesian inference is among the most
successful methods currently used in the
obviously-wrong-yet-pretty-adequate-situation (to witness,
state-of-the-art data compression methods such as
Context-Tree-Weighting \cite{WillemsST95} have a Bayesian
interpretation). Yet the present paper shows that there is a danger:
even {\em if\/} the employed model is pretty adequate (in the sense of
containing a pretty good predictor), the Bayesian machinery might not
be able to find it. The Safe Bayesian algorithm can thus be viewed as
an attempt to provide an alternative for the {\em data-analysis
  cycle\/} \citep{GelmanS12} to this, in some sense, less ambitious
setting: just like in the standard cycle, we do a model check, albeit
a very specific one: we check whether there is `cumulative
concentration of the posterior' (see Section~\ref{sec:mixability}). If
there is not, we know that we may not be learning to predict as well
as the best predictor in our model, so we {\em modify\/} our
posterior. Not in the strong sense of `going back to the drawing
board', but in the much weaker sense of making the learning rate
smaller --- we cannot hope that our model of reality has improved,
because we still employ the same model --- but we can now guarantee
that we are doing the best we can with our given model, something
which may be enough for the task at hand and which, as our experiments
show, cannot always be achieved with standard Bayes.

\paragraph{Benign vs.\ Bad Misspecification}
One might argue that the example of this paper is rather extreme, and
that in practical situations, choosing a learning rate different from
$1$ may never be a useful thing to do.  A crucial point here is that
one can have `benign' and `bad' misspecification
(Section~\ref{sec:badmisspecification}). Under benign
misspecification, standard Bayes with $\eta = 1$ will behave nicely
under weak assumptions on the prior.  While in our particular example,
after `eyeballing' the data one would probably have chosen a
different, less misspecified model, it may be the case that `bad'
misspecification (as in Figure~\ref{fig:badmix}) also occurs, at least
to some extent, in general, real-world data and is then not so easily
spotted. Since we simply do not know whether such situations occur in
practice, to be on the safe side, it seems desirable to have a theory
about when we can get away with using standard Bayesian inference for
a given prediction task even if the model is wrong, and how we can
still use it with little modification if there is bad
misspecification. Our work (esp. the theoretical counterpart to this
paper \citep{grunwald2016fast}) is a first step in this direction.
\paragraph{Towards a Theory of Bayesian Inference under
  Misspecification}
What we have in mind is a theory of Bayesian inference under
misspecification, in which the {\em goal\/} of learning plays a
crucial role. The standard Bayesian approach is very ambitious: it can
be used to solve every conceivable type of prediction or inference
task. Every such task can be encoded as a loss or utility function,
and, given the data and the prior, one merely has to calculate the
posterior, and then makes an optimal decision by taking the act that
minimizes expected loss or maximizes expected utility according to the
posterior. Crucially, one uses the same posterior, independently of
the utility function at hand, implying that one believes that one's
own beliefs are correct {\em in every possible respect}. We envision a
more modest approach, in which one acknowledges that one's beliefs are
only adequate in some respects, not in others; how one proceeds then
depends on how one's model and loss function interact. For example, if
one is interested in data-compression then, this problem being
essentially equivalent to cumulative log-loss prediction, by Barron's
(\citeyear{Barron98}) bound one can simply use the standard $(\eta =
1)$ Bayesian predictive distribution --- even under misspecification,
this will guarantee that one predicts (at least!) as well as one could
with the best element of one's model. If, on the other hand, one is
interested in any of the KL-associated inference tasks (for linear
models, these are square-loss and reliability,
Section~\ref{sec:optimalb}), then using $\eta=1$ is not sufficient
anymore, and one may have to learn $\eta$ from the data, e.g.\ in the
Safe Bayesian manner. Finally, if we are interested in an inference
task that is not KL-associated under our model (i.e., a model instance
can be good in the KL sense but bad in the task of interest), then a
more radical step is needed: either go back to the drawing board and
design a new model after all; or perhaps, the model can be changed in
a more pragmatic way so that, for the right $\eta$, $\eta$-generalized
Bayes once again will find the best predictor for the task at hand.
Let us outline such a procedure for the case that the inference talk
is simply prediction under some loss function $\ell: \cY \times \hat\cY
\rightarrow \reals$. In this case, if the $\ell$-risk is not
KL-associated this simply means that, for some data, the log likelihood is not a
monotonic function of the loss $\ell$. To get the desired association,
we may associate each conditional distribution $P_{\theta}(Y \mid X)$
in the model with its associated Bayes act $\delta_{\theta}$:
$\delta_{\theta}(x)$ is defined as the act $\hat{y} \in \hat{\cY}$
which minimizes $P_{\theta} \mid X=x$-expected loss $E_{Y \sim
  P_{\theta} \mid X=x}[\ell(y,\hat{y})]$.  We can then define a new
set of densities
\begin{equation}\label{eq:sinterklaas}
\dens^{\text{\sc new}}_{\theta,\gamma}(y \mid x) = \frac{1}{Z(\gamma)} e^{- \gamma \ell(y,\delta_\theta(x))},
\end{equation}
and perform (generalized) Bayesian inference based on these. Note that
this effectively replaces, for each $\theta$, the full likelihood by a
`likelihood' in which some information has been lost, and is thus
reminiscent of what is done in {\em pseudo-likelihood\/}
\citep{besag1975statistical} {\em substitution likelihood\/}
\citep{Jeffreys61,dunson2005approximate}, or {\em rank-based
  likelihood\/} \citep{gu2009bayesian} approaches (as a Bayesian, one
may not want to loose information, but whether this still applies in
nonparametric problems \citep{RobinsW00} let alone under
misspecification \citep{GrunwaldH04} is up to debate).

(\ref{eq:sinterklaas}) can be made precise in two ways: either one
just sets $\gamma$ and $Z(\gamma)$ to 1, and allows the
$\dens^{\text{\sc new}}_{\theta}$ to be pseudo-densities, not
necessarily integrating to $1$ for each $x$. This is a standard
approach in learning theory \cite{Zhang06loss,Catoni07}. One could
then learn $\eta$ by, e.g., the basic SafeBayes algorithm with
$\ell_{\theta}(x,y) := \ell(y,\delta_{\theta}(x))$ instead of
log-loss. Or, one could define $Z(\gamma)$ so that the densities
normalize (how to achieve this if $\int_y e^{- \gamma
  \ell(y,\delta_\theta(x))} dy$ depends on $x$ is explained by
\cite{Grunwald08b}) and put a prior on $\gamma$ as well (for linear
models, this is akin to putting a prior on the variance). This will
make the loss $\ell$ KL-associated and the KL-optimal $\tilde\theta$
will also have the reliability property, see again
\cite{Grunwald08b} for details. In this case we will get, with $z_i = (x_i,y_i)$,
$\ell_{\theta}(z_i) := \ell(y_i,\delta_{\theta}(x_i))$, and using a
prior on $\Theta$ and the scaling parameter $\gamma$, that the
$\eta$-generalized posterior becomes
\begin{equation}\label{eq:alin99}
\pi(\theta,\gamma \mid z^n,\eta) \propto \frac{1}{Z(\gamma)^{\eta n}} e^{- \eta
  \gamma \sum_{i=1}^n \ell_{\theta}(z_i) } \cdot \pi(\theta,\gamma).
\end{equation}
This idea was, in essence, already suggested by \cite[Example
5.4]{Grunwald98b} (see also \cite{Grunwald99a}) under the name of {\em
  entropification\/} (however, Gr\"unwald's  papers wrongly suggest that, by
introducing the scale parameter $\gamma$, it would be sufficient to
only consider $\eta = 1$); see also
\citep{lacoste2011approximate,quadrianto2015very}.

Now both `pure' subjective Bayesians and `pure' frequentists might
dismiss this program as severe ad-hockery: the strict Bayesian would
claim that nothing is needed on top of the Bayesian machinery; the
strict frequentist would argue that Bayesian inference was never
designed to `work' under misspecification, so in misspecified
situations it might be better to avoid Bayesian methods altogether
rather than trying to `repair' them. We strongly disagree with both
types of purism, the reason being the ever-increasing number of
successful applications of Bayesian methods in machine learning in
situations in which models are obviously wrong. We would like to
challenge the pure subjective Bayesian to explain this success, given
that the statistician is using a priori distributions that reflect
beliefs which she knows to be false, and are thus not really her
beliefs.  We would like to challenge the pure frequentist to come up
with better, non-Bayesian methods instead. In summary, we would urge
both purists not to throw away the Bayesian baby with the misspecified
bath water!

Moreover, from a prequential \citep{Dawid84}, learning theory (citations see below) and Minimum Description Length (MDL \citep{BarronRY98}) perspective, the extension
from Bayes to SafeBayes is {\em perfectly natural}.   From the prequential
perspective, SafeBayes seeks to find the largest $\eta$ at which the
generalized Bayesian predictions have a predictive interpretation in
terms of the loss of interest rather than the log-loss. The learning theory and  MDL
perspectives are further explained in the next section.
\subsection{Related Work I: Learning Theory and MDL}
\label{sec:related}
\paragraph{Learning Theory} From the learning theory perspective,
generalized Bayesian updating as in (\ref{eq:alin99}) with $Z(\gamma)$
set to $1$ can be seen as the result of a simple regularized loss
minimization procedure (this was probably first noted by
\cite{Williams1980}; see in particular \cite{Zhang06loss}), which
means that it continues to make sense if $\exp(- \gamma
\ell_{\theta})$ as in (\ref{eq:sinterklaas}) does not have a direct
probabilistic interpretation. Variations of such generalized Bayesian
updating are known as ``aggregating algorithm'', ``Hedge'' or
``exponential weights'', and often have good worst-case optimality
properties in nonstochastic settings \citep{Vovk90,CesaBianchiL06} ---
but to get these the learning rate must often be set as small as
$O(1/\sqrt{n})$. Similarly, PAC-Bayesian inference
\citep{Audibert04,Zhang06loss,Catoni07} (for a variation, see
\cite{FreundMS04}) is also based on a posterior of form
(\ref{eq:sinterklaas}) and can achieve minimax optimal rates in
e.g.\ classification problems by choosing an appropriate $\eta$,
usually also very small. From this perspective, SafeBayes can be
understood as trying to find a {\em larger\/} $\eta$ than the
worst-case optimal one, if the data indicate that the situation is not
worst-case and faster learning is possible. Finally,
\cite{bissiri2016general} give a motivation for (\ref{eq:alin99}) (with
$Z(\gamma) \equiv 1$) based on coherence arguments that are more
Bayesian in flavour.
\paragraph{MDL}
Of particular interest is the interpretation of the SafeBayesian
method in terms of the MDL principle for model selection, which views
learning as data compression. When several models for the same data
are available, MDL picks the model that extracts the most `regularity'
from the data, as measured by the minimum number of bits needed to
code the data {\em with the help of the model}.  This is an
interpretation that remains valid even if a model is completely
misspecified \citep{Grunwald07}. The resulting procedure (based on
so-called {\em normalized maximum likelihood\/} codelengths) is
operationally almost identical to Bayes factor model selection. Thus,
it provides a potential answer to the question `what does a high
posterior belief in a model really mean, since one knows all models
under consideration to be incorrect any way?' (asked by, e.g.,
\cite{GelmanS12}): even if all models are wrong, the
information-theoretic MDL interpretation stands. However, our work
implies that there is a serious issue with these NML codes: note that
any distribution $P$ in a model $\cM$ can be mapped to a code (the
{\em Shannon-Fano code\/}) that would be optimal in expectation if data were sampled
from $P$. Now, our work shows that if the data are sampled from some
$P^* \not \in \cM$, then the codes based on Bayesian predictive
distributions can sometimes compress substantially {\em better\/} in
expectation than can be done based on any $P \in \cM$ --- this is the
hypercompression phenomenon of Section~\ref{sec:hypercompression}.
The same thing then holds for the NML codes, which assign almost the
same codelengths as the Bayesian ones. Our work thus invalidates the
interpretation of NML codelengths as `compression with the help of
(and only of!) the model', and suggests that, similarly to in-model
SafeBayes one should design and use `in-model' versions of the NML
codes instead --- codes that are guaranteed not to outperform, at
least in expectation, the code based on the best distribution in the
model.

\subsection{Related Work II: Analysis of Bayesian Behavior under Misspecification}
\paragraph{Consistency Theorems}
The study of consistency and rate of convergence under
misspecification for likelihood-based and specifically Bayesian
methods go back at least to \cite{Berk66}. For recent state-of-the-art
work on likelihood-based, non-Bayesian methods see e.g.\ 
\cite{dumbgen2011approximation} and the very general
\cite{spokoiny2012parametric}. Recent work on Bayesian methods
includes \cite{KleijnV06}, \cite{de2013bayesian} and
\cite{ramamoorthi2013posterior} who obtained results in quite general,
i.i.d.\ nonparametric settings, non-i.i.d.\ settings \citep{Shalizi09},
and more specific settings \citep{sriram2013posterior}; see also
\cite{grunwald2016fast}. Yet, as explicitly remarked by
\cite{de2013bayesian}, the conditions on model and prior needed for
consistency under misspecification are generally stronger than those
needed when the model is correct. Essentially, if the data are i.i.d.\ 
both according to the model and the sampling distribution $P^*$, then
Theorem 1 (in particular its Corollary 1) of \cite{de2013bayesian}
implies the following: if, for all $\epsilon > 0$, the model can be
covered by a finite number of $\epsilon$-Hellinger balls, then the
Bayesian posterior eventually concentrates: for all $\delta, \gamma >
0$, the posterior mass on distributions within Hellinger distance
$\delta$ of the $P_{\tilde\theta}$ that is closest to $P^*$ in KL
divergence will become larger than $1-\gamma$ for all $n$ larger than
some $n_{\gamma}$.  This implies that both in the ridge regression
(finite $\nc$) and in the model averaging experiments (finite
$\pmax$), Bayes eventually `recovers' --- as we indeed see in our
experimental results. However, if $\pmax = \infty$, then the model has
no finite Hellinger cover any more for small enough $\epsilon$ and
indeed the conditions for Theorem 1 of \cite{de2013bayesian} do not
apply any more.  Our results show that in such a case we can indeed
have inconsistency if the model is incorrect.  On the other hand, even
if $\pmax = \infty$, we do have consistency in the setup of our
correct-model experiment for the standard Bayesian posterior, as
follows from the results by \cite{Zhang06density}.  
\paragraph{The Limiting $\eta = 1$}
Like several earlier results \citep{BarronC91,WalkerH02}, Zhang's
consistency results for correct models hold under very weak conditions
for generalized Bayes with any $\eta < 1$, and only under much
stronger conditions for $\eta=1$. Zhang provides an example of
inconsistency-like behavior in the well-specified case with $\eta = 1$
that automatically disappears as soon as one picks $\eta < 1$, leading
\cite{Zhang06density} to claim that in general, generalized Bayesian
methods $(\eta < 1)$ are more stable than standard Bayesian
ones. Zhang's example, and the example of Bayesian model selection
inconsistency in a well-specified model by \cite{CsiszarS00} are
closely related to ours, in that the Bayes predictive distribution for
$\eta = 1$ becomes significantly different from any distribution in
the model (see Figure~\ref{fig:badmix}). In their examples, the
problem is resolved by taking any $\eta < 1$; in our misspecification
case, $\eta$ should even be taken much smaller.

\paragraph{Anomalous Behavior and Modifications of Bayesian Posterior under Misspecification}
Anomalous behavior of Bayesian inference under misspecification was,
of course, observed before, e.g.\ (less dramatically than here) by
\cite{Yang07paradox,muller2013risk} and (as dramatically, but involving a
very artificial model) \cite{GrunwaldL07}. Presumably also related is
the `brittleness' of Bayesian inference that has been observed by
\cite{owhadi2013brittleness}.  Not surprisingly then, we are not the
first to suggest modification of likelihood-based estimators (see
e.g.\ 
\cite{white1982maximum,royall2003interpreting,kotlowski2010following})
and posteriors
\citep{royall2003interpreting,hoff2012bayesian,doucet2012robust,muller2013risk}.
The latter three approaches (that extend the first) employ the
so-called {\em sandwich posterior}, in which the covariance matrix of
the posterior is changed based on a `sandwich formula' involving the
empirical variance; \cite{muller2013risk} provides extensive
explanation and experimentation. Compared to the sandwich approach,
our proposal, besides being applicable in fully nonparametric
contexts, seems substantially more radical. This can be seen from the
regression applications in \cite{muller2013risk}, which involve a
noninformative Jeffreys' prior on the regression coefficient vector
$\beta$. With such a prior (as well as any normal prior scaled by
variance $\sigma^2$), the posterior {\em mean\/} of $\beta$, and thus
also the frequentist square-risk (which only depends on the posterior
mean) remains unaffected by the sandwich modification, so for
square-risk the method would perform like standard Bayes in our
model-wrong experiments. Thus \cite[Section 2.4]{muller2013risk}
demonstrates its usefulness on other loss functions.  Nevertheless,
both the sandwich and the safe Bayesian methods can be thought of as
methods for measuring the spread of a posterior, and it would be
useful to compare the two in detail, both in theory and practice.

\subsection{Future Work and Open Problems}
\label{sec:openproblem}
The results of this paper raise several issues and prompt the following research agenda:
\begin{enumerate}\item The misspecification in our example  would presumably  be easily
  spotted in practice. This raises the question whether `bad'
  misspecification also arises for data sets that occur in practice
  and for which it would not be easily spotted.  Currently, we know
  only of one experiment in this direction: \cite{Jansen13} applied
  the Bayesian Lasso \citep{ParkC08} to several real-world data sets,
  where the $\lambda$ (i.e.\ $1/\eta$) is taken that minimizes the
  cumulative {\em square-loss\/} whereas at the same time $\sigma^2$
  is a free parameter.  Thus it is a hybrid of $I$-square Safe Bayes
  and $I$-log SafeBayes, but equal to neither; the method was
  (somewhat) outperformed by standard Bayes on most data sets tried.
  However, we also tried this hybrid method in the model-wrong
  experiment of this paper and found that it is not competitive with
  either of the two `true' in-model SafeBayes methods either; so the
  experiment does not `really' test SafeBayes; more precise experiments are needed.
\item Our method has one major disadvantage: even if the data do not
  have a natural ordering, the $\hat\eta$ selected by SafeBayes will,
  in general, be order-dependent.  \cite{Grunwald11} suggested a very
  different (and in fact, the first) method to learn $\hat\eta$, that
  does not have this problem. However, it is only applicable to
  countable models, and has no obvious computationally efficient
  implementation, so we do not know whether it has a future. Another
  method that is clearly related to $I$-square SafeBayes is to
  determine $\eta$ using leave-one-out cross-validation based on the
  squared error. This method is also order-independent and behaves
  comparably to $I$-square SafeBayes (Appendix~\ref{sec:mainfixedsigma}), but it is not
  clear how to extend it to general misspecified models, While we show
  in the same appendix that cross-validation based on log-loss of the
  Bayes predictive distribution fails dramatically, it may be that
  cross-validation based on log-loss of the Bayes posterior {\em
    mean\/} would generally work fine, and this method can be applied
  to general misspecified models, not just linear ones. Compared to
  $I$-log-SafeBayes this {\em in-model log-loss cross-validation\/}
  would have the advantage that it is order independent, and the
  disadvantage that it cannot (at least not straightforwardly) be used
  in an online setting and/or for non-i.i.d.\ models. Also, we suspect
  that if the number of models is exponential in the covariates (as in
  variable selection), cross-validation may be prone to overfitting
  whereas SafeBayes would not be, but this is just extrapolation from the well-specified case: it
  would be useful to investigate ``in-model cross-validation'' further.
\item What exactly are relations between the sandwich posterior (see
  above) and our approach?  It would be good to test SafeBayes on the
  data sets used by \cite{muller2013risk}.
\item It would be useful to establish exactly what properties of
  Bayesian updating remain valid for generalized Bayesian updating,
  and what properties do not hold any more. For example, {\em
    telescoping\/} \citep{CesaBianchiL06} holds for the standard
  posterior, for the $\eta$-flattened, $\eta$-generalized posterior,
  but not for the (nonflattened) $\eta$-generalized posterior.
\item As discussed at the end of Section~\ref{sec:howitworks}, the
  final term in (\ref{eq:fronkie}) is lacking in the in-model versions
  of SafeBayes, and this does suggest that they should work better
  than the randomization versions --- the corresponding
  $\Delta_{\eta,\eta}$ is always smaller. Yet we have no theoretical
  results to this end, and our empirical results in this paper confirm
  this to some extent ($R$-square-SafeBayes is not competitive), but not
  fully ($R$-log-SafeBayes is competitive), so more research is needed
  here.
\item As we indicated in Section~\ref{sec:hypercompression}, hypercompression implies
  nonconcentration, but we do not know whether the reverse implication
  holds as well, so we may perhaps  have bad misspecification yet no
  hypercompression. It would give significant insight if we knew whether this indeed could happen. 
\item In light of the discussion underneath (\ref{eq:sinterklaas}),
  one would like to formulate a general theory of substitution
  likelihoods so that likelihoods can be determined based on the
  inference task of interest, so that this task becomes KL-associated,
  for {\em arbitrary\/} prediction tasks. Ideally,
  (\ref{eq:sinterklaas}) and approaches such as pseudo-likelihood and
  rank-based likelihood would all become a special case. If this can
  be done, we would have a truly generalized Bayesian method.
\end{enumerate}

\section*{Acknowledgments}
A large part of this work was done while the authors were visiting UC
San Diego. We would like to think the UCSD CS department for hosting
us. Wouter Koolen, Tim van Erven and Steven de Rooij played a crucial
role in the development of the mixability gap which underlies the Safe
Bayesian algorithm. Many thanks also to Larry Wasserman for useful
feedback and encouragement. This research was supported by NWO VICI
Project 639.073.04.
\DeclareRobustCommand{\VAN}[3]{#3 #1}
\bibliography{bib/peter,bib/MDL,bib/master}
\pagebreak
\appendix

\section{Experiments on Variations of the Prior and the Model}
\label{sec:priorvariations}
Apart from the priors on parameters given the models we used in our
main experiments, we tried several alternative prior distributions,
described in the subsections below. The first subsection describes
experiments with fixed (i.e., a degenerate prior on) $\sigma^2$.
\subsection{Experiments with Fixed $\sigma^2$}
\label{sec:mainfixedsigma}
When models with fixed $\sigma^2$ are used, our two SafeBayes methods
become $R$-square- and $I$-square-SafeBayes, as defined in
Section~\ref{sec:instantiationb}. These also have a direct
interpretation as trying to find the best $\eta$ for predicting with a
square-loss function, as was explained in that section. In this
context, the value $\eta = 1$ has no special status, so we now also
tried values $\eta > 1$ (we did experiment with varying $\eta$ in the
previous varying $\sigma^2$ experiments as well, but there it did not
make any substantial difference in the results).  Specifically, we set
${\cal S}_n$ in the Safe Bayesian algorithm to $\{2^{\kappa_{\max}},
2^{\kappa_{\max}- \kappa_{\text{\sc step}}}, 2^{\kappa_{\max}- 2
  \kappa_{\text{\sc step}}}, \ldots, 2^{- \kappa_{\max}}\}$, with
$\kappa_{\text{\sc step}}= 1/2$ and $\kappa_{\max} = 6$. All priors on
the regression coefficients $\beta$ remain as described in
Section~\ref{sec:preparing}.

\subsubsection{Model Averaging Experiment, Fixed $\sigma^2$}
\begin{figure}[htp]{\hspace*{0.2\textwidth}
\includegraphics[width=0.6\textwidth]{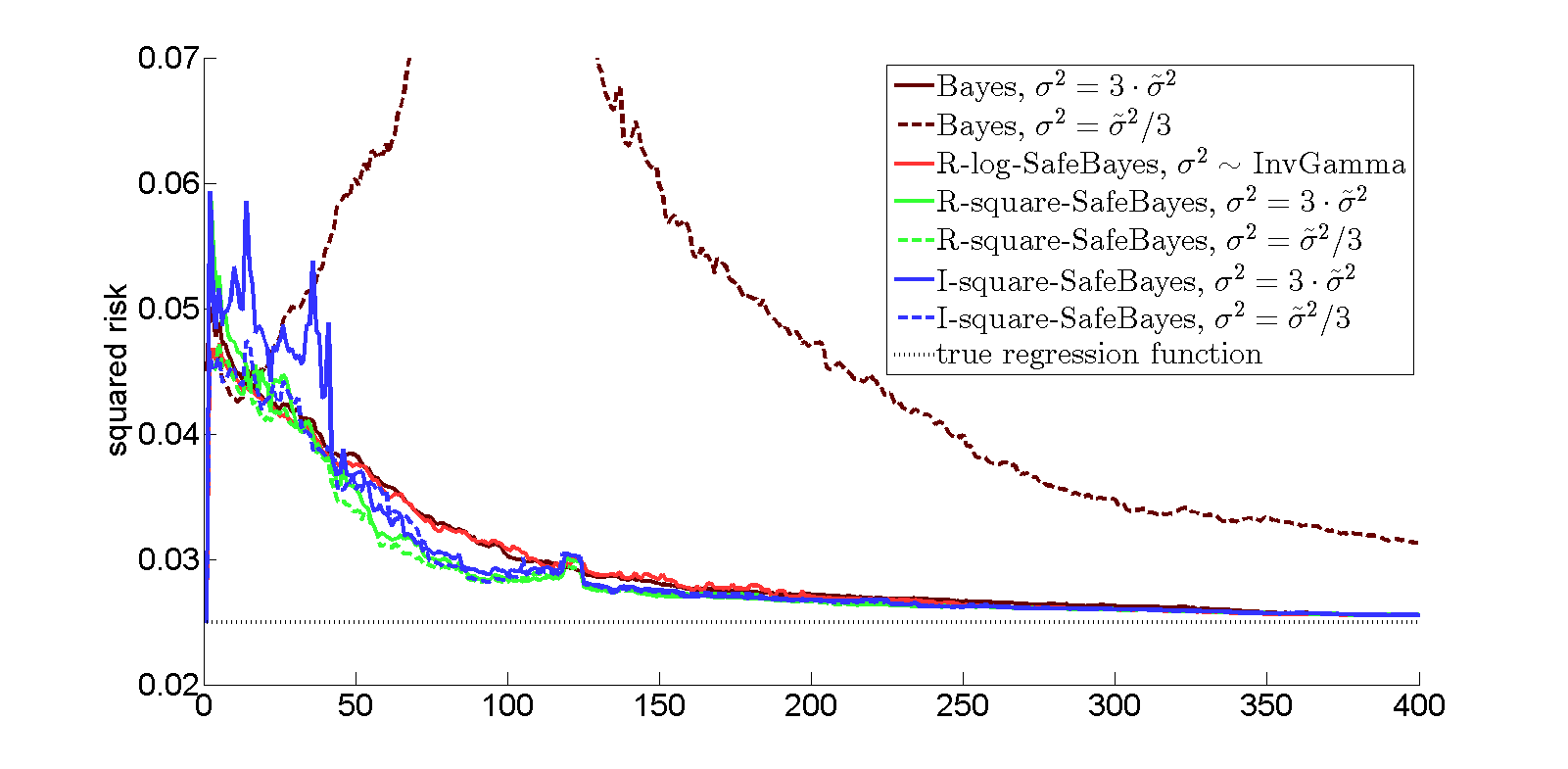} \\
\hspace*{0.2\textwidth}
\includegraphics[width=0.6\textwidth]{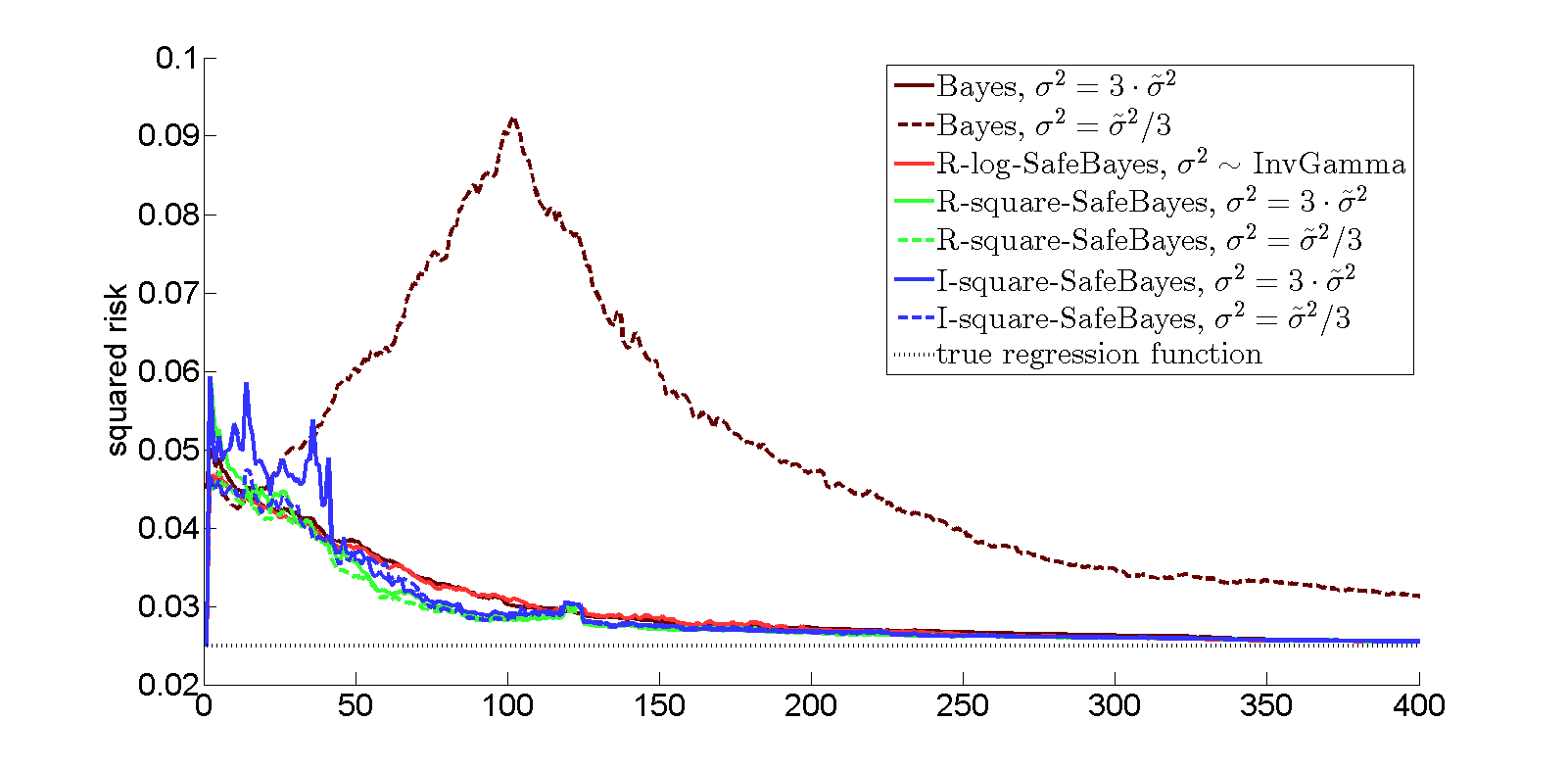} \\
\hspace*{0.2\textwidth}
\includegraphics[width=0.6\textwidth]{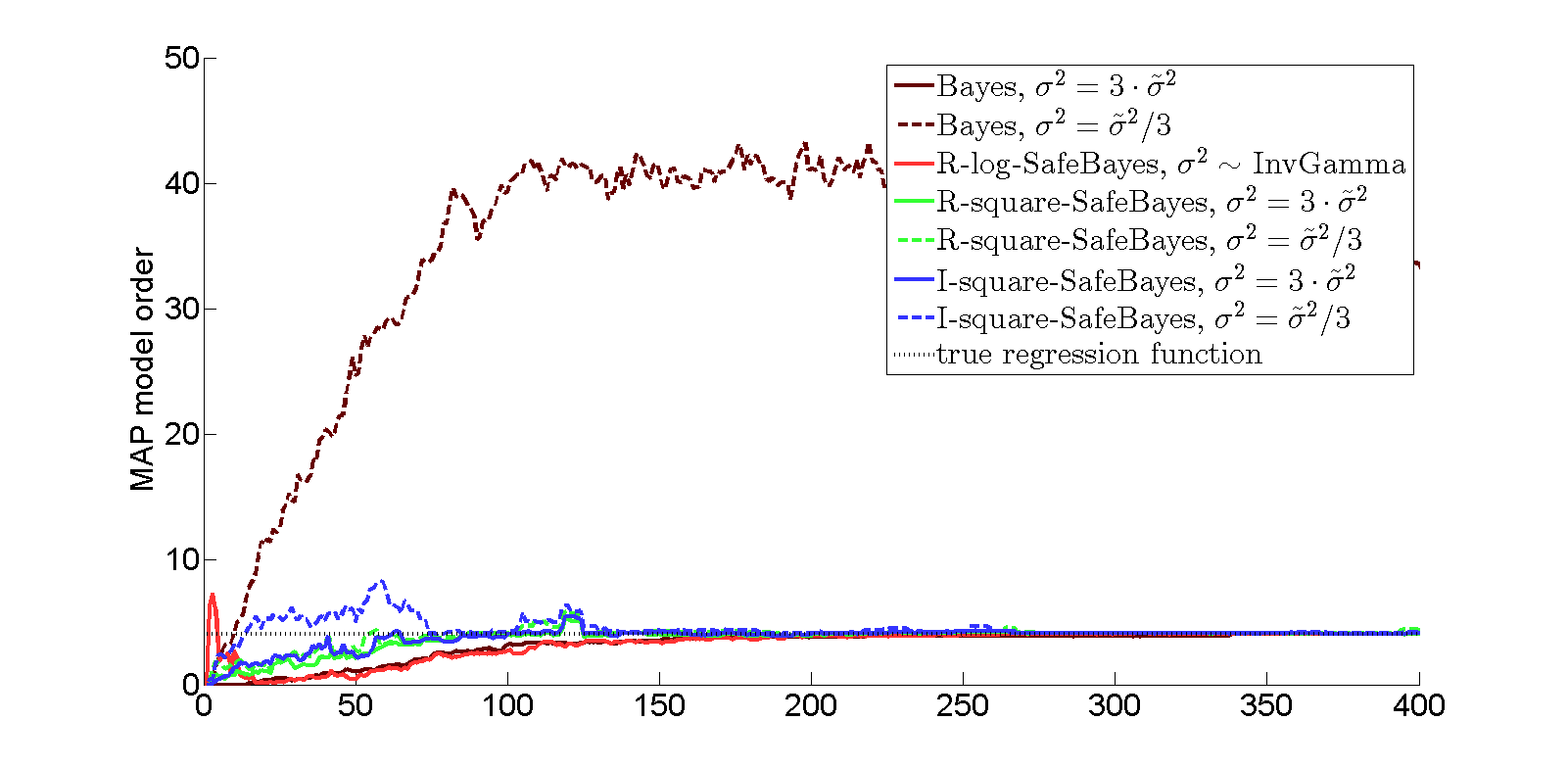} \\
\hspace*{0.2\textwidth}
\includegraphics[width=0.6\textwidth]{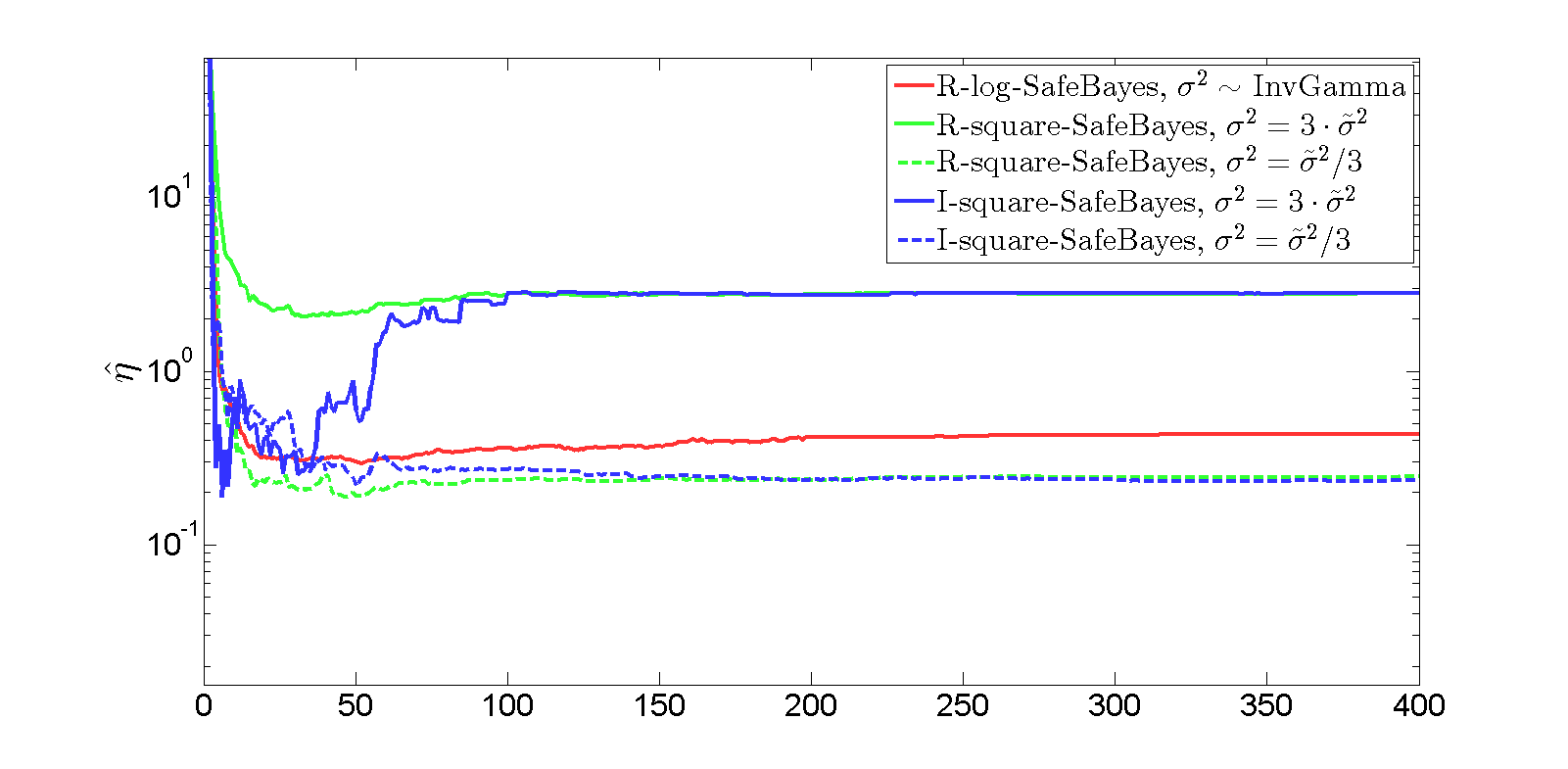}
\caption{\label{fig:mainfixedsigma} Bayesian Model Selection, fixed
  $\sigma^2$, for the model-wrong experiment of
  Figure~\ref{fig:mainexperimenta} with $\pmax = 50$. The second graph
  is a scaled version of the first. Since fixed $\sigma^2$ implies fixed
  overconfidence ratio, the overconfidence graph is not shown. 
For clarity in the $\eta$-graph we do not show standard deviations
  of the $\eta$'s. }}
\end{figure}

The model-correct experiment showed no surprises (all methods
performed well), so we only show results for the model-wrong
experiment, as described in Section~\ref{sec:preparing}, testing each
of Bayes, $R$-square- and $I$-square-SafeBayes twice: once based on a
model with variance $\sigma^2$ overly large (3 times
$\tilde\sigma^2$), and once with $\sigma^2$ overly small ($1/3$ times
$\tilde\sigma^2$) variance. To allow precise comparison with the
results in the main text, we also show behavior of $R$-log-SafeBayes
with varying variance (defined precisely as in
Figure~\ref{fig:mainexperimenta}) in Figure~\ref{fig:mainfixedsigma}.

\subsubsection{Ridge Regression Experiments, Fixed $\sigma^2$}
\label{sec:ridgefixedsigma}
Again we only show results for the model-wrong experiment.

Note that here standard Bayes --- as can be seen from
plugging $\eta = 1$ into (\ref{eq:olvg}) --- does not depend on $\sigma^2$
and thus coincides in terms of square-risk behavior with standard
Bayes in the variable $\sigma^2$ case as in
Figure~\ref{fig:ridgevarysigmaa}. Also (see below (\ref{eq:olvg}))
$I$-square-SafeBayes for fixed $\sigma^2$ does not itself depend on
$\sigma^2$ and simply minimizes the cumulative sum of squared errors.

Just as for ridge regression with variable $\sigma^2$, one may
equivalently interpret the $\eta$-generalized-posterior means
$\bar{\beta}_{i,\eta}$ as the standard, nongeneralized Bayesian
posterior means that one would get with a modified prior on $\beta$,
proportional to the original prior raised to the power $\eta^{-1}$
(see above (\ref{eq:invi}), Section~\ref{sec:ridge}). It may then once
again seem reasonable to learn $\eta$ itself in a Bayesian- or
likelihood-based way such as empirical Bayes.\footnote{In the present
  setting, learning $\eta$ by empirical Bayes has a second
  interpretation: if one fixes the variance $\sigma^2$ appearing in
  the prior on $\beta$, uses the linear model with a different
  variance $\sigma'^2$, and then learns $\sigma'^2$ by empirical
  Bayes, the result is identical to fixing $\sigma'^2 = \sigma^2$ and
  learning $\eta$ by empirical Bayes.}  Indeed, this was suggested
implicitly as early as 1999 by one of us \citep{Grunwald99a}. The
procedure described in Section 3.4.3 (`hierarchical loss') of
\cite{bissiri2016general} also arrives, via a different derivation, at
a similar prescription for finding $\eta$ (we immediately add that the
authors describe many ways for determining $\eta$, of which this is
just one). Unfortunately, just as for the empirical Bayes learning of
$\eta$ with varying $\sigma^2$, the figures below indicate that it
does not perform well at all.

\begin{figure}[htp]{\hspace*{0.15\textwidth}
\includegraphics[width=0.7\textwidth]{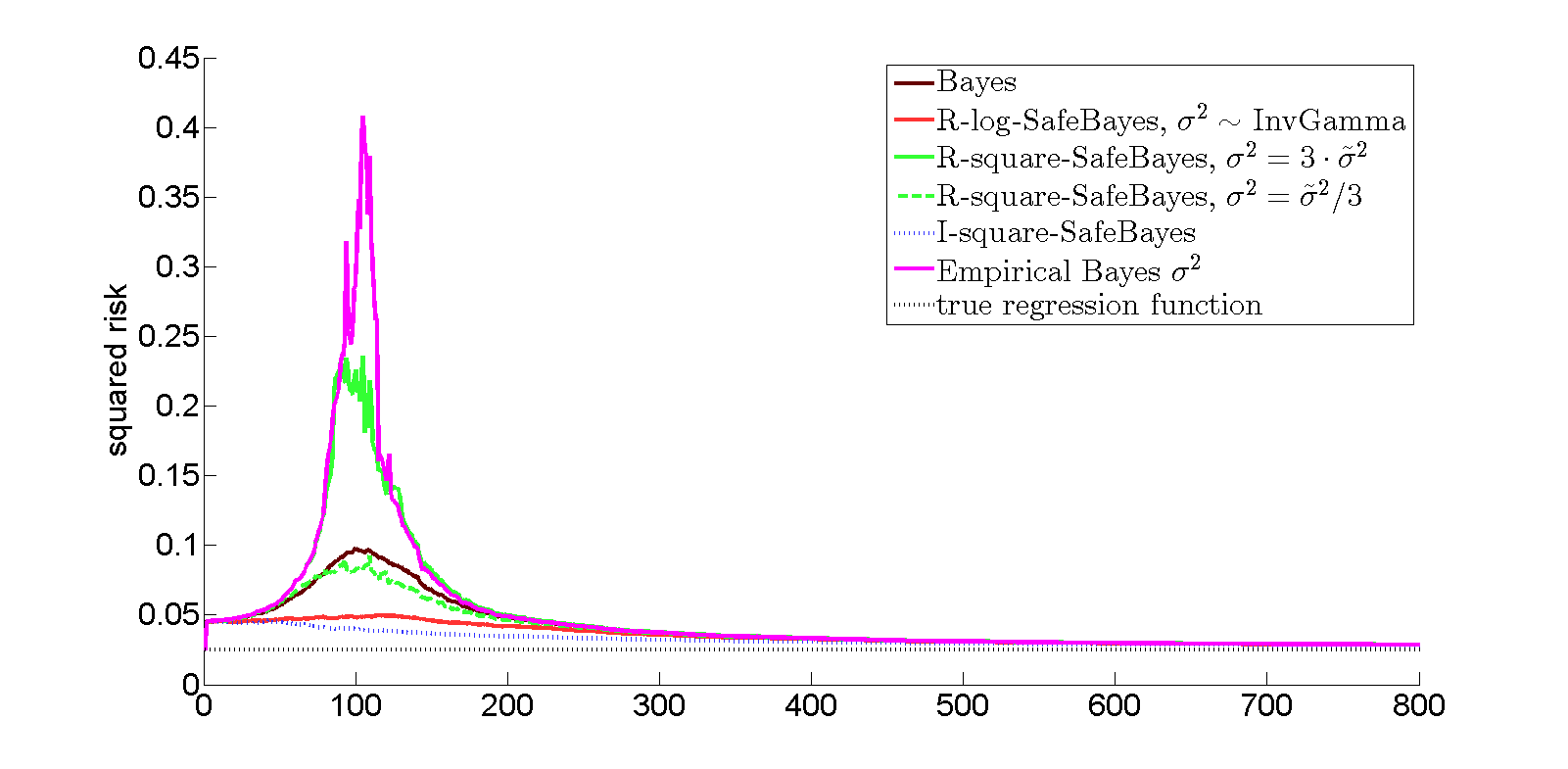} \\
\hspace*{0.15\textwidth}
\includegraphics[width=0.7\textwidth]{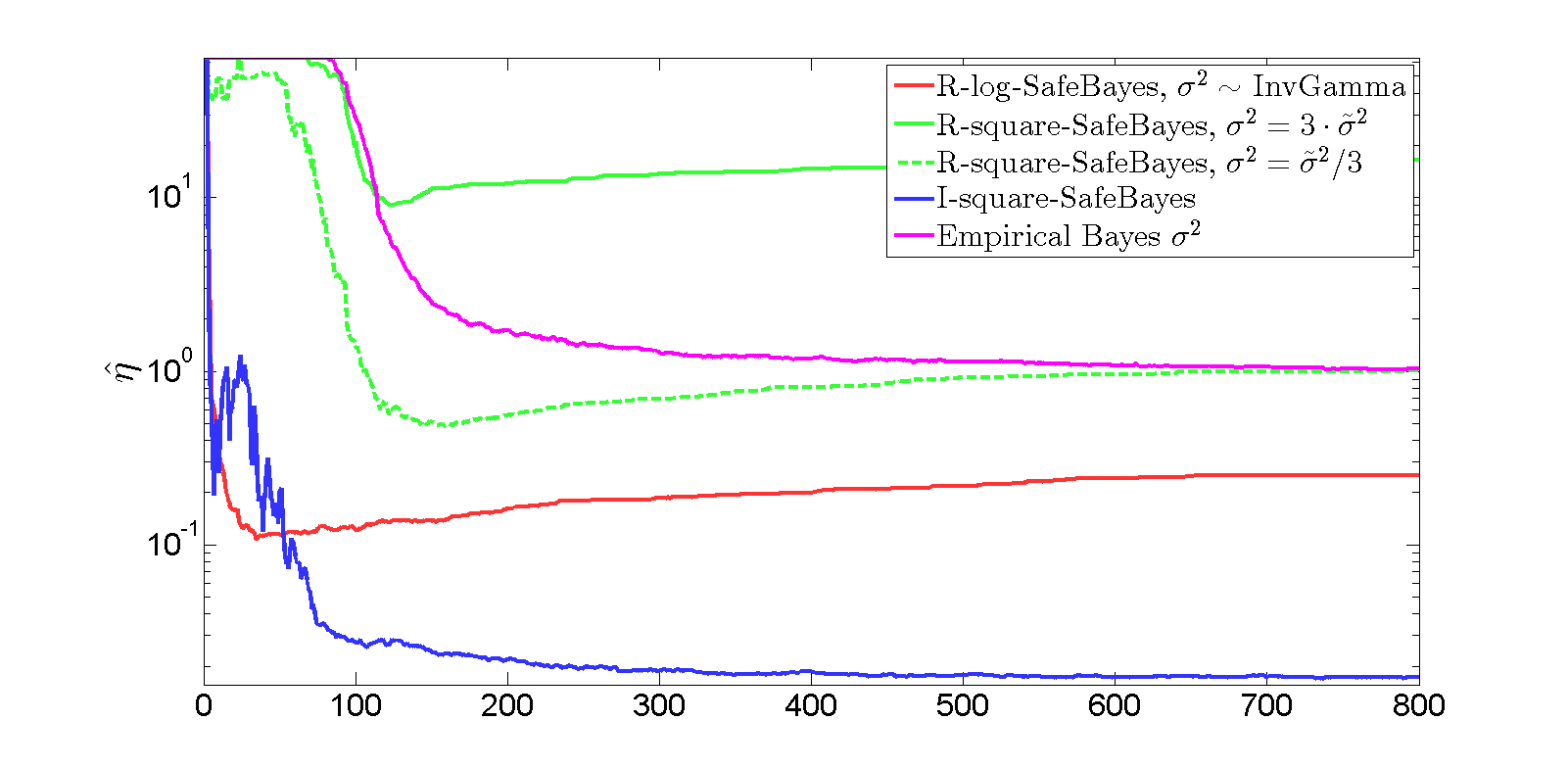}
\caption{\label{fig:ridgefixedsigma} Bayesian Ridge Regression: same
  graphs as in Figure~\ref{fig:ridgevarysigmaa},   for fixed $\sigma$ and the
  model-wrong experiment conditioned on $\nc := \pmax = 50$. Note the
  difference in scale for the risk in this figure and Figure~\ref{fig:mainfixedsigma}.
}}
\end{figure}

\paragraph{Conclusion}
Standard Bayes again performs comparably badly in both experiments
(note the difference in scale in the first graphs of
Figure~\ref{fig:mainfixedsigma} and~\ref{fig:ridgefixedsigma}).
$I$-square-SafeBayes behaves excellently in both experiments.  But now
in the ridge experiment $R$-square-SafeBayes becomes a highly
problematic method for small samples, worse even than standard Bayes.
The reason is its dependence on the specified $\sigma^2$ as can be
clearly seen from (\ref{eq:fronkie}). If $\sigma^2$ was set to be much
larger than the actual average prediction error on the sample, then
the third term in (\ref{eq:fronkie}) dominates. This term decreases
with $\eta$ and thus automatically pushes $\hat\eta$ `upward' by an
arbitrary amount.  The term also decreases with $n$, so that the
problem disappears at a large enough sample size. The problem did not
occur in the model averaging experiment; we suspect that this is
because in this experiment, there is substantial prior mass on a small
model $\nc = 4)$ containing the pseudo-truth, and for this submodel,
the final term in (\ref{eq:fronkie}) (which is approximately linear in
$p$) is much smaller than for $\nc = 50$ and does have not such a
strong influence. 

\subsection{Slightly Informative Prior}\label{sec:slightly}
Again we only consider model-wrong experiments. 
Within each model, we now use the following prior parameters: $\betao =
\mathbf{0}$ and $\Sigma_0 = 10^{3} \mathbf{I}$ for the multivariate
normal distribution on $\beta$; and $a_0 = 1$ and $b_0 = \sigma^{*2}
a_0$ (as before) for the inverse gamma distribution on $\sigma^2$
(where $\sigma^{*2}$ is the true variance of noise in our data, as
defined in Section~\ref{sec:truth}). We repeated the model-wrong
experiment of Section~\ref{sec:modsel} with $\pmax=50$ with this
slightly informative prior and obtained similar results to those
obtained using our original informative prior with $\Sigma_0 =
\mathbf{I}$: Bayes performs badly roughly between samples 90 and 130
and has some risk spikes before that so that its overall performance
is comparable to before, while $R$-SafeBayes and $I$-SafeBayes both obtain
good risks.

We also repeated the model-wrong experiment for ridge regression
(Section~\ref{sec:ridge}). Here the effect of the new prior on Bayes'
performance is similar: the square-risk peaks at a larger value, but
in a smaller range of sample sizes. However, the effect of changing
the learning rate is different in this experiment than what we have
seen before: here one can take $\eta$ {\em very\/} small and still get
good results. So in a sense, the problematic behavior of Bayes has a
trivial solution here: just pick a very small but fixed $\eta$.
$R$-log-SafeBayes was too conservative in this, $I$-log-SafeBayes did
fine. $R$-log-SafeBayes became competitive again however, if we used the
discounting version described in Section~\ref{sec:etalater} below. 

We omit the pictures corresponding to model selection/averaging
(Section~\ref{sec:modsel}) as they show no surprises; but in
Figure~\ref{fig:slightly} we do repeat
the pictures for ridge regression (Section~\ref{sec:ridge}), because
they do give additional insight:
\begin{figure}[htp]{\hspace*{0.15\textwidth}
\includegraphics[width=0.7\textwidth]{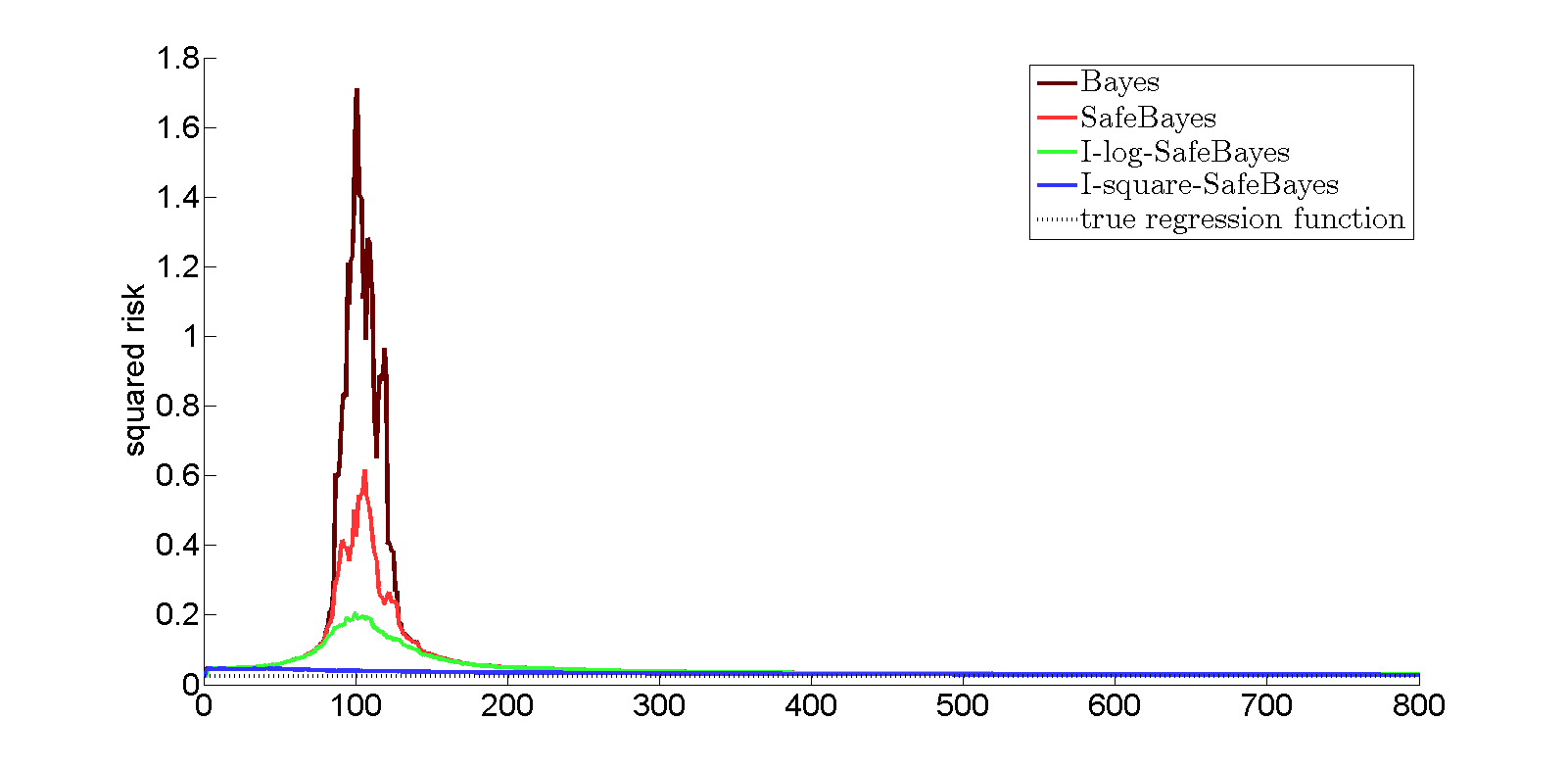} \\
\hspace*{0.15\textwidth}
\includegraphics[width=0.7\textwidth]{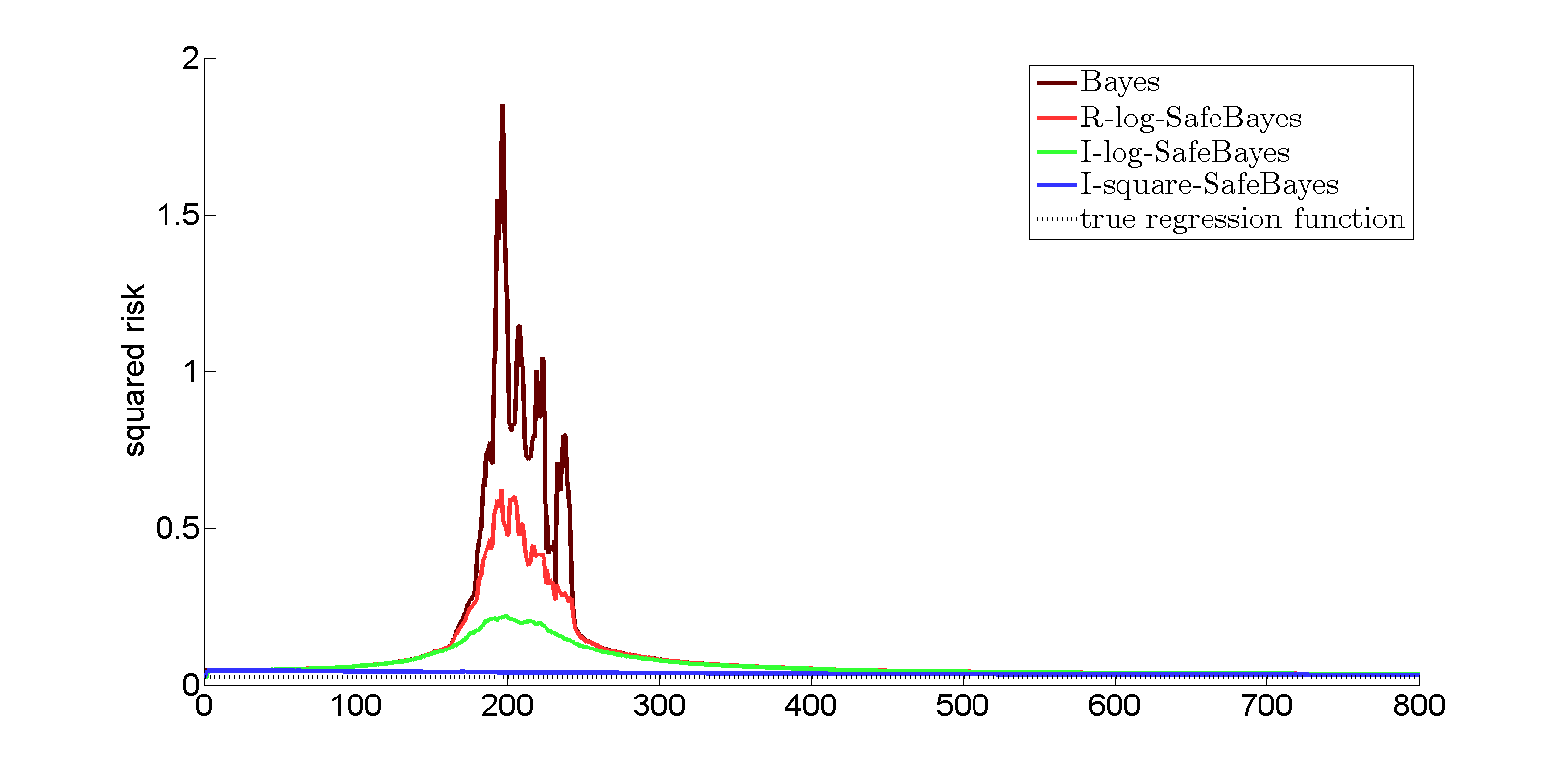} \\
\hspace*{0.15\textwidth}
\includegraphics[width=0.7\textwidth]{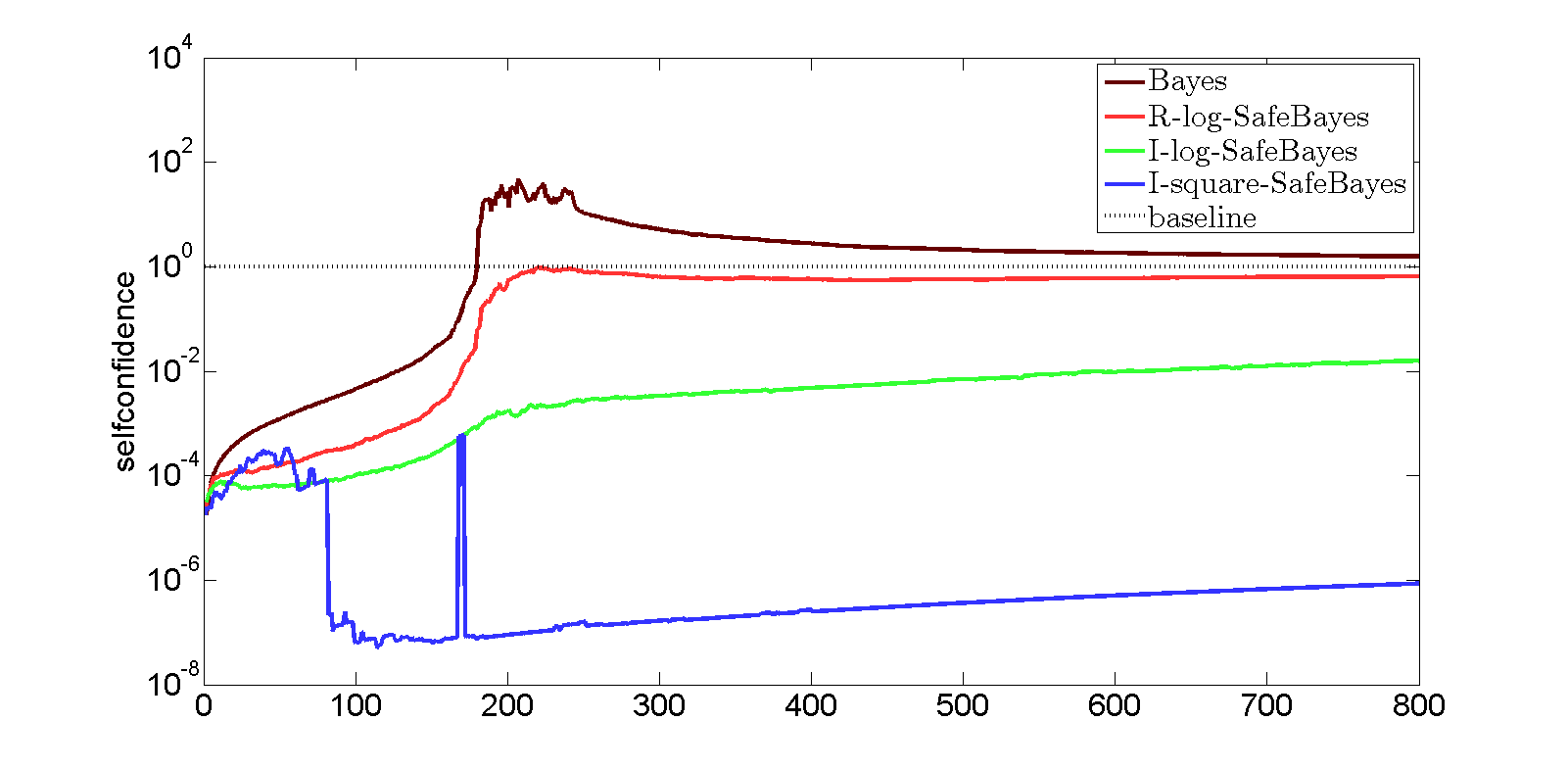}
\caption{\label{fig:slightly} Top two graphs: square-risk for two
  different ridge experiments. In both experiments the slightly
  informative prior of Section~\ref{sec:slightly} is used.  In the first experiment $\nc=50$; in the
  second $\nc = 100$; otherwise the experiments are just as the `wrong
  model experiment' of Section~\ref{sec:ridge},
  Figure~\ref{fig:ridgevarysigmaa}, but we also included performance
  of $I$-square-SafeBayes. Final graph shows
  self-confidence for the $\nc=100$ case for Bayes and SafeBayes, on a
  logarithmic scale because of the range of values involved. }}
\end{figure}
Note that the phenomenon is now much more `temporary'. In the
beginning, it seems that there is a sort of cancellation between the
influence of the irrelevant variables and standard Bayes behaves fine.
However, if we increase the number of irrelevant variables, the
problem (while starting at a later sample) takes longer to recover
from.
\subsection{Prior as advised by Raftery et al.}
\label{sec:raftery}
In \cite{raftery1997bayesian}, some guidelines for choosing priors in
regression models are given. Letting $\betao$ denote the prior mean,
one of their recommendations is that the prior densities for $\beta =
\betao$ and $\beta = \betao + \mathbf{1}$ should differ by a factor of
at most $\sqrt{10}$.
The prior density on $\beta$
marginalized over $\sigma^2$ follows a
multivariate t-distribution, and the factor in question varies with
the dimensionality of $\beta$, so that models of larger order are given less
informative priors. In our case, we find that the resulting prior is
always less informative than our original prior, and for model $\cM_{10}$
and above (i.e.~$\beta$ of dimension 11 or larger),  it becomes even less
informative than the prior introduced in the previous section.

For the prior on $\sigma^2$, \citeauthor{raftery1997bayesian} advise
that the density should vary by no more than a factor 10 in a region
of $\sigma^2$ from some small value to the sample variance of $y$. For
our choice of hyperparameters $a_0 = 1$, $b_0 = 1/40$, the mode of
$\pi(\sigma^2)$ is at $b_0 / (a_0 + 1) = 1/80$, and the density is
within a factor 10 of this maximum in the approximate region $(0.0037,
0.0941)$. For the correct model experiments, the actual variance of
$Y$ is $0.065$; for the wrong model experiments, it is $0.045$ (with a
larger variance for `good' points and zero variance for `easy'
points). For both experiments, the factor-10 condition holds between
$\textrm{Var}(Y) / 12$ and $\textrm{Var}(Y)$. We conclude that this
prior satisfies the guidelines in \citeauthor{raftery1997bayesian}
quite well.

We will refer to the prior described above as Raftery's prior (even
though it is really a different prior for each model order). Using
this prior, we found the following experimental results.

In the model-wrong experiment with model selection/averaging
(Section~\ref{sec:modsel}) with our original prior replaced by
Raftery's prior, Bayes performs somewhat {\em better\/} than
$R$-log-SafeBayes (except on very small sample sizes). However,
$I$-log-SafeBayes performs as well as Bayes, and so does the
$R$-log-SafeBayes variant that discounts half of the initial sample when
choosing the learning rate (see Section~\ref{sec:etalater}).

This might suggest that Raftery's prior could be used to accomplish
the same kind of safety against wrong models as SafeBayes provides, at
least in a model selection context. To test this, another experiment
was performed where the fraction of `easy' points was increased to
75\%. In this experiment, the misbehavior of Bayes seen in
Section~\ref{sec:modsel} returned worse than before, with risks a
factor 20 larger than before, whereas the SafeBayes methods continued
to work fine. This suggests that Raftery's prior can not be relied on
if the severeness of misspecification is unknown.

If Raftery's prior is used for model selection with a correct model,
Bayes and the SafeBayes variants perform well, and very similarly to
each other.

For ridge regression, the results with Raftery's prior for both the
correct and the incorrect model experiment are very similar to those
with the slightly informative prior, except that the peak in the risks is higher for
all methods.

\subsection{The $g$-prior}
\begin{figure}[htp]{\hspace*{0.15\textwidth}
\includegraphics[width=0.7\textwidth]{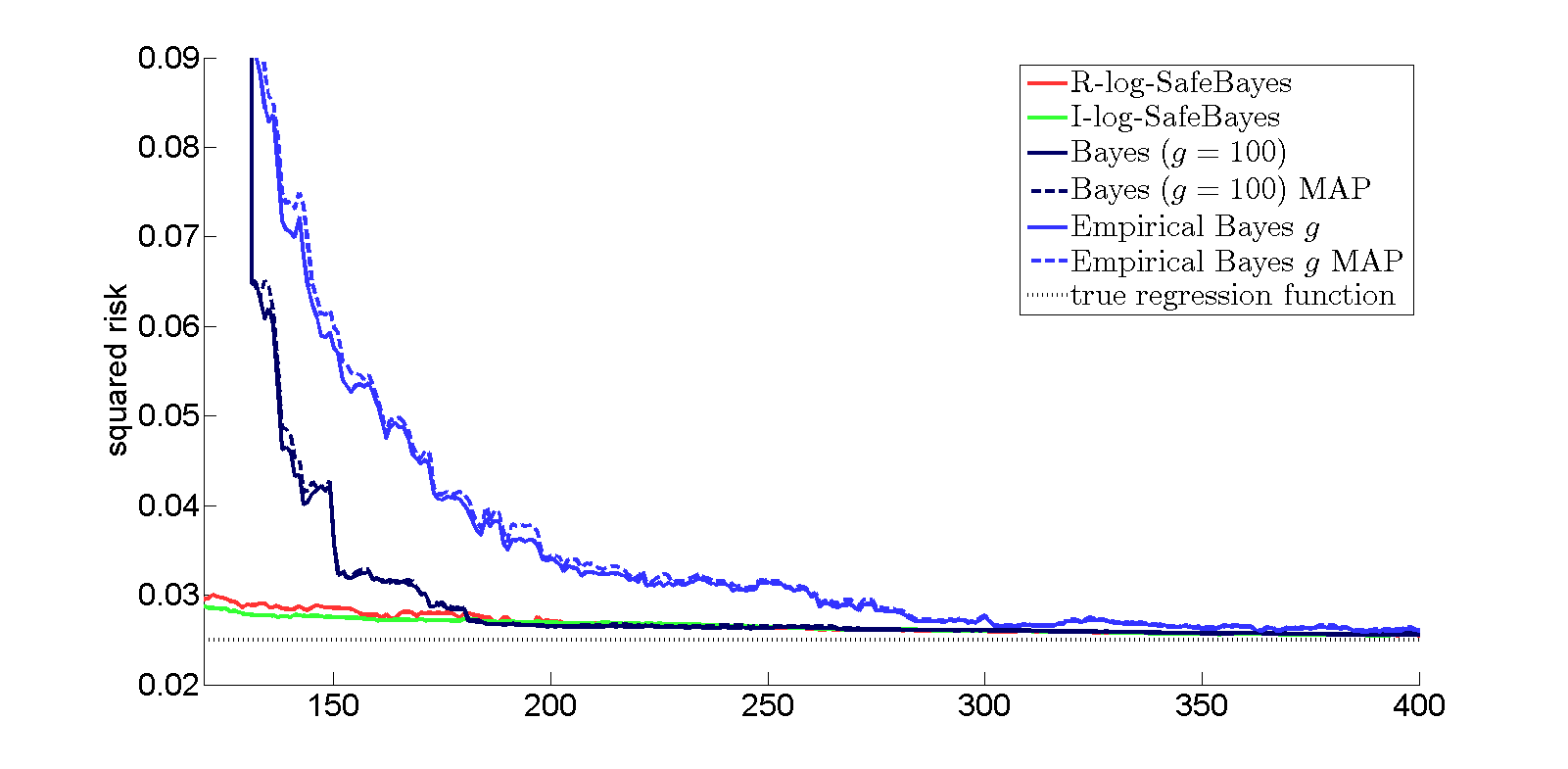} 
\caption{\label{fig:gprior}Risk as a function of sample size (starting
  at the
  first sample size at which the $g$-prior is defined) for model
  averaging and selection based on the $g$-prior in the model-wrong
  experiment of Figure~\ref{fig:mainexperimenta} both with $g = 100$
  and with $g$ chosen by Empirical Bayes at each sample size}}
\end{figure}

Another prior we experimented with was the $g$-prior, a popular choice
in model selection contexts \citep{Zellner86,liang2008mixtures}. For
all definitions we refer to the latter paper. 
In contrast to all other priors we considered, the $g$-prior depends
on the design matrix ${\bf X}^T_n {\bf X}_n$, and hence can only be
used in settings where this matrix, and hence the eventual sample size
of interest $n$, is given once and for all.  For this reason, we
decided to depict in Figure~\ref{fig:gprior}, for each value of $n$,
the risk obtained when predicting the $n$-th data point with the
posterior calculated from  the $g$-prior corresponding to the
first $n$ covariates $(x_1,\ldots, x_n)$ and observed data $y^{n-1}$. This is
subtly different from our previous graphs (e.g.\ 
Figure~\ref{fig:mainexperimenta}--\ref{fig:mainexperimentd}) that show how
the risk evolves as $n$ increases in a {\em single\/} run of the experiment,
averaged over 30 runs.

The graph is not shown starting at $n=0$, because of another
difference between the $g$-prior and the priors we used in other
experiments: 

Because of the same design dependence, with the $g$-prior, the
posterior on $\beta$ remains a degenerate distribution on an initial
segment of outcomes. For example, with $\cM_p$ for $p = 50$, the
matrix ${\bf X}_n^T {\bf X}_n$ is singular until at least 50
\emph{different} design vectors have been observed. For our
model-wrong experiment, this means that on average, about a 100
observations are required before the posterior becomes nondegenerate;
this explains why Figure~\ref{fig:gprior} starts at $n$ a little over
a 100.

The experimental results clearly indicate that the $g$-prior does not
deal with our data in a satisfactory way, regardless of the value of
$g$. Of the values of $g$ we tried (up to $10^4$), $g \approx 100$
(shown in the graph) yielded the smallest squared risk around sample
size $n=200$; for larger sample sizes, larger values of $g$ were
better, but only slightly. Furthermore, (as in fact we expected by
analogy to learning $\eta$ with Empirical Bayes), the value of $g$
found by Empirical Bayes is not optimal for dealing with our data and
only makes things worse: larger values of $g$ (which put more weight
on the data) would yield smaller risks.
\section{ Experiments on Variations on the Method}
Below we look at a number of other
more or less promising alternative approaches to modifying standard
Bayes.

\subsection{An Idea to be Explored Further: Discounting Initial Observations}
\label{sec:etalater}
Just like standard Bayes, all our SafeBayesian methods are, at heart,
{\em prequential\/} \citep{Dawid84}. All prequential methods suffer to
a greater or lesser extent from the {\em start-up problem\/}
\citep{ErvenGR07,wong2004improvement}: sequential predictions based on a model
$\cM_p$ may perform very badly for the first few samples. While they
quickly recover when the sample size gets large, the behavior on the
first few samples may dominate their cumulative prediction error for a
while, leading to suboptimal choices for moderate $n$. We can address
this issue in several ways. A very simple method to `discount' initial
observations, apparently first used (implicitly) to modify standard
Bayes factors by \cite[Chapter 6]{Lempers71}, is to only look at the
cumulative sequential prediction error on the second half of the
sample, so that the first half of the sample merely functions as a
`warming-up' sample \citep{catoni2012catching}. Without claiming that
this is the `right' method to discount initial observations, we
experimented with it to see whether it can further improve the
performance of SafeBayes; for simplicity, we concentrated on
$R$-log-SafeBayes. 

We found that in most experiments, this new method for determining
$\eta$ performed very similarly to the standard method based on the
whole sample, sometimes slightly better and sometimes slightly worse,
making it hard to say whether the new method is an improvement or not.
Still, there are two experiments in which the new method performed
substantially better, namely the experiments with less informative
priors of Section~\ref{sec:slightly} and~\ref{sec:raftery}. 
Thus we cannot just dismiss
the idea of fitting $\eta$ based on only part of the data or more
generally, discounting initial observations, and it would be
interesting to explore this further in future work: of course taking
half of the data is rather arbitrary, and better choices may be
possible.  In particular, we may try a variation of {\em switching\/}
between $\eta$'s analogously to the switch distribution
\citep{ErvenGR07} to counter the startup problem.

\subsection{Other Methods for Model Selection: AIC, BIC, (generalized) Cross-Validation}

We tested the performance of several classic model selection methods
on the same data and models as in our main model selection/averaging
experiment, Section~\ref{sec:modsel}. We associated with each model
$\cM_{\nc}$ its standard (i.e.\ $\eta = 1$) Bayes predictive
distribution under the prior described in Section~\ref{sec:preparing}
(these generally perform better than the maximum likelihood
distributions based on $\cM_{\nc}$ whose use is more standard here).
We then ran leave-one-out cross-validation, 10-fold cross-validation
and GCV based on the predictions (posterior means/MAPs
$\bar{\beta}_{i,\eta}$) made by these predictive distributions.
\begin{figure}[htp]{\hspace*{0.2\textwidth}
\includegraphics[width=0.6\textwidth]{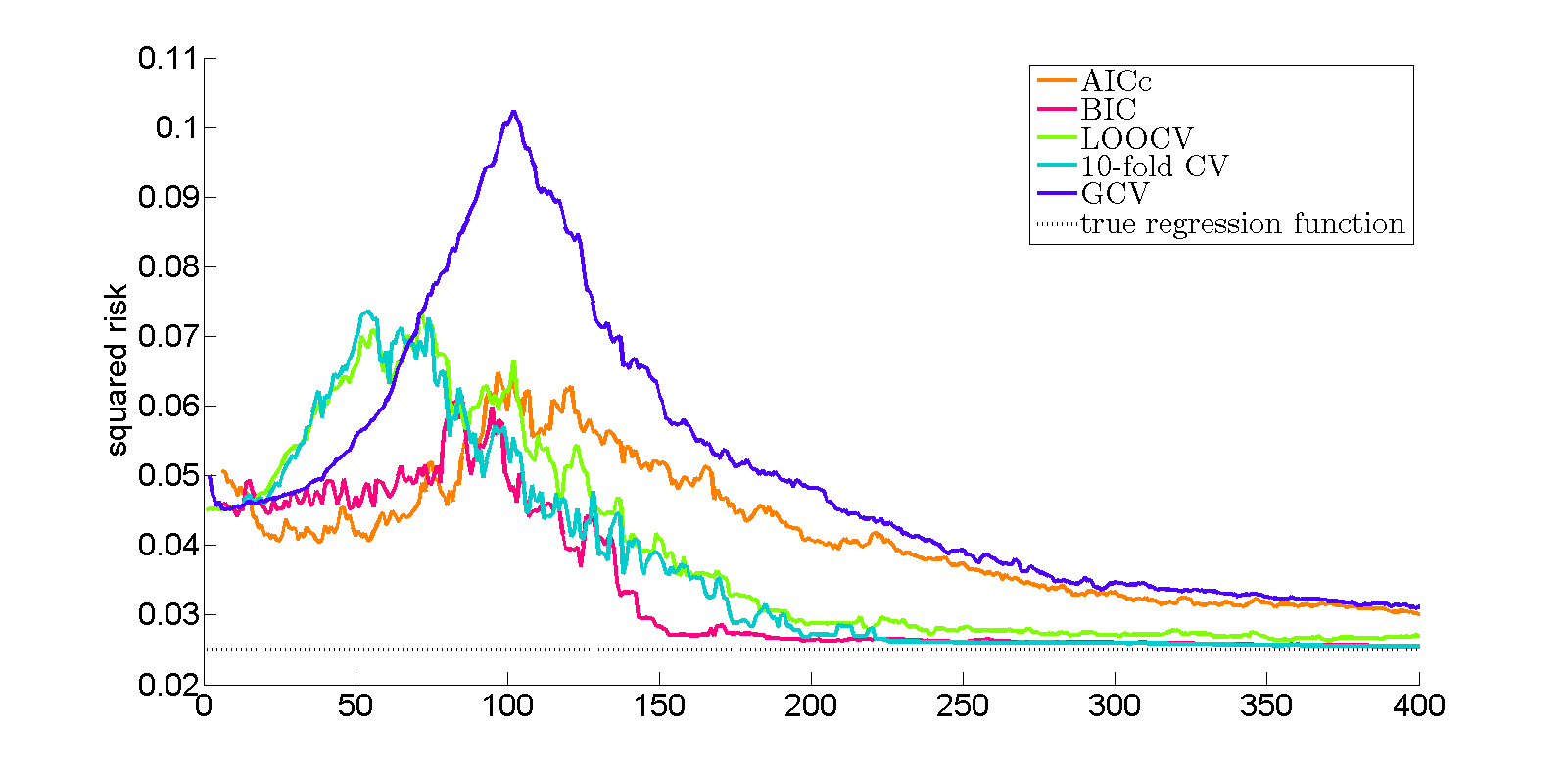} \\
\hspace*{0.2\textwidth}
\includegraphics[width=0.6\textwidth]{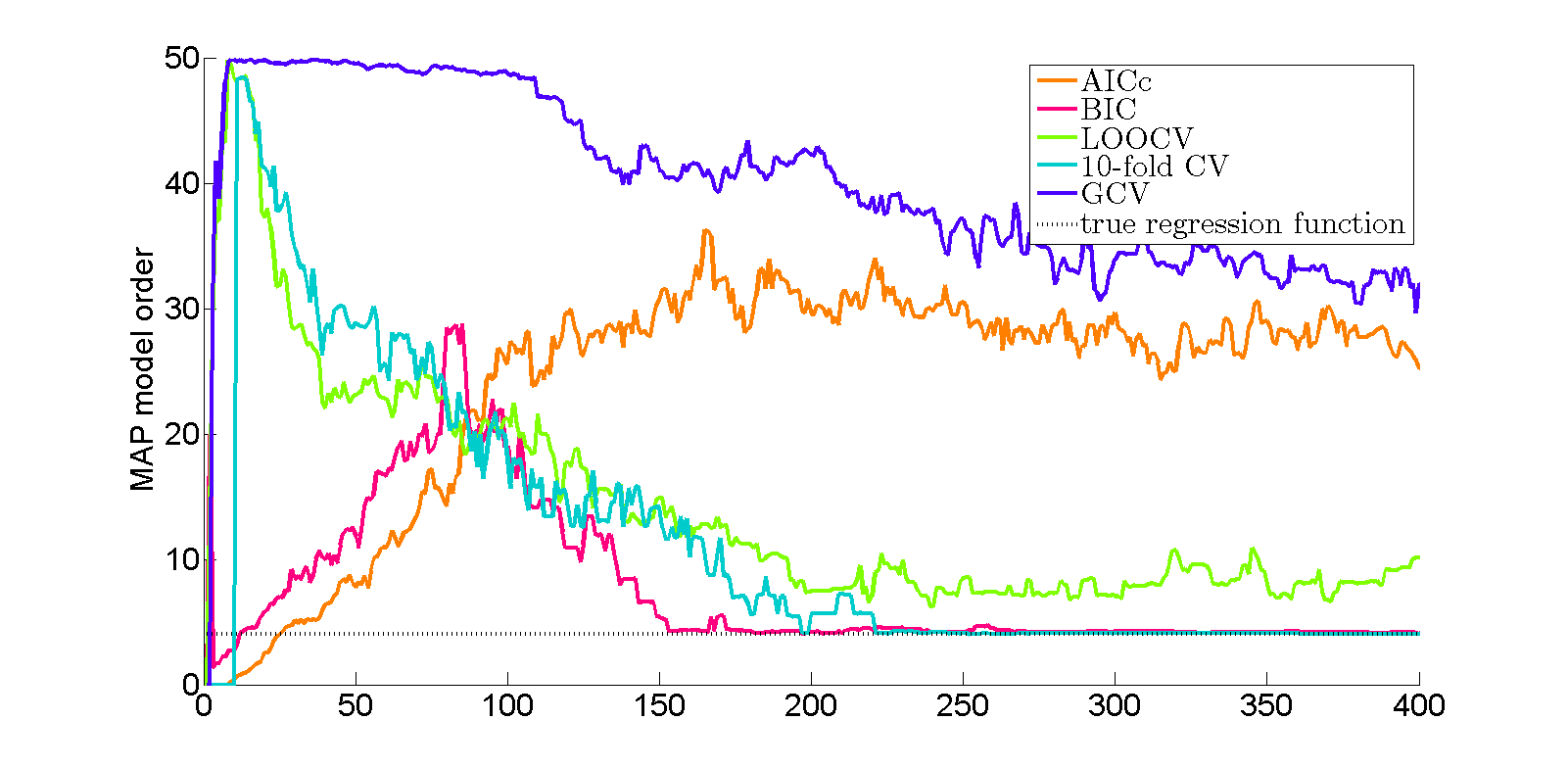}
\caption{
Squared risk and selected model order for five different model
selection methods. The risks in this graph are risks of single models
selected by each method (similar to the MAP risks shown for Bayes and
SafeBayes).
}}
\end{figure}
We also compared the models via AIC and BIC, where for AIC we used the
small-sample correction of \cite{HurvichT89}.

We see that AIC and generalized cross-validation have risks and
selected model orders similar to those of standard Bayes, though they
do not recover as well as Bayes when the sample size increases. Of the
other three methods, BIC and 10-fold cross-validation find the optimal
model and have smaller risks towards the end than leave-one-out
cross-validation, which continues to select larger-than-optimal models
with substantial probability. Note that none of the methods can
compete with SafeBayes on sample sizes below 150: SafeBayes's risk
goes down immediately after the start of the experiment while for all the
other methods it goes up first.  Also, SafeBayes finds the optimal
model quickly without first trying much larger models.

\subsection{Other Methods for Learning $\eta$: Cross-Validation on
  Log-Loss and on Squared Loss}
\begin{figure}[htp]{\hspace*{0.2\textwidth}
\includegraphics[width=0.6\textwidth]{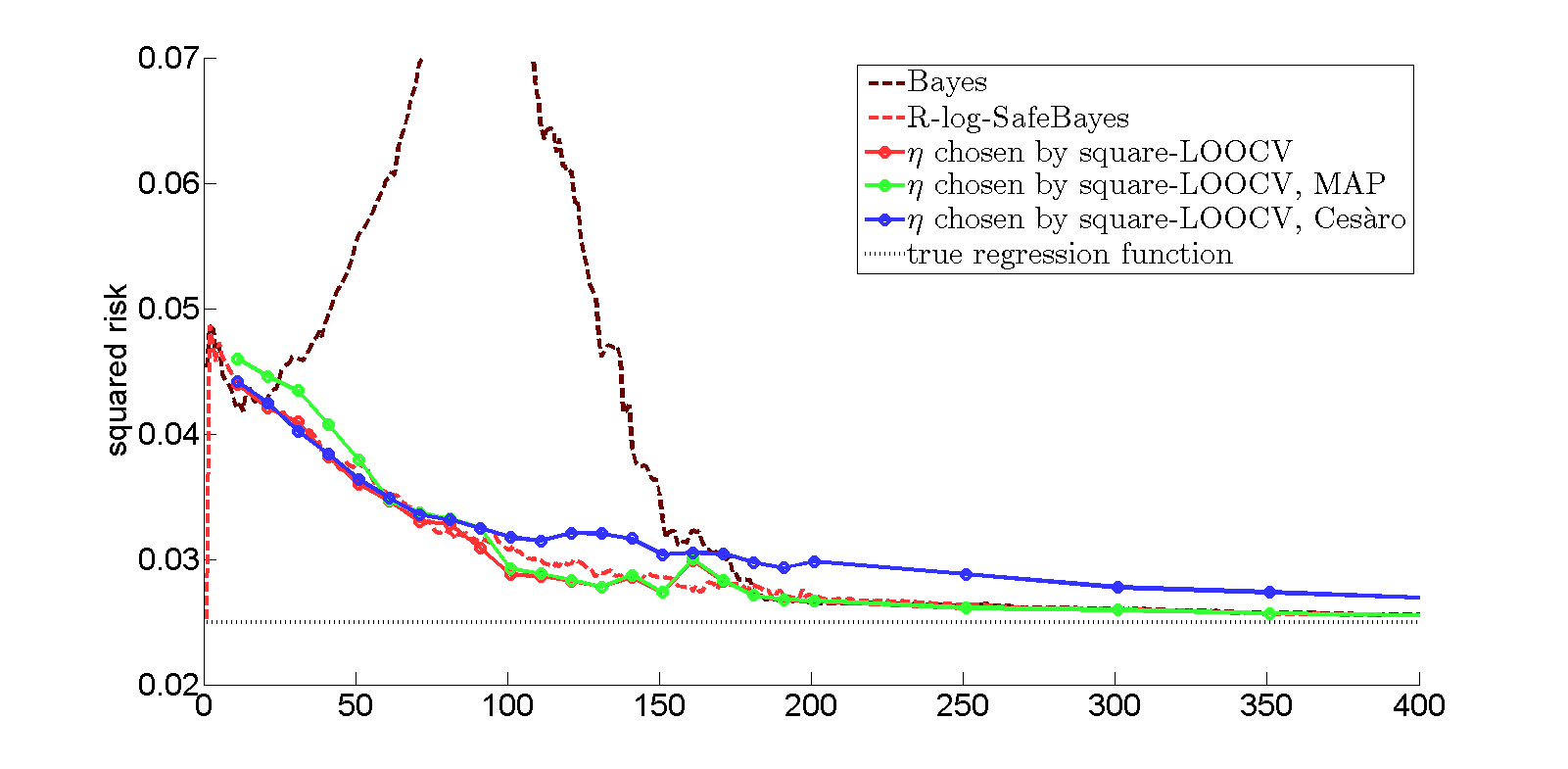} \\
\hspace*{0.2\textwidth}
\includegraphics[width=0.6\textwidth]{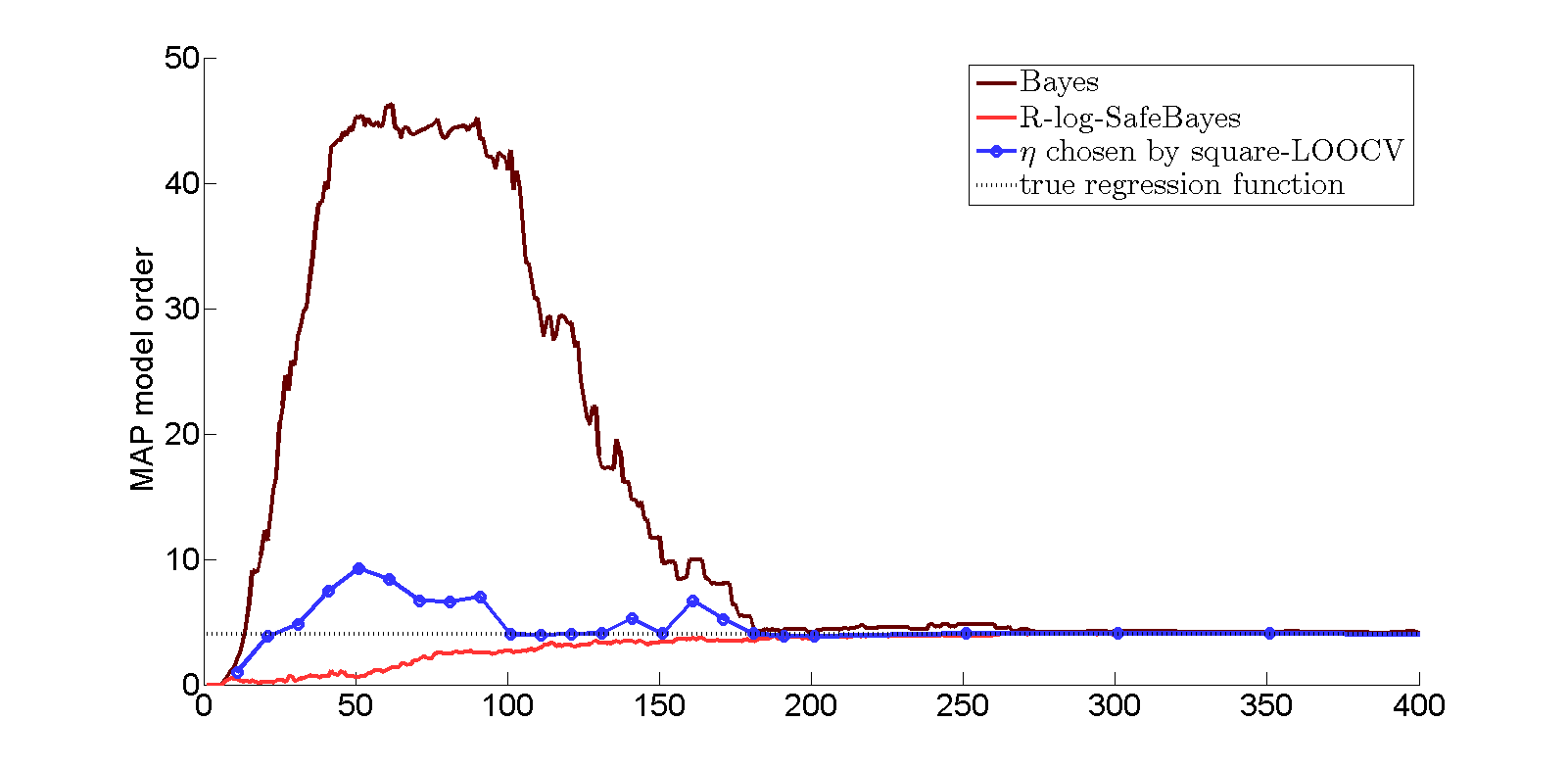} \\
\hspace*{0.2\textwidth}
\includegraphics[width=0.6\textwidth]{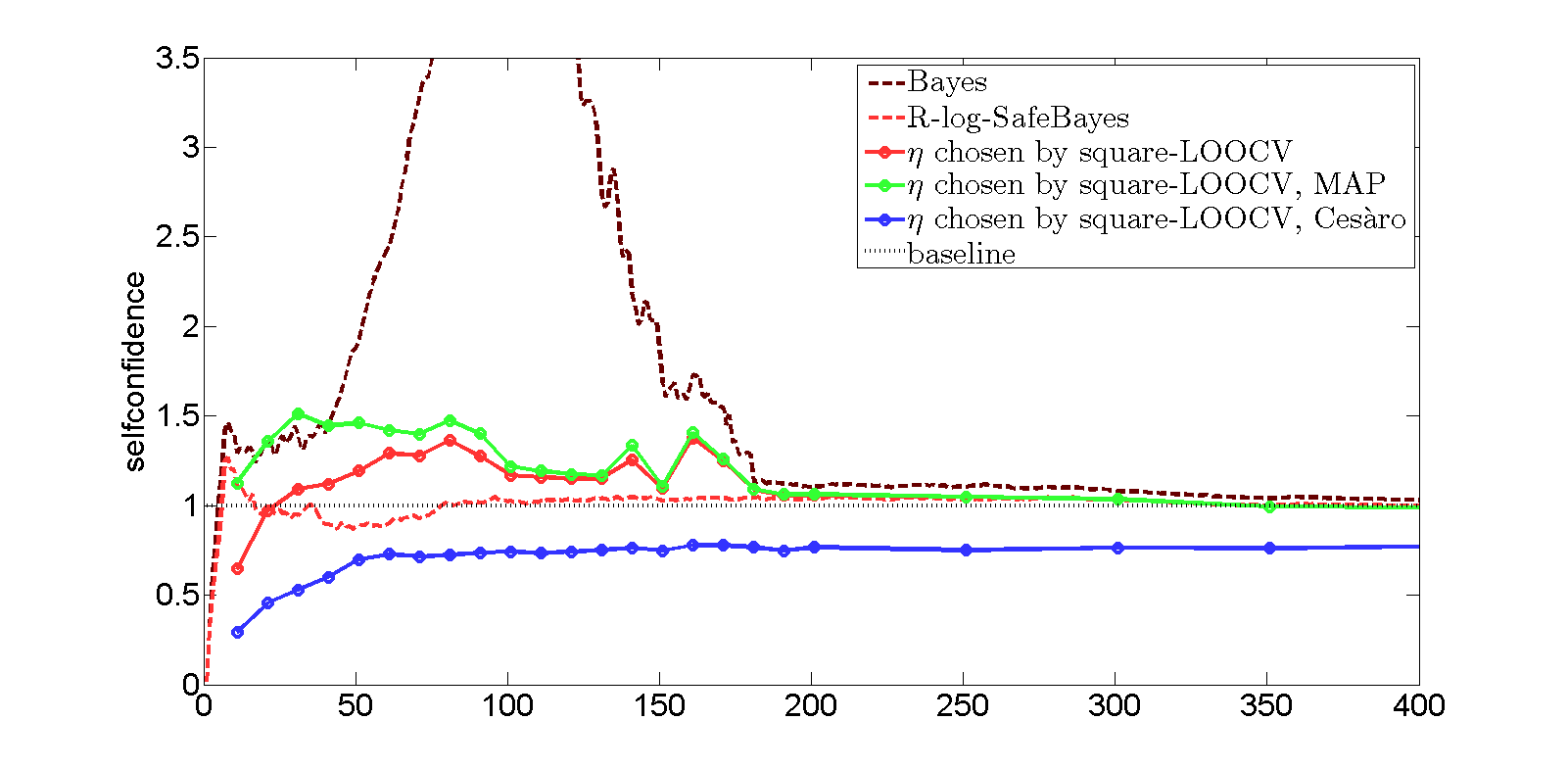} \\
\hspace*{0.2\textwidth}
\includegraphics[width=0.6\textwidth]{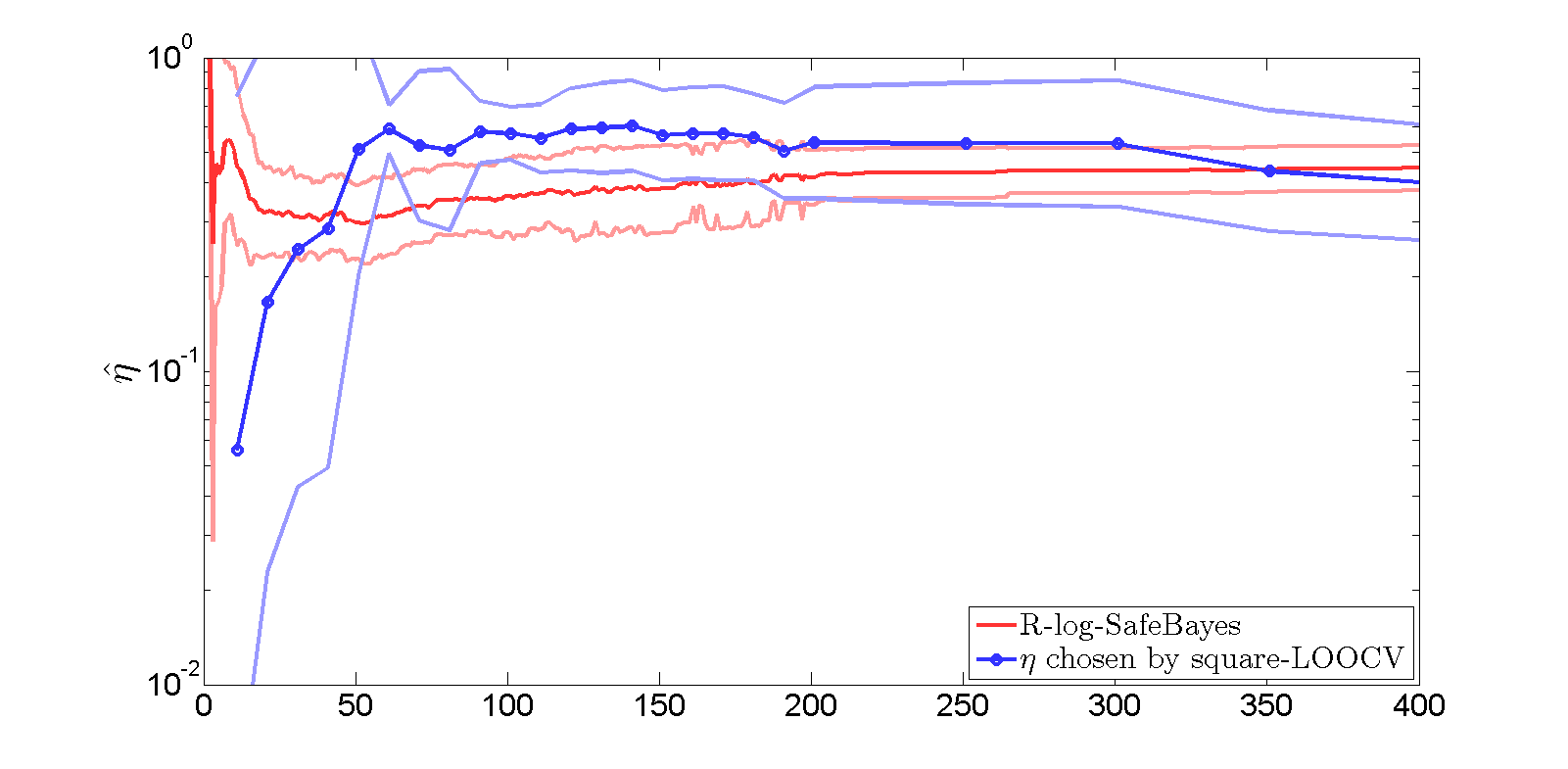}
\caption{\label{fig:canada} Analogue of Figure~\ref{fig:mainexperimenta} for determining
  $\eta$ by leave-one-out cross-validation with square-loss with the
  wrong-model experiment, $\pmax = 50$.}}
\end{figure}

As indicated in the introduction and Section~\ref{sec:instantiationb},
finding $\hat\eta$ by $I$-square-SafeBayes is somewhat similar to
finding $\hat\eta$ by leave-one-out cross-validation with the
squared-error loss, the difference being that $I$-square-SafeBayes finds
the optimal $\eta$ for predicting each point based on past data data
points rather than the optimal $\eta$ for predicting each point based
on all other data points. Since the leave-one-out method is often
employed in ridge regression, it seemed of interest to try out here as
well.  Figure~\ref{fig:canada} shows that LOO-cross validation indeed performs very
similarly to $R$-log and $I$-square SafeBayes in terms of
square-risk, but is consistently a bit worse in terms of
self-confidence; we do not have a clear explanation for this
phenomenon.

\begin{figure}[htp]{\hspace*{0.2\textwidth}
\includegraphics[width=0.6\textwidth]{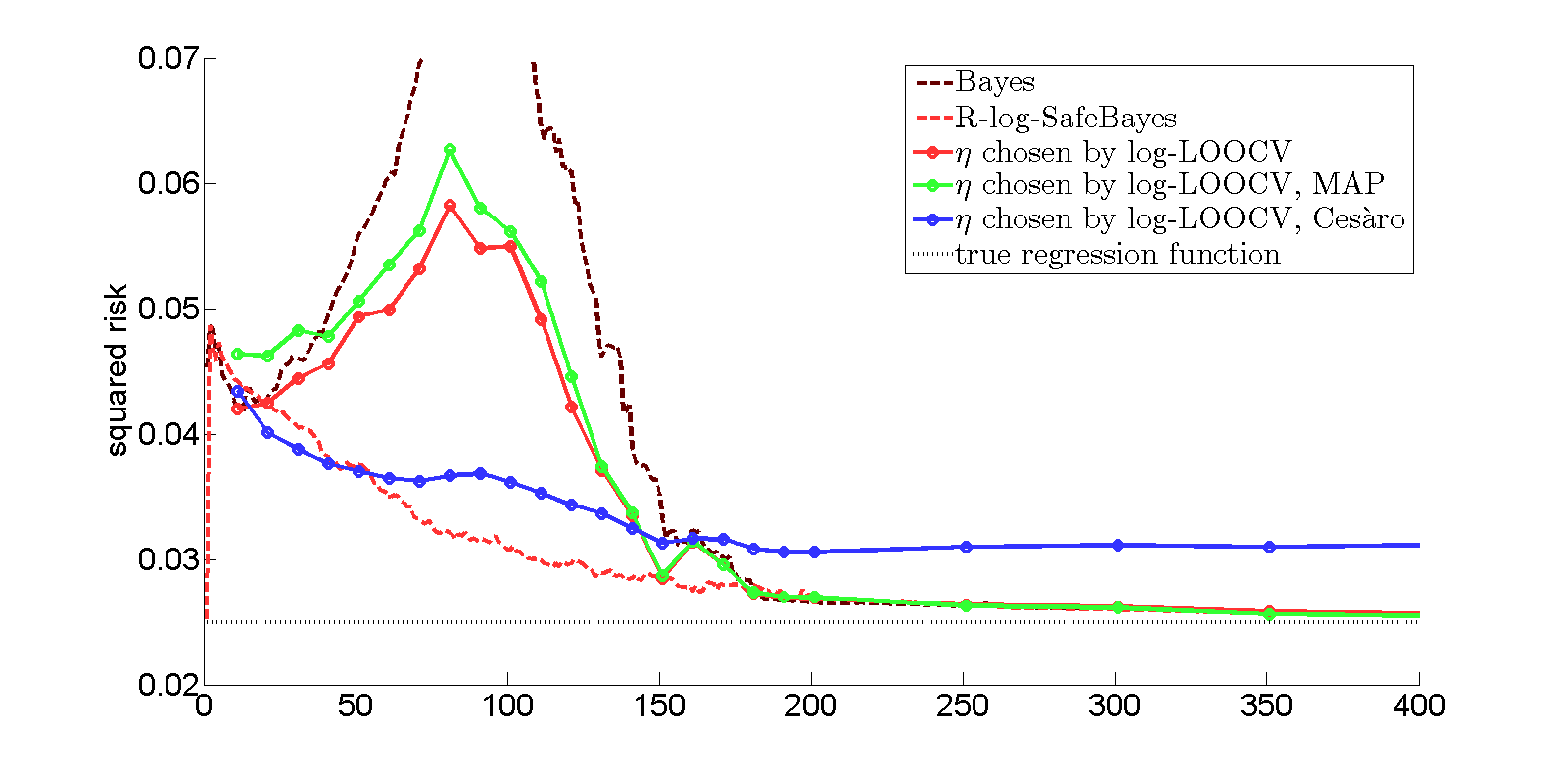} \\
\hspace*{0.2\textwidth}
\includegraphics[width=0.6\textwidth]{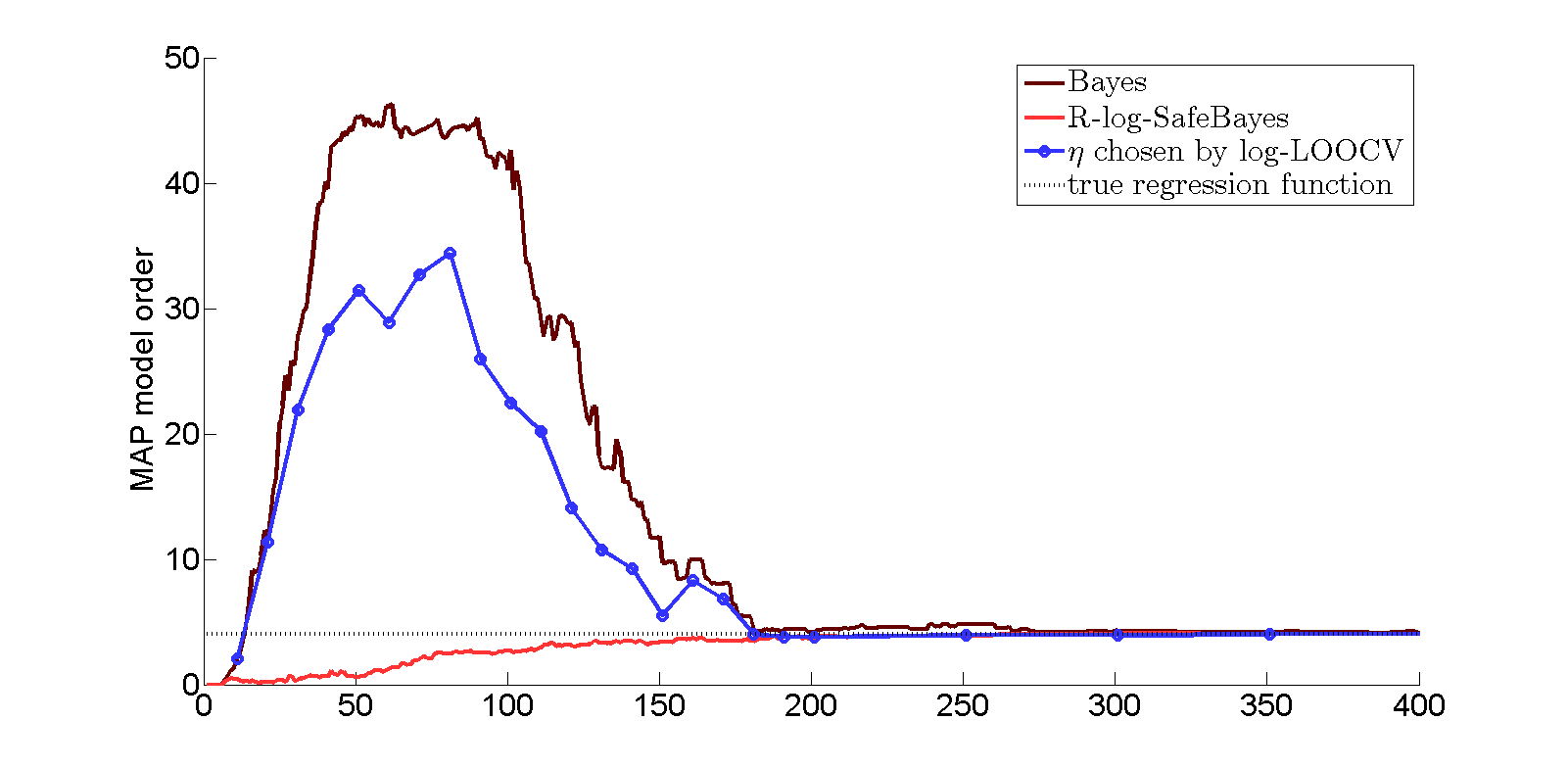} \\
\hspace*{0.2\textwidth}
\includegraphics[width=0.6\textwidth]{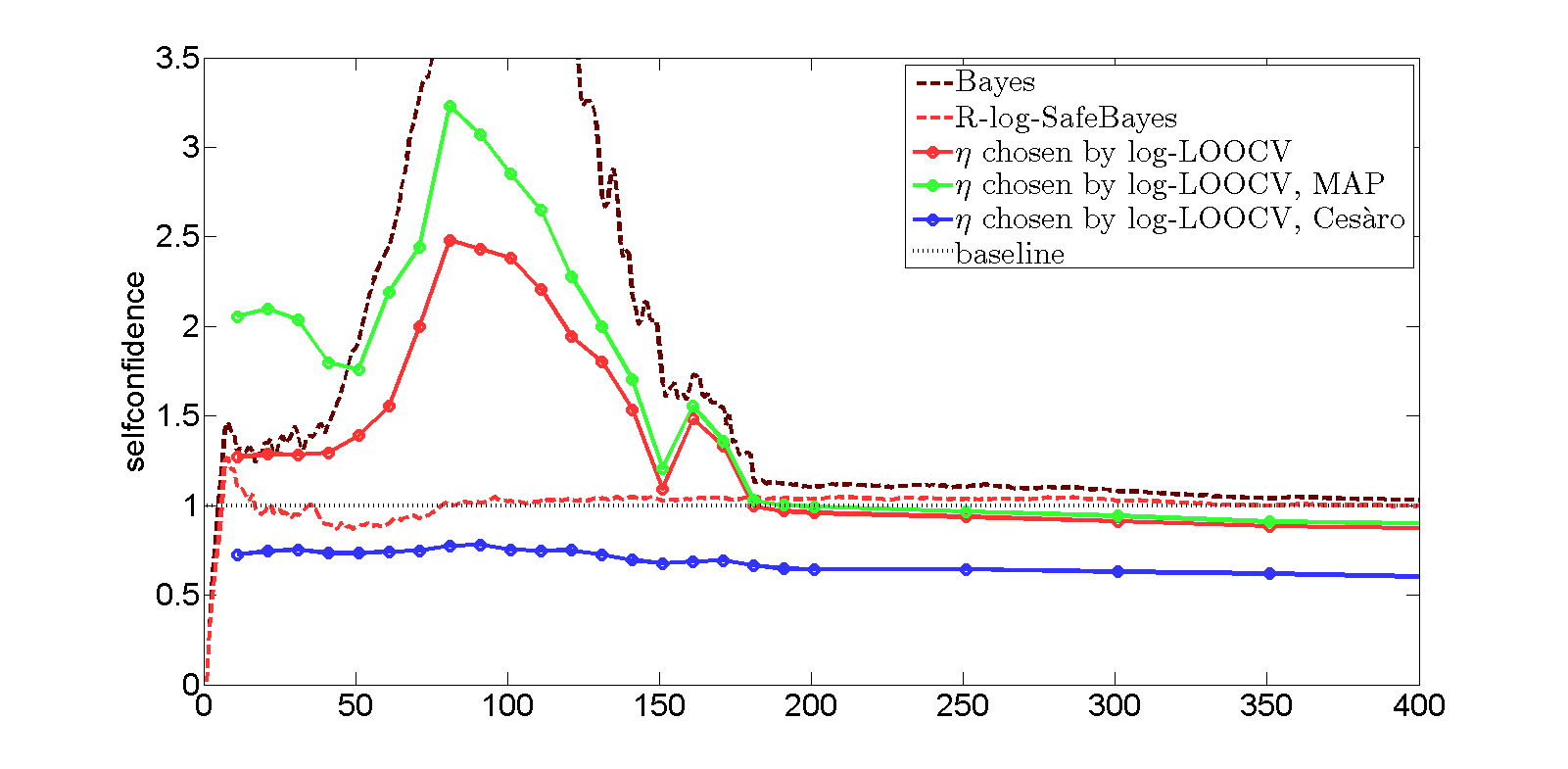}
\\
\hspace*{0.2\textwidth}
\includegraphics[width=0.6\textwidth]{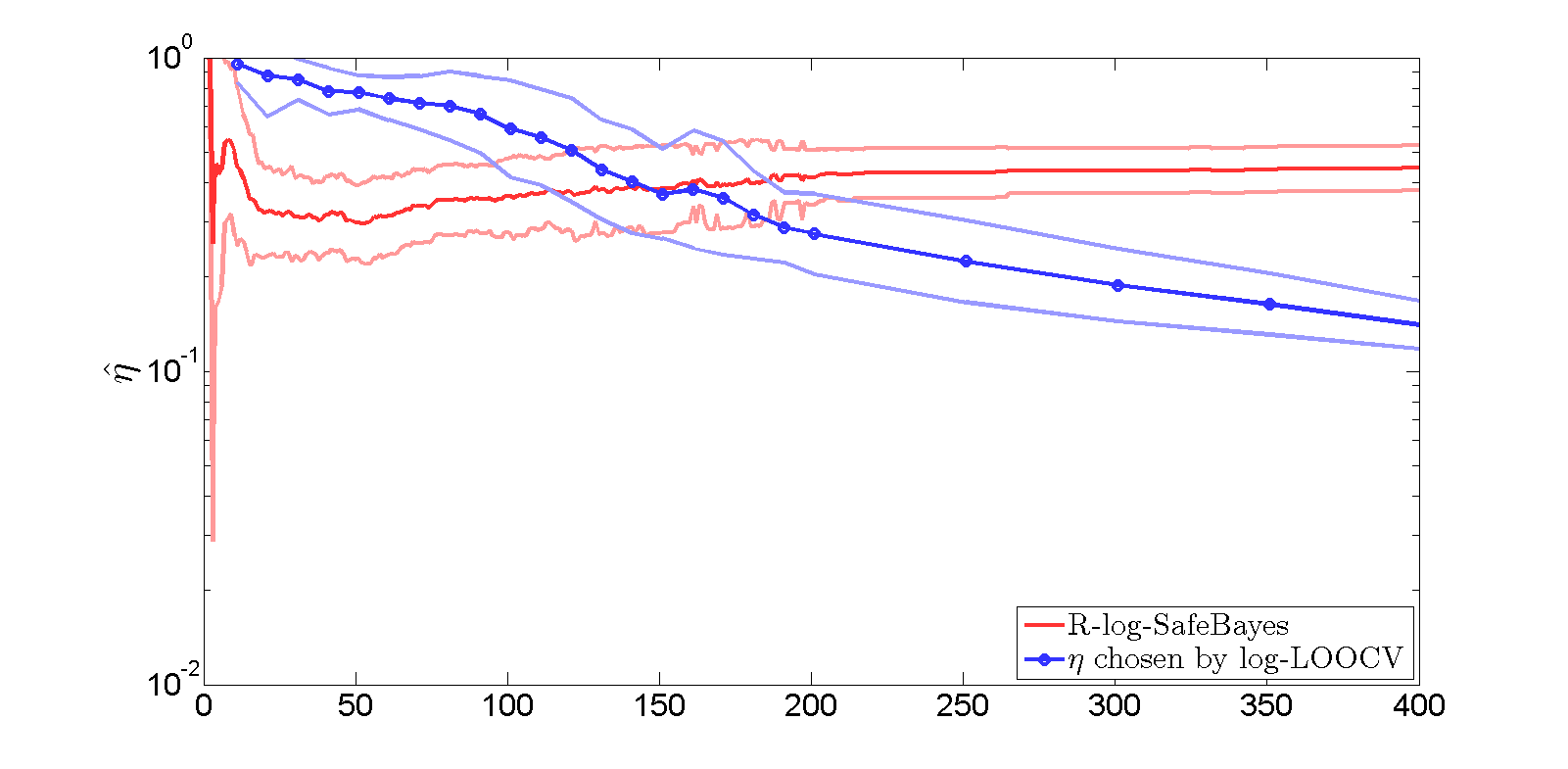}
\caption{\label{fig:canadb} Analogue of Figure~\ref{fig:mainexperimenta} for determining
  $\eta$ by leave-one-out cross-validation with log-loss.}}
\end{figure}

Perhaps more interestingly, in Figure~\ref{fig:canadb} we  show what
happens if we use LOO-cross validation based on the log-loss of the
Bayes predictive distribution, which may seem a reasonable procedure
from a `likelihoodist'  perspective. Here we see dismal behavior, the
reason being the hypercompression phenomenon of
Section~\ref{sec:hypercompression}: cross-validation will select a
model that, at the given sample size, has small log-risk, but because
of hypercompression this model can sometimes perform very badly in
terms of all the associated prediction tasks such as square-risk and
reliability. 
\pagebreak
\section{Experiments on Variations of the Truth}
\label{sec:wildvariations}
\paragraph{Other Distributions of Covariates}
In all experiments described in Section~\ref{sec:mainexperiments} and
the previous appendices, the
covariates $(X_{i1}, X_{i2}, \ldots )$ where sampled independently
from a $0$-mean multivariate Gaussian. We repeated most of our
experiments with $X_{i1}, X_{i2}, \ldots$ that were sampled
independently uniformly from $[-1,1]$, and, as already indicated in
the introduction, with polynomials, $X_{ij} =
S_i^j$ for $S_i \in [-1,1]$ uniform. This did not change the results in
any substantial way, so we do not report on it further.  

\paragraph{Fewer Easy and `Less-Easy'  Points}

If the fraction of `easy' points is reduced, one would expect the
performance of standard Bayes to improve. This is confirmed by an
experiment where each data point had a probability of only $1/4$ to be
$(0, 0)$. Here Bayes still  has some trouble finding the optimal
model, but the square-risk, MAP model order, and time taken to recover
are all much reduced compared to the original experiment in
Section~\ref{sec:modsel} where half the data points were `easy'.
SafeBayes on the other hand showed the same good performance as
before.

Two points that might be raised against the use of `easy' points in
our simulations are that they are unlikely to occur in practice, and
that if they were to occur, they would be easily detected and dealt
with another way. To address this line of argument to some extent,
another experiment was performed with a smaller contrast between
`easy' and `hard' points. Rather than being identically $(0, 0)$, the
`easy' points were random but with smaller variance than the `hard'
points. To be precise, the covariates and noise were both a factor 5
smaller (so that their variances were 25 times smaller). In this
experiment, the same phenomena as in Section~\ref{sec:modsel}
occurred, albeit again on a smaller scale (though larger than in the
previous, $1/4$-easy experiment).

\paragraph{Different Optimal Regression Functions}
We experimented with
a number of variations of the wrong-model experiment of
Section~\ref{sec:modsel}, by changing the underlying `true'
distribution $P^*$. In each variation, we still tossed, at each $i$,
an independent biased coin to determine whether $i$ would be `easy'
(still probability $1/2$) or `regular' (probability $1/2$), but in each
case we changed the definition of either the `easy' or the `regular'
instances or both. In all experiments, for the `regular' instances, only
$P^*(Y_i \mid X_i)$ was changed; the marginal distribution of the
$X_i$ was still multivariate normal as before.  Here is a list of
things we tried:
\begin{enumerate}
\item For regular instances, set $P^*(Y_i \mid X_i)$ so that $Y_i = 0 +
  \epsilon_i$ instead of (\ref{eq:pseudotruereg}), with $\epsilon_i$
  i.i.d.~normal as before; easy instances were still set to $({\bf
    0},0)$.
\item For regular instances, (\ref{eq:pseudotruereg}) was replaced by
  $Y_i = X_{i1} + X_{i2} + X_{i3} + X_{i4} + \epsilon_i$, so the
  optimal coefficients $\tilde{\beta}_1\ldots \tilde{\beta}_4$ are ten
  times as large as in the original experiment; easy instances were
  still set to $({\bf 0},0)$.
\item For regular instances, (\ref{eq:pseudotruereg}) was replaced by
  $Y_i = .1 \cdot (X_{i1} + \ldots + X_{i4}) - .04 + \epsilon_i$ (so
  the intercept is not $0$), and the easy instances were set to
  $(X_i,Y_i)$ = $({\bf .2}, .04)$, where ${\bf .2}$ represents the
  $K$-dimensional vector $(.2,\ldots, .2)$. Note that the easy points
  are on the optimal regression function.
\item For regular instances, (\ref{eq:pseudotruereg}) was replaced by
  $Y_i = .1 \cdot (X_{i1} + \ldots + X_{i4}) + .5 + \epsilon_i$ so the
  intercept was again not $0$; the easy instances were set to $({\bf
    0},.5)$.
\end{enumerate}
We explain each in turn.  For the first experiment, all the results
were comparable to the results of Experiment 1 in
Section~\ref{sec:mainexperiments}, so we do not list them.  For the
second experiment, the risks obtained by standard Bayes and SafeBayes
were  similar to each other. The model order
behaviors were similar to what they were before (with standard Bayes
selecting large model orders initially), but all methods recovered
much more quickly, converging on the optimal model shortly after
$n=50$; presumably this could happen because now the optimal
coefficients were substantially larger than the standard deviation in
the data.

The third experiment was included to see whether there would be an
effect if the `easy' points would be placed at an arbitrary point
rather than the special, fully symmetric $({\bf 0},0)$. We added the
intercept $- 0.04$ so as to make sure that, for the data we actually
observe, $\Exp_{X,Y \sim P^*}[Y_i] = (1/2).04 - (1//2).04 = 0$; thus the
$Y$-values will appear centered around $0$, which is standard both in
frequentist and Bayesian approaches to regression (for example, both
\cite{raftery1997bayesian} and \cite{HastieTF01} preprocess the data
so that $\sum_{i=1}^n Y_i = 0$). Again, we discerned no difference in
the results so did not include any further details.

Finally, the fourth experiment was included just to see what happens
if, contrary to standard methodology, we apply the method to $Y_i$
that are {\em not\/} (even approximately) centered. In this
experiment, standard Bayes did not converge to the optimal model until
after $n=150$ as in the experiment of Section~\ref{sec:modsel}, but
its risk and selected model orders were both smaller. The versions
of SafeBayes worked well as before.

\section{More on  Mix Loss}
\label{app:mixability}
\subsection{Implementing SafeBayes}
To implement the Safe Bayesian algorithm
(page~\pageref{alg:algorithm}), generalized posteriors must be
computed for different values of $\eta$, and the randomized loss
(\ref{eq:kibbeling}) must be computed for each sample size. For linear
models with conjugate priors as considered in our experiments, all
required quantities can be computed analytically. We have already seen
how to do this for models $\M_{\nc}$ with fixed dimension $\nc$. For
unions of such models, it turns out that the mix-loss is a helpful
tool. 

\paragraph{Role of mix loss in generalized posterior over models}

The generalized
posterior \emph{across} a discrete set of models is given by
(\ref{eq:genpostb}), which, writing $\tau = (\beta,\sigma^2)$, is,
via (\ref{eq:modelposteriorb}) and (\ref{eq:genpostd}), equivalent to 
\begin{align}\label{eq:vrijdag}
\pi(\nc \mid z^n, \eta)
&= \int_{\Theta_\nc} \pi(\nc,\tau \mid z^n,\eta) \,d\tau \nonumber \\
&\propto \int (\/ \dens(y^n \mid x^n, \tau, \nc) \/)^{\eta}
  \pi(\tau \mid \nc) \,d\tau \, \pi(\nc).
\end{align}
Here $\propto$ means `proportional to' when $p$ is varied and $z^n$
and $\eta$ are fixed.  In practice we prefer to calculate this quantity
incrementally: the posterior for $z^{n+1}$ with prior $\Pi$ is equal
to the posterior for a single data point $z_{n+1}$ when the posterior
for $z^n$ is used as prior (in this sense the generalized posterior
behaves like the standard posterior): using this to further rewrite the second
line of (\ref{eq:vrijdag}) gives
\begin{align*}
\pi(\nc \mid z^{n}, \eta)
&\propto \int (\/ \dens(y^n \mid x^n, \tau, \nc) \/)^{\eta}
  \pi(\tau \mid \nc) \,d\tau \, \pi(\nc) \\
&= \int (\/ \dens(y_n \mid x_n, \tau, \nc) \/)^{\eta}
 \cdot
(\/ \dens(y^{n-1} \mid x^{n-1}, \tau, \nc) \/)^{\eta}
   \pi(\tau \mid \nc) \,d\tau \, \pi(\nc) \\
& = 
\int (\/ \dens(y_n \mid x_n, \tau, \nc) \/)^{\eta}
 \cdot  \left( 
\pi(\tau \mid  z^{n-1}, \nc, \eta)  \cdot \int (\/ \dens(y^{n-1} \mid
x^{n-1},  \tau' )^{\eta}  \pi(\tau' \mid p)  d\tau' \,
\right) \,d\tau \, \pi(p)  \\
&\propto \int (\/ \dens(y_n \mid x_n, \tau, \nc) \/)^{\eta}
\cdot  \pi(\tau \mid  z^{n-1}, \nc, \eta) \,d\tau \cdot  \pi(p \mid z^{n-1},\eta),
\end{align*}
where in the third inequality we used the definition of the
generalized posterior and in the last we used (\ref{eq:vrijdag}).

The integral appearing in both the cumulative and the step-wise
expression equals the expectation in (\ref{eq:flattened}) from the
$\eta$-flattened $\eta$-generalized Bayesian predictive density for
$n$ and $1$ outcome respectively; $-\log [(\cdot)^{1/\eta}]$ of this
quantity is the mix loss of model $\nc$. We will now derive formulas
for this quantity.

\paragraph{Model with fixed variance}
Use the notation of Section~\ref{sec:instantiation}.
Write $\sigma^2_{\text{mix}} = \sigma^2 (1/\eta + x_{n+1} \Sigma_n
x_{n+1}^T)$. Then the mix loss for predicting one new data point
$y_{n+1}$ is
\begin{multline*}
-\log \pbayes(y_{n+1} \mid x_{n+1}, z^n, \langle\eta\rangle ; \eta)\\
= \frac{1}{\eta} \left[ \frac{1}{2}(\eta - 1) \log (2 \pi \sigma^2)
+ \frac{1}{2} \log \eta
+ \frac{1}{2} \log(2 \pi \sigma^2_{\text{mix}})
+ \frac{1}{2\sigma^2_{\text{mix}}}(y_{n+1} - x_{n+1} \beta_n)^2 \right]
\end{multline*}

\paragraph{Model with conjugate prior on variance}
Using the notation of Section~\ref{sec:instantiation}, 
the mix loss is given by
\begin{multline*}
-\log \pbayes(y_{n+1} \mid x_{n+1}, z^n, \langle\eta\rangle ; \eta)
= \frac{1}{\eta} \biggl[
\frac{1}{2} \eta \log \pi
+ \frac{1}{2} \log(1 + \eta x_{n+1} \Sigma_n x_{n+1}^T)\\
+ a_{n+1} \log(2 b_n +
  \frac{(y_{n+1} - x_{n+1} \beta_n)^2}{1/\eta + x_{n+1} \Sigma_n x_{n+1}^T})
- a_n \log 2 b_n
- \log \frac{\Gamma(a_{n+1})}{\Gamma(a_n)}
\biggr],
\end{multline*}

\subsection{Belief in Concentration (proof of Theorem~\ref{thm:concentration})}
\label{sec:concentration}
For simplicity, we only give the proof for the unconditional case, in
which the $\theta$ represent distributions $P_{\theta}$ on $z \in
\cZ$; extension to the conditional case is straightforward.
\newcommand{\badset}{\ensuremath{\bar{\Theta}_{\eta,\epsilon}}}
\renewcommand{\Pbayes}{\ensuremath{\Pi}} For $0 < \eta < 1$, let
$d_{\eta}(\theta^* \| \theta)$ denote the R\`enyi divergence of order
$1- \eta$ \citep{ErvenH14}, i.e.\ $d_{\eta}(\theta^* \| \theta) = -
\frac{1}{\eta} \log \Exp_{Z \sim \theta^*} \left(\frac{
    \dens_{\theta}(Z)}{\dens_{\theta^*}(Z)} \right)^\eta$.  We first
state a lemma, proved further below.  In the lemma, as in the
remainder of the proof, $(\theta^*,Z^n)$ is the random variable
distributed according to the Bayesian distribution $\Pi$. 

\begin{lemma}\label{lem:emma} Let $\Theta$, $\Pi$ and $\pi$ be as in the statement of Theorem~\ref{thm:concentration}. For every $1/2 \leq \eta < 1$, $\epsilon > 0$,
let $\badset := \{ \theta \in \Theta : d_{\eta}(\theta^* \| \theta) >
  \epsilon \}$. For every $b > 0$ and every sample size $n$ and
  setting $\epsilon :=  (b \log n)/( n\eta)$ and $c_{\eta} = (1- \eta)/(1 + \eta(1- \eta))$, we have:
$$
\Pbayes\left(\  \Pi(\badset \mid Z^n) \geq n^{- b c_{\eta}}\ \right) \leq 2 \left(\sum_{\theta \in
  \Theta} \pi(\theta)^{\eta}\right) \cdot n^{- b c_{\eta}}.
$$
In particular, if $\pi$ is summable for some $\eta < 1$, then using $b
= 1/c_{\eta}$, we get that the Bayesian probability that the posterior
probability of the set of $\theta$ farther than $b (\log n)/n$ from
$\theta^*$ exceeds $1/n$, is $O(1/n)$.
\end{lemma}
We proceed to prove Theorem~\ref{thm:concentration} using this lemma.
By the information inequality, we have for every probability density $\dens \neq
\dens_{\theta^*}$ that
$$
D(\theta^* \| \theta) = \Exp_{Z_n \sim P_{\theta^*}} [- \log \dens_{\theta}(Z_n) + \log \dens_{\theta^*}(Z_n) ] \geq \Exp_{Z_n \sim P_{\theta^*}} [- \log \dens_{\theta}(Z_n) + \log \dens(Z) ].
$$
In particular this holds with $\dens = \pbayes \mid Z^{n}$, the Bayes
predictive distribution based on the sample seen so far.
It then follows from (\ref{eq:mixgap}) that 
\begin{equation}\label{eq:rai}
\bar{\delta}_{n} \leq \Exp_{\theta \sim \Pi \mid Z^n} [D(\theta^* \| \theta)]
\end{equation}
Since $\pi^{\eta}$ is decreasing in $\eta$, we may
without loss of generality assume that the $\eta$ mentioned in the
theorem statement is at least $1/2$.  Now note \cite[Theorem
16]{ErvenH14} that for every $1/2 < \eta < 1$, $d_{1/2}(\theta^* \|
\theta) \leq (\eta/(1-\eta)) \cdot d_{\eta}(\theta^* \| \theta)$. We
also know from \cite[Lemma 4]{YangB99} that the KL divergence
$D(\theta^* \| \theta)$ satisfies $D(\theta^* \| \theta) \leq (2 +
\log v) d_{1/2}(\theta^* \| \theta)$.  
Since trivially  $d_{\eta}(\theta^* \| \theta) \leq \log v$, we have, with $C = \frac{\eta}{1- \eta} \cdot
(2+ 2 \log v)$, for every $\epsilon > 0$, using (\ref{eq:rai}),
\begin{multline}
\bar{\delta}_{n} \leq 
C \cdot \Exp_{\theta \sim \Pi \mid Z^n} [ d_{\eta}(\theta^* \|
\theta)]  \\ 
\leq C \Pi \left(d_{\eta} > \epsilon \mid Z^n \right) \log v + C
\left( 1- \Pi \left(d_{\eta} > \epsilon \mid Z^n \right) \right)\epsilon
\leq C \left( 
\Pi \left(d_{\eta} > \epsilon \mid Z^n \right) \log v + \epsilon\right), \nonumber
\end{multline}
so that $\Pi \left(d_{\eta} > \epsilon \mid Z^n \right) \geq (C^{-1}
\bar{\delta}_n - \epsilon)/(\log v)$ and by Lemma~\ref{lem:emma}, we have
for $\epsilon = b (\log n)/(n \eta)$ as in the lemma, that
$$
\Pbayes\left(\  
\frac{C^{-1} \bar{\delta}_n - \epsilon}{\log v}
\geq n^{- b c_{\eta}}\ \right) \leq 2 \left(\sum_{\theta \in
  \Theta} \pi(\theta)^{\eta}\right) \cdot n^{- b c_{\eta}}.
$$
Rewriting this expression, plugging in the value of $\epsilon$ and using $\eta \geq 1/2$, gives 
\begin{equation}\label{eq:schiphol}
\Pbayes\left(\  
\bar{\delta}_n
\geq C \left((\log v) n^{- b c_{\eta}} +  \frac{2 b (\log n)}{n}\right) \ \right)\leq 2 \left(\sum_{\theta \in
  \Theta} \pi(\theta)^{\eta}\right) \cdot n^{- b c_{\eta}}.
\end{equation}
The first part of the result follows by setting $b = a/c_{\eta}$. 
For the second result, note that the first result implies (take $a =
2$), by the union bound over sample sizes $1, \ldots, n$, that the
Bayesian probability that $\Exp_{Z^n \sim \theta^*} [\Delta_n]$
exceeds $C_0 \sum_{i=1}^n (\log i)/i \asymp (\log n)^2$ is $O(1/n)$. Thus there exists $C', C'_0$ such that the Bayesian probability that $\Exp_{Z^n \sim \theta^*} [\Delta_n]$
exceeds $C'_0 (\log n)^2$ is bounded by
$C'/n$. Thus for the probability in
(\ref{eq:heathrow}) we have
\begin{align*}
\Pi\left(\  
{\Delta}_n
\geq C_2  \cdot n^{a'} \ \right) =& 
\Pi\left(\  
{\Delta}_n
\geq C_2  \cdot n^{a'}, \Exp_{Z^n \sim \theta^*} [\Delta_n]
\geq C'_0 (\log n)^2   \ \right) + \\ &
\Pi\left(\  
{\Delta}_n  \geq C_2  \cdot n^{a'}, \Exp_{Z^n \sim \theta^*} [\Delta_n]
< C'_0 (\log n)^2   \ \right) \\ \leq &
\Pi\left(\  \Exp_{Z^n \sim \theta^*} [\Delta_n]
\geq C'_0 (\log n)^2   \ \right) + \\ &
\Pi\left(\  
{\Delta}_n  \geq C_2  \cdot n^{a'}, \Exp_{Z^n \sim \theta^*} [\Delta_n]
< C'_0 (\log n)^2   \ \right)
\\ \leq &  
\frac{C'}{n} + \frac{C'_0 (\log n)^2 
}{C_2 n^{a'}},
\end{align*}
where in the final step we used Markov's inequality. The second result follows.
\paragraph{Proof of Lemma~\ref{lem:emma}}
Fix $A > 0$ and $\gamma > 0$. We have
\begin{multline}\label{eq:verwarming}
\Pbayes\left(\Pi(\badset \mid Z^n) \geq  A
 \right) = 
\Pbayes\left(\frac{\sum_{\theta \in \badset} \pi(\theta) \cdot
    \dens_{\theta}(Z^n)}{
\sum_{\theta \in \Theta} \pi(\theta) \cdot
    \dens_{\theta}(Z^n)
} \geq  A
 \right)
= \\ 
\Pbayes\left(\frac{\sum_{\theta \in \badset} \pi(\theta) \cdot
    \dens_{\theta}(Z^n)}{\dens_{\theta^*}(Z^n)} \cdot \frac{\dens_{\theta^*}(Z^n)}
{\sum_{\theta \in \Theta} \pi(\theta) \cdot
    \dens_{\theta}(Z^n)} \geq  A
 \right) \leq  \\
\Pbayes\left( \frac{\sum_{\theta \in \badset} \pi(\theta) \cdot
    \dens_{\theta}(Z^n)}{\dens_{\theta^*}(Z^n)} \geq A^{1+ \gamma}
\right) + 
\Pbayes\left( 
\frac{\dens_{\theta^*}(Z^n)}
{\sum_{\theta \in \Theta} \pi(\theta) \cdot
    \dens_{\theta}(Z^n)}
\geq A^{-\gamma}
\right),
\end{multline}
where we used the union bound.
The first term is equal to, and can be further bounded as
\begin{align*}
= & \Pbayes\left( \frac{\left(\sum_{\theta \in \badset} \pi(\theta) \cdot
    \dens_{\theta}(Z^n)\right)^{\eta}}{(\dens_{\theta^*}(Z^n))^{\eta}} \geq A^{\eta(1+ \gamma)}
\right) 
\leq   \Pbayes\left( \frac{\sum_{\theta \in \badset} \pi(\theta)^{\eta} \cdot
    (\dens_{\theta}(Z^n))^{\eta}}{(\dens_{\theta^*}(Z^n))^{\eta}} \geq A^{\eta(1+ \gamma)}
\right) \\
= & \sum_{\theta^*} \pi(\theta^*) P_{\theta^*} \left( \frac{\sum_{\theta \in \badset} \pi(\theta)^{\eta} \cdot
    (\dens_{\theta}(Z^n))^{\eta}}{(\dens_{\theta^*}(Z^n))^{\eta}} \geq A^{\eta(1+ \gamma)}
\right) \\ \leq & \sum_{\theta^* \in \Theta} \pi(\theta^*) \Exp_{Z^n
  \sim P_{\theta^*}} \left[ \frac{\sum_{\theta \in \badset} \pi(\theta)^{\eta} \cdot
    (\dens_{\theta}(Z^n))^{\eta}}{(\dens_{\theta^*}(Z^n))^{\eta}}
\right] \cdot A^{- \eta(1+ \gamma)} \\ = & 
\sum_{\theta^* \in \Theta} \pi(\theta^*) \sum_{\theta \in
    \badset} \pi(\theta)^{\eta} \cdot \left( \Exp_{Z
  \sim P_{\theta^*}} \left[ \frac{ 
    (\dens_{\theta}(Z))^{\eta}}{(\dens_{\theta^*}(Z))^{\eta}}\right] \right)^n
\cdot A^{- \eta(1+ \gamma)} \\ \leq & \left( \sum_{\theta \in
    \badset} \pi(\theta)^{\eta}\right) e^{-n \eta \epsilon} \cdot A^{- \eta(1+ \gamma)}.
\end{align*}
where the first inequality follows by differentiation to $\eta$ (or
equivalently, by monotonicity of $\ell^p$-norms), the second is
Markov's, and the third is the definition of R\`enyi divergence.

The second term in (\ref{eq:verwarming}) can be bounded as
\begin{align*}
\leq & 
\Pbayes\left( 
\frac{\dens_{\theta^*}(Z^n)}
{\pi(\theta^*) \cdot
    \dens_{\theta^*}(Z^n)}
\geq A^{-\gamma}
\right) = \Pbayes(\pi(\theta^*)^{-1+ \eta} \geq A^{-(1-\eta) \gamma} )
\leq  \Exp_{\theta^* \sim \Pbayes} [ \pi(\theta^*)^{-1+ \eta}] A^{
  \gamma(1- \eta)} \\ & = 
\sum_{\theta^*} \pi(\theta^*)^{\eta} A^{
  \gamma(1- \eta)}.
\end{align*}
Combining the upper bounds on the two terms on the right in
(\ref{eq:verwarming}), we get: 
$$
\Pbayes\left(\Pi(\badset \mid Z^n) \geq  A
 \right) \leq 
\left( \sum_{\theta \in
    \badset} \pi(\theta)^{\eta}\right) \left( 
e^{-n \eta \epsilon} \cdot A^{- \eta(1+ \gamma)} + A^{
  \gamma(1- \eta)} \right). 
$$
Now we plug in the chosen value of $\epsilon = (b \log n)/(n \eta)$ and we
set $A = n^{- b/(\gamma + \eta)}$. With these values the second factor
on the right becomes
$$
e^{-n \eta \epsilon} \cdot A^{- \eta(1+ \gamma)} + A^{
  \gamma(1- \eta)} = n^{-b} n^{b(\eta(1 + \gamma))/(\gamma + \eta)}+
n^{-b \gamma(1-\eta)/(\gamma+ \eta)} = 2 n^{- b \cdot \gamma \cdot
  \frac{1- \eta}{\gamma + \eta}}.
$$
Since this holds for al  $\gamma >0$, it also holds for $\gamma = 1/(1-\eta)$,
and the result follows. 
\end{document}

%% file: declarations.tex
\usepackage{epsf}
\usepackage{verbatim}
\usepackage{makeidx}
\usepackage{graphicx}
\usepackage{amssymb}
\usepackage{amsmath}
\usepackage{latexsym}
\usepackage{amsbsy}
\usepackage{epsfig}

\newtheorem{theorem}{Theorem}

\newtheorem{lemma}{Lemma}

\newcommand{\commentout}[1]{}
\newcommand{\prs}{\par}

\newcounter{myenumeratecounter}

\newtheorem{examplehidden}{Example}
\newenvironment{example}{
\begin{examplehidden}
\em
}
{\end{examplehidden}}

\newenvironment{ownquote}
               {
\list{}{\rightmargin\leftmargin}%
                \small \item[]}
               {\endlist}

\newcommand{\Pbayes}{\ensuremath{P}_{\text{B}}}

\newcommand{\cA}{\ensuremath{\mathcal A}}

\newcommand{\cM}{\ensuremath{\mathcal M}}

\newcommand{\cS}{\ensuremath{\mathcal S}}

\newcommand{\cX}{\ensuremath{\mathcal X}}
\newcommand{\cY}{\ensuremath{\mathcal Y}}
\newcommand{\cZ}{\ensuremath{\mathcal Z}}

\newcommand{\reals}{{\mathbb R}}

\newcommand{\M}{\ensuremath {{\cal M}}}

\newlength{\ownboxwidth} \newlength{\ownleftboxwidth}
  \newlength{\ownrightboxwidth} 
\newcommand{\longp}[1]{}
\newcommand{\shortp}[1]{}

\newcommand{\pbayes}{\ensuremath{p_{\text{\rm Bayes}}}}
\newcommand\Exp{\ensuremath{\mathbf E}}